\documentclass[11pt,a4paper,leqno]{article}

\usepackage{amsmath,amscd,amsfonts,amsthm,amssymb,dsfont,bbm,manfnt}
\usepackage{url,color,appendix,eucal,comment,tabu,float}
\usepackage[usenames,dvipsnames]{xcolor}
\usepackage[nottoc,numbib]{tocbibind}
\usepackage{multicol}
\usepackage[OT2,T1]{fontenc}
\usepackage[utf8]{inputenc}
\usepackage[all]{xy}
\usepackage{lscape,rotating} 
\usepackage{graphicx}
\usepackage[nolabel]{showlabels}
\usepackage[hypertexnames=false]{hyperref} 

\newenvironment{pf}
{\medskip\noindent {\it Proof.  }}
{\hfill\nobreak $\Box$ \par\bigbreak}

\newcommand{\isomo}{\overset{\sim}{\rightarrow}}

\newcommand{\ps}{\par \smallskip}
  
\newcommand{\N}{\mathbb{N}}
\newcommand{\Z}{\mathbb{Z}}

\newcommand{\Q}{\mathbb{Q}}

\newcommand{\R}{\mathbb{R}}    
\renewcommand{\C}{\mathbb{C}}

\newcommand{\AAA}{\mathbb{A}}

\newcommand{\sff}{\sffamily\selectfont}

\newcounter{introcounter}
\newtheorem{thm}[subsection]{Theorem}
\newtheorem{lemma}[subsection]{Lemma}
\newtheorem{propdef}[subsection]{Proposition-Definition}
\newtheorem{remark}[subsection]{Remark}
\newtheorem{cor}[subsection]{Corollary}

\newtheorem{prop}[subsection]{Proposition}
\newtheorem{example}[subsection]{Example}

\newtheorem{thmintro}{Theorem}
\newtheorem{defintro}[introcounter]{Definition}
\newtheorem{exampleintro}[introcounter]{Example}
\newtheorem{corintro}[thmintro]{Corollary}
\newtheorem{remarkintro}[introcounter]{Remark}

\newtheorem{definition}[subsection]{Definition}

\newtheorem{fact}[subsection]{Fact}

\usepackage{titlesec}

\titleformat{\section}{\bf \normalsize}{\arabic{section}.}{1 em}{}

\titleformat{\subsection}
{}% apparence commune au titre et au numXro
{\bf \arabic{section}.\arabic{subsection}.}% apparence du numXro
{0.5 em}% espacement numXro/texte
{\normalsize}% apparence du titre

\titleformat{\subsubsection}[runin]
{\small \bf}% apparence commune au titre et au numXro
{}% apparence du numXro
{}% espacement numXro/texte
{}% apparence du titre []
\hypersetup{colorlinks=true,linkcolor=black,citecolor=black, urlcolor=cyan,pdfborderstyle={/S/U/W 1}}

\makeatletter

\@addtoreset{equation}{section}
\makeatother

\makeatletter

\@addtoreset{table}{section}
\makeatother

\title{Unimodular hunting}

\author{Ga\"etan Chenevier\thanks{\'Ecole Normale Sup\'erieure-PSL, D\'epartement de Math\'ematiques et Application, 45 rue d'Ulm, 75230 Paris Cedex, France. During this work, the author has been supported by the C.N.R.S. and by the project ANR-19-CE40-0015-02 COLOSS.}}

\begin{document}
\maketitle
%\ps\ps
%\begin{center}{With an appendix by Bill Allombert and Ga\"etan Chenevier}\end{center}
%\ps\ps

\begin{abstract} 
We develop a method initiated by Bacher and Venkov, and based on a study of the Kneser neighbors of the standard lattice $\Z^n$, which allows to classify the integral unimodular Euclidean lattices of rank $n$. As an application, of computational flavour, we determine the isometry classes of unimodular lattices of rank $26$ and $27$. 
\end{abstract}

\section{Introduction}

\subsection{The classification of unimodular lattices} 
\label{subsectintro:classuni}
Consider the standard Euclidean space $\R^n$, with inner product $x.y\,=\,\sum_i \,x_i y_i$. 
Recall that a lattice $L \subset \R^n$ is called {\it integral} if we have $x.y \in \Z$ for all $x,y \in L$, 
and {\it unimodular} if its covolume is $1$ (see~\S\ref{notations} for the basics on integral lattices).  
We denote by $\mathcal{L}_n$ the set of all integral unimodular lattices in $\R^n$ and by 
${\rm X}_n$ the (finite) set of isometry classes of such lattices. 
The most trivial element of $\mathcal{L}_n$, 
that will nevertheless play a major role here, 
is the {\it standard} or {\it square} lattice: 
\begin{equation}\label{defin}  {\rm I}_n:=\Z^n. \end{equation}
For $n\leq 7$ the isometry class of ${\rm I}_n$ is the unique element of ${\rm X}_n$, a well-known fact with famous contributions from Lagrange, Gauss, Hermite and Minkowski. 
For $n=8$ there is also the ${\rm E}_8$ lattice, which is the unique other isometry class in ${\rm X}_8$ by Mordell,
and the first example of an {\it even} unimodular lattice.
Thanks to the works of many authors, including Witt, Kneser, Niemeier, Conway $\&$ Sloane and Borcherds, see {\it e.g.} \cite{kneser16, Ni, conwaysloane, borcherds}, representatives of ${\rm X}_n$ have then been determined before this work up to $n = 25$.

\begin{table}[H]

\tabcolsep=5pt
{\scriptsize \renewcommand{\arraystretch}{1.8} \medskip
\begin{center}
\begin{tabular}{c|c|c|c|c|c|c|c|c|c|c|c|c|c|c|c}
 $n$ & $1-7$ & $8-11$ & $12-13$ & $14$ & $15$ & $16$ & $17$ & $18$  & $19$ & $20$ & $21$ & $22$ & $23$ & $24$ & $25$  \\
\hline
$|{\rm X}_n|$ & $1$ & $2$ & $3$ & $4$ & $5$ & $8$ & $9$ & $13$  & $16$ & $28$ & $40$ & $68$ & $117$ & $297$ & $665$ 
\end{tabular} 
\end{center}
} 
\caption{The size of ${\rm X}_n$ for $n\leq 25$.}
\label{tab:cardcarn} 
\end{table}
%OEIS: 

The subset ${\rm X}_n^{\emptyset} \subset {\rm X}_n$ of classes of lattices 
with no vector of norm\footnote{Following a standard abuse of language, the {\it norm} of an element $v$ is defined as $v.v$.} $\leq 2$ is especially interesting.
By the aforementioned classifications, we know $|{\rm X}_n^{\emptyset}|=0$ for $n<23$ or $n=25$, $|{\rm X}_{23}^{\emptyset}|=1$ (Short Leech lattice), $|{\rm X}_{24}^{\emptyset}|=2$ (Leech and Odd Leech lattices).
Moreover, Borcherds showed
$|{\rm X}_{26}^{\emptyset}|=1$ in \cite[Thm. 3.6]{borcherds}, and Bacher and Venkov proved $|{\rm X}_{27}^{\emptyset}|=3$ and $|{\rm X}_{28}^\emptyset|=38$ in \cite{bachervenkov}. 
Our first main result in this paper is the following:

\begin{thmintro} \label{X2627}
We have $|{\rm X}_{26}|=2566$ and $|{\rm X}_{27}|=17059$.
%Moreover, representatives for the isometry classes in ${\rm X}_{26}$ and ${\rm X}_{27}$ 
%are listed in {\rm \cite{cheweb}}.
\end{thmintro}

Let us mention that general lower bounds on $|{\rm X}_n|$ may be obtained using the Minkowski-Siegel-Smith mass formula, although they are very bad for ``small'' $n$ as in our range. Much better lower bounds have been obtained by O. King in \cite{king} in the case $n\leq 32$. He proved in particular $|{\rm X}_{26}| \geq 2307$ and $|{\rm X}_{27}|\geq 14179$. Our computations show that King's estimate were not too far from the actual values.
As we shall see later, King's computations also played an important role in our search.
See Sect.~\ref{sect:newlatt} for a brief analysis of the isometry groups of the lattices of Theorem~\ref{X2627}.
\ps

\subsection{The cyclic $d$-neighbors of ${\rm I}_n$}
\label{subsecintro:cyclicdnei}

Our goal in proving Theorem~\ref{X2627} is actually 
not only to classify all the aforementioned lattices, 
but also to provide constructions of all of them as {\it cyclic neighbors} of 
the simplest lattice of all, namely of ${\rm I}_n$. 
The cyclic $d$-neighbors of a unimodular lattice $L \in \mathcal{L}_n$ 
are the unimodular lattices $N \in \mathcal{L}_n$ with $L/(N\cap L) \simeq \Z/d$. 
This is a fairly classical variant of Kneser's original definition of neighbor lattices~\cite{kneser16,bacher,bachervenkov,scharlau,chlannes,voight}, 
see \S \ref{sectkneser} for some background on this notion.\footnote{For our purposes, it will be important to allow $d$ to be an arbitrary integer. 
Many references only treat in details the case where $d$ is odd, or assume the lattices to be even for $d$ even.} 
In the sequel, we usually omit the adjective {\it cyclic} and just talk about $d$-neighbors. We now recall the concrete construction of the $d$-neighbors of ${\rm I}_n$. \ps
Fix $d\geq 1$ an integer and $x=(x_i)$ in ${\rm I}_n$ with ${\rm gcd}(d,x_1,\dots,x_n)=1$. 
Then the image of $x$ in ${\rm I}_n \otimes \Z/d$ generates a {\it line}, 
{\it i.e.} a cyclic subgroup $l:=l(x)$ of order $d$. 
The orthogonal of this subgroup in ${\rm I}_n$ is the lattice 
\begin{equation} 
\label{defmx} {\rm M}_d(x) := {\rm M}_d(l) := \{ v \in {\rm I}_n \, \, |\, \, \sum_{i=1}^n x_i v_i \equiv 0 \bmod d \}
\end{equation}
and it satisfies ${\rm I}_n/{\rm M}_d(x) \simeq \Z/d$. 
Any subgroup $M \subset {\rm I}_n$ with ${\rm I}_n/M \simeq \Z/d$ 
has the form $M={\rm M}_d(l)$ for a unique line $l \subset {\rm I}_n \otimes \Z/d$. 
Set $e=1$ if $d$ is odd and $e=2$ otherwise.
We say that $x$ (or $l$) is {\it $d$-isotropic} if we have 
\begin{equation} \label{defdiso}\sum_{i=1}^n x_i^2 \equiv 0 \bmod ed. \end{equation}
It is a fact that there is $N \in \mathcal{L}_n$ with $N \cap {\rm I}_n = {\rm M}_d(x)$ if, 
and only if, $x$ is $d$-isotropic, and if so, there are exactly $e$ such $N$. 
The following formula, in which we choose $x'  \in {\rm I}_n$ 
arbitrarily with $x' \equiv x \bmod d$ and $x'.x' \equiv 0 \bmod d^2$, 
defines the $1$ or $2$ possible unimodular lattices $N$ with $N\cap {\rm I}_n = {\rm M}_d(x)$:
\begin{equation}\label{defldx} {\rm N}_d(x')  = {\rm M}_d(x) + \Z \frac{x'}{d}. \end{equation}
It is easily checked that Formula \eqref{defldx} indeed defines an integral unimodular lattice.
For $d$ odd, the lattice ${\rm N}_d(x')$ does not depend on the choice of $x'$ as above, 
and we simply denote it by ${\rm N}_d(x)$. For $d$ even, 
we temporarily denote by ${\rm N}_d(x)^{\pm }$ the two possibilities for ${\rm N}_d(x')$, 
postponing to \S \ref{subsect:dneiin} the discussion of how to distinguish them. 
The lattices ${\rm N}_d(x)^+$ and ${\rm N}_d(x)^-$ are not isometric in general, 
but they are isometric if we have $x_i \equiv d/2 \bmod d$ for some $i$. \ps

\subsection{Theoretical exhaustion}
\label{subsectintro:theoexaust}
Before giving examples, we mention a key fact, 
proved by Hsia and J\"ochner in \cite[Cor. 4.1]{HJ}, asserting that given any (say) odd $L,L' \in \mathcal{L}_n$, there are infinitely many primes $p$ such that $L'$ is isometric to a $p$-neighbor of $L$.
In the companion paper \cite{chestat}, we proved several quantitative variants of this result (by very different methods). We give here yet another variant of these results.
Assume $d$ odd to simplify. Then any odd element in $\mathcal{L}_n$ has the same (explicit) number ${\rm c}_n(d)$ of $d$-neighbors, and for $n>2$ we have ${\rm c}_n(d) \sim d^{n-2}$ for $d \rightarrow +\infty$. If $L$ is an Euclidean lattice we denote by ${\rm O}(L)$ 
its (finite) isometry group and define its {\it mass} by ${\rm mass}(L) = \frac{1}{|{\rm O}(L)|}$. 
Also, we denote by ${\rm m}_n^{\rm odd}$ the mass of the genus of odd unimodular lattices of rank $n$ (see \S~\ref{subsect:minkosiegsmith}).

\begin{thmintro} 
\label{thmi:stat1}
Let $L,L' \in \mathcal{L}_n$ be any odd unimodular lattices of rank $n$. For an integer $d\geq 1$, denote by ${\rm n}_d(L,L')$ the number of $d$-neighbors of $L$ which 
are isometric to $L'$. Then we have 
$$\frac{{\rm n}_d(L,L')}{{\rm c}_n(d)} \rightarrow \frac{{\rm mass}(L')}{{\rm m}_n^{\rm odd}}\, \, \, \, {\rm for}\, \, d\, \, {\rm odd}\, \,{\rm and}\,\, d \,\,\rightarrow \, \, +\infty.$$
\end{thmintro} 
\ps

In particular, {\it for odd $d \rightarrow \infty$, any odd $L' \in \mathcal{L}_n$ 
appears as a $d$-neighbor of ${\rm I}_n$ with a probability proportional to its mass ${\rm mass}(L')$}.
Theorem~\ref{thmi:stat1} follows from~\cite[Thm. A]{chestat} if one assumes $d$ prime in its statement : 
see Sect.~\ref{sect:pfthmb} for the general case.
Let us stress that this result, however, does not say anything about the smallest integer $d$ such that 
a given $L \in \mathcal{L}_n$ is isometric to a $d$-neighbor of ${\rm I}_n$. 
This quantity is one of the most difficult to predict, and deserves a definition.

\begin{defintro} \label{def:farness} For $L \in \mathcal{L}_n$, the {\rm farness} of $L$,  denoted by ${\rm far}(L)$, is the smallest integer $d\geq1$ such that $L$ is isometric to a (cylic) $d$-neighbor of ${\rm I}_n$. 
\end{defintro}

\subsection{First examples}
\label{subsect:firstexamples}

We now give a few interesting examples. 
The element $1^n=(1,1,\dots,1)$ of ${\rm I}_n$ is $2$-isotropic for $n \equiv 0 \bmod 4$,
so we have a unimodular lattice ${\rm N}_2(1^n)^{\pm}$ for such an $n$.
We must have ${\rm N}_2(1^4)^{\pm} \simeq {\rm I}_4$ 
but recognize ${\rm N}_2(1^8)^{\pm} \simeq {\rm E}_8$. 
Also, ${\rm N}_2(1^{12})^{\pm}$ is the unique rank $12$ unimodular lattice 
with no vector of norm $1$, and for $n \equiv 0 \bmod 8$ and $n>8$ 
the lattice ${\rm N}_2(1^n)^{\pm }$ is even with root system ${\rm D}_n$ 
and is sometimes denoted by ${\rm E}_n$ or ${\rm D}_n^+$ in the literature. 
The {\it Leech lattice} also has the following beautiful description 
due to Thompson (see \cite{conwaysloane} p. lvi) as a $94$-neighbor of ${\rm I}_{24}$:
\begin{equation}
\label{eq:LeechThompson} 
{\rm Leech} \simeq {\rm N}_{94}(x)^{\pm}, \, \, x=(1,3,5,7,...,47) \in \Z^{24}.
\end{equation}
In the same vein, 
we have $1^2+2^2+\cdots+n^2=\frac{1}{6}n(n+1)(2n+1) \equiv 0 \bmod 2n+1$ 
for $n \not \equiv 1 \bmod 3$, hence a unimodular lattice 
\begin{equation}
\label{eq:oddshortborchleech} 
{\rm N}_{2n+1}(1,2,3,\dots,n) \in \mathcal{L}_n,\, \, \, n \not \equiv 1 \bmod 3.
\end{equation}
It may be shown that this lattice has no vector of norm $\leq 2$ for $n\geq 23$ 
(see Sect.~\ref{sec:ex12n} for a study of those lattices). For $n=23$, $24$ and $26$ 
we recover this way respectively the short Leech lattice, 
the odd Leech lattice and Borcherds' lattice in ${\rm X}_n^\emptyset$ !  
It is hard to think of a simpler definition for these lattices than those ones.
All of our lattices will be given in this form:

\begin{thmintro} \label{thm:listuninei}
A \href{http://gaetan.chenevier.perso.math.cnrs.fr/unimodular_lattices/unimodular_lattices.gp}{list} of $(d,x)$ such that the ${\rm N}_d(x')$ are representatives for all the unimodular lattices of Theorem~\ref{X2627} is given in~{\rm \cite{cheweb}}.
\end{thmintro}

This extends previous work by Bacher~\cite{bacher} in the case $n\leq 24$, as well as the study of ${\rm X}_n^\emptyset$ for $n\leq 28$ in \cite{bachervenkov}.
This (partially aesthetic) wish of giving all of our lattices 
in the form ${\rm N}_d(x')$ added in practice a number of extra difficulties, 
and forced us to find neighbor constructions of some lattices 
more naturally defined in other ways (see e.g. Sect.~\ref{sec:excsec} \&~\ref{sect:lattconstr}). 

\begin{remarkintro} {\rm
An further question is to determine, for each lattice $L={\rm N}_d(x')$ given in our list,
the farness of $L$. We obviously have ${\rm far}\, L \leq d$, and a neighbor form $(d,x')$ will be called {\it optimal} if we have ${\rm far}({\rm N}_d(x'))=d$.
As an example, it is easy to see that we have ${\rm far}({\rm Leech}) \geq 94$ (see Example \ref{ex:emptyrs}), so that Thompson's construction~\eqref{eq:LeechThompson} is optimal. Many neighbor forms in our lists are actually optimal, but certainly not all, and we leave as an open question to determine optimal forms for all of them.}
% (this question is not as futile as it may seem!).}
\end{remarkintro}
\ps

\subsection{The general method}
\label{subsect:generalalgo}
	
	\indent 	The basic strategy we follow to prove Theorem~\ref{X2627}, which of course heavily relies on computer calculations,\footnote{In all this work, we heavily used the open-source computer algebra system \texttt{PARI/GP}.} 
is well-known:  if we are able to produce non-isometric lattices $L_1,\dots,L_h$ in $\mathcal{L}_n$ such that $\sum_{i=1}^h {\rm mass}(L_i)$ coincides with the total mass of $\mathcal{L}_n$, whose exact value is known thanks to the mass formulae (see e.g. \cite[Chap. 16]{conwaysloane}  p. 409), then $L_1,\dots,L_h$ must be representatives of ${\rm X}_n$. 
This strategy requires at least three ingredients:\ps

{\bf (M1) Mass computations}. 
We obviously use the Plesken-Souvignier algorithm \cite{pleskensouvignier}, 
which efficiently computes $|{\rm O}(L)|$ from a given Gram matrix of $L$ with small diagonal. 
This algorithm turned out to work well for unimodular lattices of rank $\leq 27$, 
which tend to be generated by vectors of norm $\leq 3$, especially when combined with root system arguments introduced in~\cite{chniemeier} : 
see Remark~\ref{pleskensouvignier} for a discussion about this point.\ps

{\bf (M2) Finding invariants}. Although the Plesken-Souvignier algorithm also allows to check 
whether two lattices are isometric or not,
it is unrealistic to rely on such an isometry test in situations 
like ours where millions of lattices will have to be compared. 
Instead, {\it we define a few ad hoc (and easy enough to compute) invariants, 
and we bet in our search that they be enough 
to distinguish the elements in ${\rm X}_n$ for our specific $n$}. 
The most obvious invariants of a lattice $L$ are its {\it configuration of vectors of norm $i$} or {\it $\leq i$}
\begin{equation}
\label{eq:fms}
{\rm R}_i(L) = \{ v \in L\, |\, \, v.v=i\}\, \,\,\, {\rm and}\, \, \,\,{\rm R}_{\leq i}(L) = \{ v \in L\, |\, \, v.v \leq i\},
\end{equation}
say viewed as a finite (Euclidean) metric spaces (see Sect.~\ref{sec:euclideansets}). 
For later use we also set ${\rm r}_i(L):=|{\rm R}_i(L)|$.
As is well-known, for any integral lattice $L$  there is a unique decomposition
\begin{equation}\label{decnorm1} L = A \perp B \,\,\,{\rm with}\,\, \, \, A \simeq {\rm I}_m\, \, {\rm and}\, \, {\rm r}_1(B) = 0\end{equation}
(and we have $2m={\rm r}_1(L)$). So we may and do restrict to classify the lattices $L \in \mathcal{L}_n$ with ${\rm r}_1(L) = 0$. 
The Euclidean set ${\rm R}_2(L)$ is called the {\it root system} of $L$, and is a disjoint union of classical
${\rm ADE}$ root systems (see~\S\ref{subsect:rootsystems}). It is very easy to determine in practice (see Remark~\ref{rem:rscalc}). \ps

A well-known but curious fact about unimodular lattices of rank $\leq 23$, and about Niemeier lattices, 
is that they are uniquely determined by their ${\rm R}_{\leq 2}$. 
This is however only obtained as a by-product of the classification in these cases.
It does not hold anymore for ${\rm X}_{24}$, which contains for instance two odd lattices with empty ${\rm R}_1$
and same root system $8 {\bf A}_1 \,4 {\bf A}_3$. 
It is natural to study ${\rm R}_{\leq 3}$ to go further. 
It turns out that most lattices in our range are indeed spanned (over $\Z$) by vectors of norm $3$.
The main difficulty is that contrary to the case 
${\rm R}_{\leq 2}$ we are not aware of any existing study or classification of configurations of vectors of norm $\leq 3$. 
We refer to~\S \ref{subsect:norm3}, as well as the companion paper \cite{allche},
for a few invariants of ${\rm R}_{\leq 3}$ that proved useful in our situations.  As a by-product of Theorem~\ref{X2627}, we obtain the following result.

\begin{corintro} 
\label{cori:invnorm3}
Let $L, L' \in \mathcal{L}_n$ with $n\leq 27$. Then $L$ and $L'$ are isometric if, and only if, 
the Euclidean sets ${\rm R}_{\leq 3}(L)$ and ${\rm R}_{\leq 3}(L')$ are isometric.
\end{corintro}
\ps
{\bf (M3) Neighbors enumeration}. 
The basic idea is to enumerate, with the computer, and for increasing integers $d=2,3,4,\dots,$ all the $d$-isotropic 
vectors $x$ in ${\rm I}_n$ and study their associated lattices ${\rm N}_d(x')$, by computing their invariant and, if they are new, their mass, and so on until the total mass of $\mathcal{L}_n$ is exhausted. 
In order to exclude
trivial isometries between neighbors induced by ${\rm O}({\rm I}_n)$ 
(permutations and sign changes of the coordinates of $x$), we may restrict the enumeration to $d$-isotropic vectors $x\in \Z^n$
satisfying\footnote{This step of dealing with the orbits, immediate here, gives a clear advantage to this method compared to Kneser's original one, 
consisting of computing successive $2$-neighbors until exhaustion. Another important advantage is that ${\rm I}_n$
has a both a large root system and a large automorphism group, which will allow to efficiently biase our search.}
\begin{equation} \label{mainineq} 0 \leq x_1 \leq x_2 \leq \cdots \leq x_d \leq d/2.\end{equation}
As two $d$-isotropic vectors generating the same $\Z/d$-line give rise to the same neighbor, we also usually further assume $x_1=1$ and that $x_1$ has the maximal multiplicity among the $x_i$'s 
(this is not a restriction at all for $d$ prime). \ps

\subsection{The splitting root system by root system}
\label{subsect:splittingrs}
${}^{}$\indent
Most of the computation time in the algorithm (M3) above is spent computing our invariants of the lattices ${\rm N}_d(x')$ found in the enumeration. These invariants always include the much faster computation of the root system of ${\rm N}_d(x')$. 
Also, as will be clear later when discussing the visible root system, the algorithms starts by finding the lattices with the biggest root systems (and, often, quickly finds all of them). It is thus highly desirable to split our search root system by root system and not to compute the full invariant when all the lattices with a given root system have already been found. \par
This is fortunately permitted by the aforementioned work \cite{king}. Indeed, as explained {\it loc. cit.}, although the main computation there concerns even unimodular lattices of rank $32$, it allows to determine, for any $n\leq 30$ and any root system of rank $n$, the mass of the subset of $\mathcal{L}_n$ consisting of lattices with root system isomorphic to $R$. The details of this step are not fully given in \cite{king}: we give another point of view on it in \S~\ref{sect:mass}, and explain as well how to deduce similar results for several other genera of interest. The table below compares the number of possible root systems, given by King, to the actual number of lattices in known cases (Table~\ref{tab:cardcarn} and Theorem~\ref{X2627}), in dimension $n \leq 30$:

\tabcolsep=2.5pt

\begin{table}[H]
{\scriptsize 
\renewcommand{\arraystretch}{1.8} \medskip
\begin{center}
\begin{tabular}{c|c|c|c|c|c|c|c|c|c|c|c|c|c|c|c}
$n$ & $12$ & $14\!\!-\!\!17$ & $18$ & $19$ & $20$ & $21$ & $22$ & $23$ & $24$ & $25$ & $26$ & $27$ & $28$ & $29$  & $30$ \\
\hline
$\texttt{rs}$ & $1$ & $1$ & $4$ & $3$ & $12$ & $12$ & $28$ & $49$ & $149$ & $327$ & $1086$ & $2797$ & $4722$ & $11\,085$  & $18\,220$ \\

$\texttt{min}$ & $1$ & $1$ &  $4$ & $3$ & $12$ & $12$ & $28$ & $49$ & $156$ & $360$ & $1626$ & $11671$ & $312\,287$ & $37\,604\,456$  & $20\,131\,670\,647$\\

${\bf \sharp}$ & $1$ & $1$ & $4$ & $3$ & $12$ & $12$ & $28$ & $49$ & $156$ & $368$ & $1901$ & $14493$ & $?$ & $?$ & $?$
\end{tabular} 
\end{center}
}
\caption{\small{Number $\texttt{rs}$ of isometry classes of root systems of unimodular lattices in $\mathcal{L}_n$ with no norm $1$ vectors \cite{king}, lower bound $\texttt{min}$ for the number of such classes by the method in \cite{king}, compared to the actual number ${\bf \sharp}$ of isometry classes.}}
\label{tab:nbrs}
\end{table}
% $357\,003$ in dimension 28

\begin{remark} 
\label{remintro:redmass}
{\rm (Reduced mass)} 
Assume we are interested in the groupoid $\mathcal{G}$ of unimodular lattices $L$ having a given rank and root system $R$. For each such $L$, then $|{\rm W}(R)|$ divides $|{\rm O}(L)|$, where ${\rm W}(R)$ denotes 
the Weyl group of $R$. Thus it is often more natural to multiply the mass of $L$ (and $\mathcal{G}$) by $|{\rm W}(R)|$ :
we call this the {\it reduced mass} of $L$ or $\mathcal{G}$ (see\,\S \ref{subsect:rootsystems}).
\end{remark}

\subsection{The visible part of a $d$-neighbor of ${\rm I}_n$}
\label{subsect:visiblepart}
	
	Although this search root system by root system is necessary, it is still by far not enough to find all unimodular lattices. Indeed, the number of ${\rm O}({\rm I}_n)$-orbits of $d$-isotropic lines in ${\rm I}_n$ is $\geq \frac{d^{n-2}}{n! 2^{n-1}}$, and in practice it is very lengthy to run over all $d$-isotropic lines already for $d$ about $30$ in our range, whereas the farness of many unimodular lattices is much bigger than $30$.
In any case, it would be ridiculous to enumerate naively all isotropic lines !  Indeed, 
certain ''{\it visible}'' properties of the neighbors, in the sense that they can be directly read off from the isotropic lines defining them, 
substantially biase our search and suggest more clever choices of isotropic lines. This is one of the main topics of this paper. \ps	

%This way we immediately get at out disposal 
%a large quantity of unimodular lattices, which is quite satisfying, but certainly not all, and in practice 
%the last $d$ such that it is reasonable to study all $d$-neighbors of ${\rm I}_n$ is much smaller than the $d$ 
%we need to reach to find all of them. 
%\begin{definition} 
%Let $L \in \mathcal{L}_n$ and let $N$ be a $d$-neighbor of $L$.  Set $M=L \cap N$.
%We define the {\rm visible root system} of $N$ as ${\rm R}_2(M)={\rm R}_2(L) \cap {\rm R}_2(N)$, and {\it the visible isometry group} of $N$ as ${\rm O}(L) \cap {\rm O}(M)$.
%
%\end{definition}
%\begin{definition} \label{def:vispart} 
%For $L \in \mathcal{L}_n$, the {\it visible part} of a $d$-neighbor $N$ of $L$ is the lattice 
%$M:=N\cap L$; it is a sublattice of index $d$ in both $L$ and $N$. 
%\end{definition}

For $i\geq 1$ and $N$ a $d$-neighbor of ${\rm I}_n$, there is usually a part of mystery in the Euclidean 
set ${\rm R}_i(N)$, but what we do control is its subset ${\rm R}_i(M) = {\rm R}_i({\rm I}_n) \cap {\rm R}_i(N)$.
We have $M={\rm M}_d(x)$ for some ($d$-isotropic) $x \in \Z^n$ by Formula~\eqref{defmx}, and in the canonical basis $\varepsilon_1,\dots,\varepsilon_n$ of $\R^n$, we also have 
\begin{equation}
\label{eq:r1r2in}
\begin{array} {cccccc}
{\rm R}_1({\rm I}_n) & = & \{ \pm \varepsilon_i\,\,|\,\, 1\leq i \leq n\} & & ,\\
{\rm R}_2({\rm I}_n) & = & \{\pm \varepsilon_i \pm \varepsilon_j, \, \, 1\leq i<j\leq n\} & \simeq & {\bf D}_n.
\end{array}
\end{equation}
For instance, we see ${\rm R}_1({\rm M}_d(x))= \{ \pm \varepsilon_i \,\,|\,\, x_i \equiv 0 \bmod d\}$.
In particular, we have ${\rm r}_1({\rm M}_d(x))=0$ if, and only if, $x_i \not\equiv 0 \bmod d$ for all $i$,\footnote{If $x$ is $d$-ordered, {\it i.e.} as in \eqref{mainineq}, this is equivalent to $x_1 \geq 1$.} an assumption that we shall always make since we are only interested in neighbors $N$ with ${\rm r}_1(N)=0$.
The root system ${\rm R}_2({\rm M}_d(x))$ is also {\it visible}, in the sense that it immediately follows from 
an inspection of the $(i, j)$ with $1 \leq i \leq j \leq n$ and $x_j \equiv \pm x_i \bmod d$ : see Section~\ref{sec:visiblers}. It has the form
{\small
\begin{equation} \label{eq:visibrs} 
{\rm R}_2(M) \simeq {\bf A}_{n_1-1} {\bf A}_{n_2-1} \cdots {\bf A}_{n_s-1} {\bf D}_m, \, \, \, {\rm with}\, \, \, n=n_1+n_2+\dots+n_s+m.
\end{equation}
}
\begin{definition} \label{def:visrs} 
For $L \in \mathcal{L}_n$, and $N$ a $d$-neighbor of $L$, the {\rm visible root system} of $N$ is 
${\rm R}_2(M)$ with $M=L \cap N$. It is a sub root system of ${\rm R}_2(N)$, namely
${\rm R}_2(N) \cap L$. 
\end{definition} 

Properties of the visible root system are studied in Section~\ref{sec:visiblers}.
Other visible objects will be studied in this paper and play some important roles, such as {\it visible isometries} 
in Section~\ref{sect:visisometry}, and {\it visible exceptional vectors} in Section~\ref{sec:excsec}.

\begin{exampleintro} 
\label{ex:emptyrs}
{\rm (Empty root system)} 
As a trivial example, consider the problem of searching for 
$d$-neighbors $N={\rm N}_d(x')$ of ${\rm I}_n$ 
with ${\rm R}_{\leq 2}(N) =\emptyset$. 
For such an $N$, the visible root system has of course to be empty, 
which forces $x_i \neq 0$ and $x_i \neq \pm x_j \bmod d$ for all $1\leq i <j\leq n$. 
In particular we have $d \geq 2n+1$. 
Better, this shows that for $d=2n+1$ the unique possibility 
up to isometry is the lattice ${\rm N}_{2n+1}(1,2,\dots,n)$ 
for $n \not \equiv 1 \bmod 3$ already introduced in \eqref{eq:oddshortborchleech}. 
A similar reasoning, taking into account that for ${\rm N}_d(x')$ to be even 
we must have $d$ even and all coordinates of $x$ odd, 
immediately leads to Thompson's construction of {\rm Leech}. 
\end{exampleintro}

\subsection{The Biased Neighbor Enumeration algorithm}
\label{subsect:biasedalgo}	
	
	Consider now the problem of searching for unimodular lattices in $\mathcal{L}_n$ with a given (arbitrary) root system $R$. The basic idea would be to restrict the enumeration (M3) to $d$-isotropic lines $x$ such that the visible root system of ${\rm N}_d(x')$ is $R$. However, this cannot not work in general for at least two different reasons. \ps

-- First of all, $R$ may not be of the form of the right hand side of \eqref{eq:visibrs} at all. This happens either if it contains some component of type ${\bf E}$, of several components of type ${\bf D}$, or if the union of its ${\bf A}$ components have a too large total rank. For instance, in dimension $n$ there is no visible root system of type $k\,{\bf A}_1$ with $k>n/2+1$, although of course there may be unimodular lattices with such a root system. \ps

-- Worse, even if $R$ may occur as a visible root system in dimension $n$, it may be the case that certain unimodular lattices of rank $n$ and root system $R$ cannot be obtained as a neighbor with visible root system $R$. 
One reason for this is that the visible root system of $N={\rm N}_d(x')$ is always saturated in ${\rm I}_n$, hence closed to be so in $N$. More precisely, the visible root system of $N$ is {\it a $d$-kernel} of ${\rm R}_2(N)$ in our terminology : see \S~\ref{subsect:satprop}. For this strategy, it becomes important to classify all $d$-kernels of ${\bf ADE}$ root systems, we do so in \S~\ref{subsect:dkerdesc} using properties of affine Weyl groups. \ps

	The good news is that the two obstacles above are essentially the only constraints. 
Indeed, this is a special case of \cite[Thm. 7.1]{chestat}, which extends Thm.~\ref{thmi:stat1}, and that we now try to state in its simplest form. For $R$ a root system and $L' \in \mathcal{L}_n$ we denote by ${\rm emb}(R, L')$ the set of isometric embeddings\footnote{We denote by ${\rm Q}(R)$ the root lattice generated by $R$.} ${\rm Q}(R) \rightarrow L'$ with saturated image and odd orthogonal complement. We also define ${\rm m}_n^{\rm odd}(R)$ as the mass of the groupoid of pairs $(L,\iota)$ with $L \in \mathcal{L}_n$ and $\iota \in {\rm emb}(R,L)$. The following result follows from the special case $A={\rm Q}(R)$ of \cite[Thm. 7.10]{chestat}, as well as Corollary 7.12 {\it loc. cit.}

\begin{thmintro} 
\label{thmi:stat2} 
Let $L,L'$ be odd unimodular lattices of rank $n$. Assume  $R \subset {\rm R}_2(L)$ 
is a rank $r$ saturated sub root system whose orthogonal complement in $L$ contains a sublattice isometric to ${\rm I}_3$.
For $p$ a prime, let ${\rm n}_p(L, L', R)$ be the number of $p$-neighbors of $L$ isometric to $L'$ 
and with visible root system containing $R$. Then we have 
$$\frac{{\rm n}_p(L, L', R)}{p^{n-r-2}} \,\longrightarrow\, \frac{ |{\rm emb}(R,L')|\,\, {\rm mass}\, (L')}{{\rm m}_n^{\rm odd}(R)}\,\,\,{\rm when}\,\,\,\, p \,\,\rightarrow \,\,+\infty.$$
\end{thmintro}

In other words, by prescribing the visible root system to contain $R$ we biase the statistics of Theorem~\ref{thmi:stat1} exactly by the factor $|{\rm emb}(R,L')|$. 
We apply this result to $L={\rm I}_n$.
Prescribing the visible root system in the enumeration is then immediate : it just amounts to restrict to $d$-isotropic lines having certain equal coordinates. In practice, this method is extremely efficient, and allows in only a few seconds to find a unimodular lattice with given root system, and even all of those lattices $L'$ having the largest possible $|{\rm emb}(R,L')|\,\, {\rm mass}\, (L')$. As an example, {\it it allows to reproduce the full list of all unimodular lattices of rank $\leq 25$ in just a few minutes}. It is sometimes delicate to understand which visible root system is the best to use in our search : see \S~\ref{subsect:safe} \& \ref{subsect:dkerdesc} for this question, which often boils down to exercises in root systems and coding theory. \ps

We give a concrete example of application of this Biased Neighbor Enumeration algorithm in dimension $26$ in the next section, and others in \S~\ref{sect:exdim26}. We also refer to \S 4 of the companion paper~\cite{allche} for a more formal exposition of this algorithm, and to {\it loc. cit.} \S 5 for several other examples. \ps

\subsection{ An example : the root system $10\,{\bf A}_1$ in ${\rm X}_{26}$}
\label{subsect:example10A1}	

Let us consider the problem of finding all unimodular lattices of rank $26$ with root system $10\,{\bf A}_1$ (and 
with no norm $1$ element, we will not repeat this condition). By King, the reduced mass of this groupoid of lattices is \scalebox{.8}{$4424507/116121600$} (see Remark~\ref{remintro:redmass}). Moreover, we may show that any lattice with root system $10\,{\bf A}_1$ contains a saturated $8\,{\bf A}_1$ (see \S~\ref{subsect:safe} and especially Example~\ref{ex:safena1}), 
so it looks safe to search for our lattices with such a visible root system. This means we restricts to enumerating $d$-isotropic $x \in \Z^{26}$
satisfying \eqref{mainineq} and with exactly {\it $8$ pairs of equal coordinates} modulo $d$. In particular, this forces $d\geq 2 \cdot(8+10)=36$. \ps
	For $d=36$ we do instantly find $108$ such isotropic lines, and $49$ of them happen to lead to a neighbor with root system $10 {\bf A}_1$ (and no norm $1$ element). This high ratio $49/108$ attests that our bias is successful : we chosed very well the visible root system. Those $49$ lattices happen to contain $4$ different isometry classes, with respective reduced mass \scalebox{.8}{$1/64$}, \scalebox{.8}{$1/96$}, \scalebox{.8}{$1/96$} and \scalebox{.8}{$1/640$} : see Table~\ref{tab:10A1dim26} for representatives. \ps
	For $d=37$ and $38$ the ratio of $d$-neighbors with root system $10\, {\bf A}_1$ are respectively \scalebox{.8}{$41/73$} and \scalebox{.8}{$458/1095$}, but no new lattices are found. For $d=39$, the ratio is \scalebox{.8}{$820/1821$} and we find the fifth lattice below, with reduced mass \scalebox{.8}{$1/12288$} : this later finding fits the fact that the probability to find this lattices was smaller, as so is its mass, by Theorem~\ref{thmi:stat2}. \ps Unfortunately, for $40 \leq d \leq 49$, a systematic enumeration of more than $2$ millions of isotropic lines do not lead to any new lattice. This is not really a surprise. Indeed, the remaining reduced mass is
$$\scalebox{.8}{4424507/116121600} - \scalebox{.8}{1/64} -\scalebox{.8}{1/96}-\scalebox{.8}{1/96}-\scalebox{.8}{1/640}-\scalebox{.8}{1/12288} = \scalebox{.8}{17/116121600}.$$
This is only about $4 \cdot 10^{-6}$ of the initial mass, hence extremely small, so unless we have at our disposal many cores (and time to waste), it is not reasonable to search for the remaining lattices just by pursuing the enumeration in the ``coupon collector'' style (see~\S\ref{subsect:coupon}).  It is thus highly desirable to have other methods to find the remaining lattices. We propose two methods here : one that we call the {\it adding ${\bf D}_m$ method}, and that we will explain in details in Sect.~\ref{sect:lattconstr} (see especially~Prop.\ref{prop:adddn} and Remark~\ref{rem:addn}), and another one based on visible isometries, studied in Sect.~\ref{sect:visisometry}. \par

The first idea is to observe that if $L$ is a rank $26$ unimodular lattice with root system $10\,{\bf A}_1$ and no norm $1$ vector, then the orthogonal of any $2 {\bf A}_1$ in $L$ has root system $8 {\bf A}_1$, rank $24$, and index $2$ in some unimodular lattice. For our $L$ of interest, this lattice has a presumably high farness, and we do see in rank $24$ an odd unimodular lattice $L_0$ with root system $8\,{\bf A}_1$ appearing in the end of our lists. It has the reduced mass \scalebox{.8}{$1/20643840$}, which is close enough to \scalebox{.8}{$1/116121600$} since the quotient is \scalebox{.8}{$45/8$},
so this is promising. Actually, we can immediately discover this lattice using the visible root system $7 {\bf A}_1$ in rank $24$ : we have {\it e.g.}  $L_0 \simeq {\rm N}_{35}(x)$
for 
$$x=\scalebox{0.8}{({\color{magenta} 1, 1}, 2, 3, {\color{magenta}4, 4}, 5, 6, {\color{magenta}7, 7}, 8, {\color{magenta}9, 9}, 10, 11, {\color{magenta}12, 12}, 13, 14, {\color{magenta}15, 15}, 16, {\color{magenta}17, 17})}\in \Z^{24}.$$ This suggests to look at 
the $2$-neighbors of $L_0 \oplus {\rm I}_2$, which are specific $70$-neighbors of ${\rm I}_{26}$. 
By restricting to those with same visible root system $7 {\bf A}_1$, there are now only $2^{16}$ lines to check. 
We find $16384$ isotropic lines, $1816$ leading to a neighbor with root system $10 {\bf A}_1$. Among those, two lattices are found, with reduced masses \scalebox{.8}{$1/7372800$} and \scalebox{.8}{$1/92897280$} respectively, which fulfills the mass : see Table~\ref{tab:10A1dim26} below.  

\tabcolsep=2.5pt
\begin{table}[H]
{\scriptsize 
\renewcommand{\arraystretch}{1.8} \medskip
\begin{center}
\begin{tabular}{c|c|c|c}
$d$ & $x \in \Z^{26}$ & ${\rm reduced \,mass}$ & $|{\rm emb}(8{\bf A}_1,-)|$ \\ \hline
$36$ & $\scalebox{0.7}{({\color{magenta}1, 1}, 2, 3, 4, 5, {\color{magenta}6, 6}, {\color{magenta}7, 7}, {\color{magenta}8, 8}, 9, 10, 11, {\color{magenta}12, 12}, {\color{magenta}13, 13}, {\color{magenta}14, 14}, 15, {\color{magenta}16, 16}, 17, 18)}$ & $1/64$ & $45$ \\
$36$ & $\scalebox{0.7}{({\color{magenta}1, 1}, 2, 3, 4, 5, {\color{magenta}6, 6}, 7, 8, 9, {\color{magenta}10, 10}, {\color{magenta}11, 11}, 12, {\color{magenta}13, 13}, {\color{magenta}14, 14}, 15, {\color{magenta}16, 16}, 17, {\color{magenta}18, 18})}$ & $1/96$ & $45$\\
$36$ &  $\scalebox{0.7}{({\color{magenta}1, 1}, 2, 3, 4, 5, 6, 7, {\color{magenta}8, 8}, 9, 10, {\color{magenta}11, 11}, {\color{magenta}12, 12}, {\color{magenta}13, 13}, {\color{magenta}14, 14}, 15, {\color{magenta}16, 16}, 17, {\color{magenta}18, 18})}$ & $1/96$ & $45$\\
$36$ & $\scalebox{0.7}{({\color{magenta}1, 1}, {\color{magenta} 2, 2}, 3, 4, 5, {\color{magenta}6, 6}, 7, {\color{magenta}8, 8}, {\color{magenta}9, 9}, {\color{magenta}10, 10}, 11, 12, 13, {\color{magenta}14, 14}, 15, 16, {\color{magenta}17, 17}, 18)}$ & $1/640$ & $45$\\
$39$ & $\scalebox{0.7}{({\color{magenta}1, 1}, 2, 3, 4, {\color{magenta} 5, 5}, 6, {\color{magenta}7, 7}, {\color{magenta}8, 8}, 9, {\color{magenta}10, 10}, 11, 12, 13, {\color{magenta}14, 14}, 15, {\color{magenta}16, 16}, 17, {\color{magenta}19, 19})}$ & $1/12288$ & $45$\\
$70$ & $\scalebox{0.7}{({\color{magenta}1, 1}, 33, 3, {\color{magenta} 4, 4}, 5, 29, {\color{magenta}7, 7}, 27,{\color{magenta}9, 9}, 25, 11, {\color{magenta}12, 12}, 13, 21, {\color{magenta}15, 15}, 19, {\color{magenta}18,  18}, {\color{magenta} 35, 35})}$ & $1/7372800$ & $25$ \\
$70$ & $\scalebox{0.7}{({\color{magenta}1, 1}, 33, 3, {\color{magenta} 4, 4}, 5, 29, {\color{magenta}7, 7}, 27,  {\color{magenta} 9, 9}, 25, 11, {\color{magenta} 23, 23}, 13, 21, {\color{magenta}15, 15}, 19, {\color{magenta}17, 17}, {\color{magenta} 35, 35})}$ & $1/92897280$ & $9$ \\
\end{tabular} 
\end{center}
}
\caption{{\small The $7$ lattices with no norm $1$ elements and root system $10 {\bf A}_1$ in ${\rm X}_{26}$.}}
\label{tab:10A1dim26}
\end{table}

	Another way to find the last lattices in Table~\ref{tab:10A1dim26} amounts to use the theory of visible isometries. 
Indeed, we  have $116121600=2^{13}\,3^4\,5^2\,7$, so we knew for instance that some lattice had an order $7$ automorphism (namely, the last in the table above).  Applying the method of Section~\ref{sect:visisometry} for searching for elements of ${\rm  X}_{26}$ with root system $10\,{\bf A}_1$ and a visible automorphism with characteristic polynomial $\Phi_7^3 \Phi_1^{8}$ immediately leads to many constructions of this lattice, such as 
one\footnote{The corresponding $x \in \Z^{26}$ is congruent to $(1^{14}, 3^7, 4, 8, 9, 10, 11)$ mod $27$ and to $(1, 16, 24, 7, 25, 23, 20, 1, 16, 24, 7, 25, 23, 20, 1, 16, 24, 7, 25, 23, 20, 0^5)$ mod $29$.} for $d=29 \cdot 27=783$. Here $29$ is the prime $\equiv 1\bmod 7$. 

\subsection{A coupon collector's problem and choices of invariants}
\label{subsect:coupon}	
${}^{}$\indent Assume we want to find all unimodular lattices
	of given rank $n$ by running through all isotropic lines in ${\rm I}_n \otimes \Z/d$, with $d=2, 3, ...$.
%	In other words, we want to complete our ``Panini album'' of unimodular lattices, and
	Theorem~\ref{thmi:stat1} shows that we are in the situation of a coupon collector's problem (or "Panini album").
	Indeed, studying a new line (which costs some time to the computer) corresponds to buying a new box, and the mass of a lattice is proportional to the probability of find it in a box. Recall that in a uniform situation with $N$ coupons, then on average we need to buy $N\,\, {\rm log}\, N$ boxes to collect all coupons (see e.g. \S 2.4.1 in  \cite{mitzenmacherupfal}). In our non-uniform situation this can be much worse, since certain lattices have much smaller masses than others. By Theorem~\ref{thmi:stat2}, the idea of fixing the visible root system substantially helps reducing this non uniformity, 
and it does work in most cases. Sometimes, however, it is not enough as shown by the Example in Sect.~\ref{subsect:example10A1}, and it is preferable to use other methods to find the lattices with smallest mass. \par
	A bit surprisingly, those statistic arguments have also been very helpful in order to discover that some invariants were not fine enough.
Indeed, when searching for lattices with a given root system using a certain invariant, if start finding new lattices with very small mass compared to the remaining mass for this root system, it is a strong indication that we missed the most likely ones : our chosen invariant is not fine enough. See {\it e.g.} the discussion of the case $R \simeq 2\,{\bf A}_1\,2{\bf A}_2\,2{\bf A}_3\,2\,{\bf A}_4$ in dimension $26$ in Sect.~\ref{sect:exdim26} for an example of such a sitation.\par

\subsection{Proofs of Theorems~\ref{X2627} and~\ref{thm:listuninei} : the full lists}
\label{subsect:fulllists}
${}^{}$ 
	Although most of the lattices in the lists of Theorem~\ref{thm:listuninei} have been found 
by applying the naive algorithm described in (M3) of \S\ref{subsect:generalalgo}, 
the main work was to deal with the remaining ones.  
A large part of them were found by using the biased algorithm of \S\ref{subsect:biasedalgo}, 
using for each remaining root system clever choices of visible root systems as explained in Sect.~\ref{sec:visiblers}. 
The most resisting (and interesting!) lattices were then dealt with more specific methods, 
some of them already encountered in the introduction : enumerating the $2$-neighbors of well-chosen lattices,  ``addition of ${\bf D}_m$'' method (see Sect.~\ref{sect:lattconstr}), separate study of exceptional lattices (see Sect.~\ref{sec:excsec})...\footnote{Actually, when we first made these computations in 2020 we had not yet discovered the theory of visible isometries explained in~\S\ref{sect:visisometry}, and mostly used instead the ideas of \S\ref{subsect:neismallmasshassmallmass} to construct lattices with small reduced mass. Many such lattices can also be found faster using the arguments of~\S\ref{sect:visisometry}  (although with a worst farness!). } \par
 The whole computation required so many case by case considerations that it would be not be reasonable to list them here. For instance, several hundreds of the $2797$ possible root systems in dimension $27$ had to be treated separately (as in Example~\ref{subsect:example10A1}). Instead, our expository choice in this paper was to explain the theoretical aspects underlying each method that we used, and to only provide a few detailed examples as illustrations. We refer to the companion paper~\cite{allche} for more examples of our method, including some emphasize (and improvements) on some computational aspects only briefly discussed in this paper.  We mention that all of our computations have been made using \texttt{PARI/GP} and a processor \texttt{Intel(R) Xeon(R) CPU E5-2650 v4 @ 2.20GHz} with  \texttt{65 GB} of memory; the total CPU time was about $1$ month in dimension $26$ and $1$ year in dimension $27$.\ps
	Actually, if our goal is just to prove Theorems~\ref{X2627} and~\ref{thm:listuninei}, rather than understanding the mathematical ideas involved, it is not really necessary to explain how our lists in \cite{cheweb} were discovered  ! Indeed, it may be checked independently and a posteriori that these lists are complete : it is enough to check that all the given lattices have distinct invariants and that the sum of their masses equals the mass formula. See~\cite{cheweb} for the relevant \texttt{PARI/GP} source code for this check. This is of course much shorter : it only requires $5$ hours in dimension $26$, and $40$ hours in dimension $27$. As a consequence, this also provides an independent verification of King's computations. \ps
	In the companion paper \cite{allche}, in collaboration with Bill Allombert, we pursue the ideas of this paper and determine in particular ${\rm X}_{28}$ and ${\rm X}_{29}^\emptyset$. \ps \ps\ps

{\sc Aknowledgments.} We warmly thank the developers of \texttt{PARI/GP} for their work, as well as the LMO for sharing their machine \texttt{pascaline}. We also thank Bill Allombert and Olivier Ta\"ibi for useful discussions.

\ps\ps
\begin{comment}

\begin{remark}\label{remanormdeven} {\rm (Normalisation for $d$ even)} For a $d$-isotropic $x$ with $x_1=1$ and $d$ even, 
we define ${\rm N}_d(x)^+$ {\rm (}resp. ${\rm N}_d(x)^-${\rm )} using $x' = x - (x.x/2 \,+\, t \, d^2/2)\varepsilon_1$ with $t=0$ {\rm (}resp. $t=1${\rm )} in Formula \eqref{defldx}. These lattice are denoted $[d, [x_1;x_2;\dots;x_n], t]$ in the lists {\rm \cite{cheweb}}, with the convention $t=0$ for $d$ odd. 
\end{remark}
\ps\ps

 A second important motivation for this work was the relationship between genera of integral lattices and  the geometric $\ell$-adic Galois representations of ${\rm Gal}(\overline{\Q}/\Q)$ (or pure motives over $\Q$), which follows from the general yoga and point of views of Langlands and Arthur on automorphic representations. The genera of even unimodular lattices, which yield Galois representations of Artin conductor $1$ in this yoga, has been studied in details from this point of view in \cite{chlannes}. 
As was shown in the recent Ph.D. thesis of Lachauss\'ee \cite{lachaussee}, the genera of unimodular lattices of odd rank yield Galois representations of Artin conductor $2$, and the classification of these representations was a motivation to obtain the case of dimension $27$ in the theorem above. Actually, in a joint work with O. Ta\"ibi, we even hope to use the computations here in a quite indirect way to push further the dimension formulas for level one Siegel modular cuspforms proved in \cite{taibisiegel,chetai}, and to obtain new information on $|{\rm X}_{32}|$.\ps

\end{comment}
  
  {\scriptsize
\tableofcontents
}
  
\section{General conventions and notations}\label{notations}\ps
%\begin{center}{\sc General notations and conventions}\end{center}  
	In this paper, group actions will be on the left. We denote by $|X|$ the cardinality of the set $X$. For $n\geq 1$ an integer, we denote by ${\rm S}_n$ the symmetric group on $\{1,\dots,n\}$, by ${\rm Alt}_n \subset {\rm S}_n$ the alternating subgroup, and we also denote by $\Z/n$ the cyclic group $\Z/n\Z$.\ps
	(i) If $V$ is an Euclidean space, we usually denote by $x \cdot y$ its inner product. A {\it lattice} in $V$ is a subgroup generated by a basis of $V$, or equivalently, a discrete subgroup $L$ with finite covolume, denoted ${\rm covol}\, L$. Its {\it dual lattice} is the lattice $L^\sharp$ defined as $\{ v \in V\, \, |\, \, v \cdot x \in \Z ,\,\, \, \forall x \in L\}$.  Recall that $L$ is {\it integral} if we have $L \subset L^\sharp$.  An integral lattice is called {\it even} if we have $x \cdot x \in 2\Z$ for all $x \in L$, {\it odd} otherwise. The orthogonal group of $V$ is denoted by ${\rm O}(V)$, and we also denote by ${\rm O}(L)=\{ \gamma \in {\rm O}(V), \, \, \gamma(L)=L\}$ the isometry group of $L$ (a finite group).  \ps
	(ii) Assume $L \subset V$ is an integral lattice. The finite abelian group ${\rm res}\, L := L^\sharp/L$ (sometimes called the {\it discriminant group} \cite{nikulin}, the {\it glue group} \cite{conwaysloane} or the {\it residue} \cite{chlannes}) is equipped with a non-degenerate $\Q/\Z$-valued symmetric bilinear form, defined by $(x,y) \mapsto x.y \bmod \Z$. We have $({\rm covol}\, L)^2 \,=\,  |{\rm res}\, L|$. This integer, also denoted ${\rm det}\, L$, is also the {\it determinant} of the {\it Gram matrix} ${\rm Gram} (e) = (e_i \cdot e_j)_{1 \leq i,j \leq n}$ of any $\Z$-basis $e=(e_1,\dots,e_n)$ of $L$. \ps
	(iii) A subgroup $I \subset {\rm res}\, L$ is called {\it isotropic}, if we have $x.y  \equiv 0$ for all $x,y \in I$, and a {\it Lagrangian} if we have furthermore $|I|^2=|{\rm res}\, L|$. The map $\beta_L : M \mapsto M/L$ defines a bijection between the set of integral lattices containing $L$ and the set of isotropic subgroups of ${\rm res}\, L$. In this bijection, $M/L$ is a Lagrangian if and only if $M$ is unimodular. If $I$ is finite abelian group, we denote by ${\rm H}(I)$ the {\it hyperbolic} symmetric bilinear space $I \oplus I^\ast$, with $I^\ast={\rm Hom}(I,\Q/\Z)$, defined by $(x+\phi).(x'+\phi')=\phi(x')+\phi'(x)$.\ps
	(iv) Assume furthermore that $L$ is an even lattice. Then the finite symmetric bilinear space ${\rm res}\, L$ 
has a canonical quadratic form ${\rm q} : {\rm res}\, L \rightarrow \Q/\Z$ such that ${\rm q}(x+y)-{\rm q}(x)-{\rm q}(y) \equiv x.y$,
defined by ${\rm q}(x) = \frac{x.x}{2} \bmod \Z$. In the bijection $\beta_L$ above, the even lattices $M$ correspond to the {\it quadratic} isotropic subspace $I \subset {\rm res}\, L$, {\it i.e} with ${\rm q}(I)=0$. 
We also denote by ${\rm qm} : {\rm res}\, L \rightarrow \Q_{\geq 0}$ the {\it Venkov map}, defined by ${\rm qm}(x) \, =\, {\rm Min}_{ y \in x+ L}\, \, \frac{y.y}{2}$. It satisfies ${\rm qm}(x) \equiv {\rm q}(x) \bmod \Z$, ${\rm qm}(x)={\rm qm}(-x)$ and ${\rm qm}(x) >0$ for $x \not \equiv 0$.\ps

     (v) ({\it Standard lattices}) Here $\R^n$ denotes the standard Euclidean space, for $n\geq 0$, with canonical basis $\varepsilon_1,\dots,\varepsilon_n$. We set
${\rm D}_n=\{ x \in \Z^n\, \, |\, \, \sum_i x_i \equiv 0 \bmod 2\}$ and ${\rm A}_n = \{ x \in \Z^{n+1} \, \, |\, \, \sum_i x_i  = 0\}$. The ${\rm E}_8$ {\it lattice} is ${\rm D}_8 + \Z e$ with $e=\frac{1}{2}\sum_{i=1}^8 \varepsilon_i$, the 
${\rm E}_7$ (resp. ${\rm E}_6$) {\it lattice} is the orthogonal of $\varepsilon_7+\varepsilon_8$ (resp. of $\varepsilon_7-\varepsilon_6$ and $\varepsilon_7+\varepsilon_8$) in ${\rm E}_8$. All these lattices are even. Their Venkov map are well-known (theory of {\it minuscule weights}),  with nonzero values given by Table \ref{tab:venkov}:
\begin{table}[H]
{\scriptsize \renewcommand{\arraystretch}{1.8} \medskip
\begin{center}
\smallskip
\begin{tabular}{c|c|c|c|c|c|c}
$L$ &  ${\rm A}_n$ & ${\rm D}_{n}$, $n>0$\,\,{\rm even} & ${\rm D}_{n}$, $n$\,\,{\rm odd}& ${\rm E}_6$ & ${\rm E}_7$ & ${\rm E}_8$ \cr
\hline
${\rm res}\, L$ & $\Z/(n+1)$ & $\Z/2\times \Z/2$ & $\Z/4$ & $\Z/3$ & $\Z/2$ & $0$ \cr
%$x$ & $i \bmod n+1$ & $\omega, 1, \overline{\omega}$ & $1,2,3 \bmod 4$ & $1,2 \bmod 3$ & $1 \bmod 2$ & \cr
\hline
${\rm qm}$ & $\frac{i(n+1-i)}{2(n+1)}$ {\rm with}\, $1\leq i \leq n$ & \multicolumn{2}{c|}{$\frac{n}{8},\frac{1}{2},\frac{n}{8}$} & $\frac{2}{3}, \frac{2}{3}$ & $\frac{3}{4}$ & \cr
\end{tabular} 
\end{center}
} 

\caption{The nonzero values of ${\rm qm}$ on ${\rm res}\, L$ (with multiplicities).}
\label{tab:venkov}
\end{table}
More precisely, for $L={\rm A}_n$ there is a group isomorphism $\phi : \Z/(n+1) \isomo {\rm res}\, {\rm A}_n$ with 
${\rm qm}(\phi(i \bmod n+1))=\frac{i(n+1-i)}{2(n+1)}$ for $1 \leq i \leq n$. \ps

        (vi) A subgroup $A$ of a lattice $L$ is called {\it saturated} if the abelian group $L/A$ is torsion free, 
or equivalently, if $A$ is a direct summand of $L$ as $\Z$-module. 
The {\it saturation} of $A$ in $L$, defined as ${\rm Sat}_L(A)=L \cap (A \otimes \Q)$, 
is the smallest saturated subgroup $S$ of $L$ containing $A$.\ps

%   (vii) ({\it Euclidean sets}) Define an {\it Euclidean set} as a set $X$ equipped with 
%an injection $j : X \hookrightarrow V$ into some Euclidean space $V$. 
%We then denote by $\R X$ the Euclidean subspace of $V$ generated by $j(X)$.
%Euclidean sets form a category if we define a morphism $(X,j) \rightarrow (X',j')$ as 
%a map $f: X \rightarrow X'$ {\it induced} by a linear isometric  embedding
%$\widetilde{f} : \R X \rightarrow \R X'$, {\it i.e.} verifying\footnote{Note that such an $\widetilde{f}$ is unique if it exists.} $\widetilde{f} \circ j = j' \circ f$.\ps

% It make sense to talk about the {\it scalar product} $x.y$ of two elements $x,y$ of an Euclidean set $(X,j)$, namely $x.y=j(x).j(y)$, to define the {\it rank} of $X$ as $\dim_\R \R X$, 
%the lattice $\Z X$ about the {\it lattice generated by $X$} ({\it i.e.} $\Z {\rm j}(X)$) ...\ps

\section{Cyclic Kneser neighbors of unimodular lattices}\label{sectkneser}

\subsection{Definitions and notations} Let $V$ be an Euclidean space, $L$ a unimodular integral lattice in $V$, and $d\geq 1$ an integer. 
A {\it $d$-neighbor} of $L$ is a unimodular integral lattice $N \subset V$ such that we have a group isomorphism $$L/(L \cap N) \simeq \Z/d.$$
Such lattices $N$ are sometimes called {\it cyclic} $d$-neighbors of $L$, but from now on we will omit the adjective cyclic for short. A few remarks are in order:\ps\ps
	
(Na) If $N$ is a $d$-neighbor of $L$, then $L$ is a $d$-neighbor of $N$. Indeed, if we set  $M=L \cap N$, then we have $M^\sharp = N^\sharp+L^\sharp=N+L$, and thus ${\rm res}\, M = N/M \oplus L/M$ with $N/M$ and $L/M$ Lagrangians. 
As a consequence, the pairing of ${\rm res}\, M$ identifies $N/M$ with ${\rm Hom}(L/M,\Q/\Z) \simeq \Z/d$. \ps\ps

(Nb) As $L$ is unimodular, the subgroups $M \subset L$ with $L/M \simeq \Z/d$ are the 
$${\rm M}_d(L; x) := \{ m \in L\,\,| \, \, m . x \equiv 0 \bmod d\}$$
% \,\, \hspace{.5 cm}\,\,(\,\,=  ``\,\, x^{\perp \bmod d}\,\,{\textrm{''}})$$
where $x \in  L$ is a {\it $d$-primitive} vector. By this mean an element of $L$ 
whose image in $L/dL$ generates a subgroup of order $d$. For $d$-primitive $x,x' \in L$,
we have ${\rm M}_d(L; x) ={\rm M}_d(L; x')$ if, and only if, 
$x$ and $x'$ generate the same subgroup in $L/dL$.
We thus denote as well by ${\rm M}_d(L ; \ell)$ 
the lattice ${\rm M}_d(L; x)$, if $\ell \simeq \Z/d$ is the subgroup in $L/dL$ generated by $x$.\ps\ps

%For $d \geq 1$, $\ell \subset L/dL$ a {\it $\Z/d$-line}, {\it i.e.} a cyclic subgroup of order $d$, and $x \in L$ 
%generating $\ell$ mod $dL$, the element $x.x \bmod ed$ does not depend on the choice of $x$, and will be denoted by $\ell\cdot \ell$. 
% A more precise description is given by the following.

\begin{propdef}\label{classdnei} Let $L$ be a unimodular integral lattice in $V$, $d\geq 1$ an integer, $x \in L$ a $d$-primitive element. Set $M={\rm M}_d(L; x)$, and set also $e=1$ for $d$ odd and $e=2$ for $d$ even.
% and choose $y \in L$ with $x.y \equiv 1 \bmod d$. 
\begin{itemize}
\item[(i)] If $x.x \equiv 0 \bmod ed$, there are exactly $e$ cyclic $d$-neighbors $N$ of $L$ with $N \cap L = M$, and 
none otherwise. These $d$-neighbors are the 
\begin{equation}\label{ndxprime} M + \Z \, \frac{\tilde{x}}{d}\end{equation}
where $\tilde{x}$ is any element of $L$ with $\tilde{x}\equiv x \bmod dL$ and $\tilde{x}.\tilde{x} \equiv 0 \bmod d^2$.\ps
\item[(ii)]
For $d$ odd, the lattice \eqref{ndxprime} does not depend on the choice of $\tilde{x}$ and we denote it ${\rm N}_d(L; x)$.
For $d$ even, it only depends on the element $\epsilon \in \Z/2$ defined by $\tilde{x}.x\equiv \frac{x.x}{2}+ \epsilon \frac{d^2}{2} \bmod d^2$, and we denote it ${\rm N}_d(L; x; \epsilon)$.
\end{itemize}
\end{propdef}

\begin{pf} 
As $x$ is $d$-primitive and $L$ is unimodular, we may and do choose $y \in L$ with $y.x \equiv 1 \bmod d$. 
We clearly have $L \subset M^\sharp$ and $x/d \in M^\sharp$.
Using $L/M \,= \,\Z/d\, \,y$, $|{\rm res}\, M|=d^2$ and 
$x.y \equiv 1 \bmod d$, we obtain ${\rm res}\, M = \, \Z/d \,y\,  \oplus \, \Z/d \,\,\frac{x}{d}$. \par
The $d$-neighbors $N$ of $L$ with $L\cap N = M$
are in bijection with the Lagrangians $I \subset {\rm res}\, M$ 
which are {\it transversal} to $L/M$, {\it i.e.} with $I \cap (L/M) =\{0\}$. 
As $|I|=d$, each such Lagrangian is necessarily $\simeq \Z/d$ and generated by 
a unique element of the form $x/d \,-\, ry $ with $r \in \Z/d$. 
But such an element is isotropic if, and only if, we have $\frac{x.x}{d} \equiv 2r  \bmod d\Z$.
For $d$ odd, there is a unique possibility for $r$, and for $d$ even there are none if $x.x \not\equiv 0 \bmod 2d$,
and exactly two otherwise, namely the
$r_\epsilon = \frac{x.x}{2d} + \epsilon \frac{d}{2} \bmod d$ with $\epsilon \in \Z/2$. 
In the latter case, note that we have 
$$x.(x/d \,-\, r_\epsilon y) \equiv \frac{x.x}{d}- r_\epsilon \equiv 
\frac{x.x}{2 d} + \epsilon \frac{d}{2} \bmod d\Z.$$
This proves the first assertion of (i), as well as the second once we observe that 
for any $\tilde{x} \in x+dL$ we have $\tilde{x}/d \in M^\sharp$, and thus 
$\tilde{x}.\tilde{x} \equiv 0 \bmod d^2$ if and only if $\tilde{x}/d$ is isotropic in ${\rm res}\, M$.
%For $d$ even, the analysis above shows that for such an $\widetilde{x}$ we have
%$\widetilde{x} /d \equiv  x/d- r_\epsilon y \bmod M$, and proves assertion (ii).
\end{pf}

\begin{remark} \label{remdefline}
{\rm \begin{itemize}
\item[(i)] By the proof, we may define ${\rm N}_d(L, x)$ (resp. ${\rm N}_d(L; x; \epsilon))$ by taking 
$\tilde{x} \,=\, x\, + r d y$ in \eqref{ndxprime} with $r=\frac{d+1}{2d}\, x.x$ (resp. 
$r=-\,\frac{x.x}{2d} + \epsilon \frac{d}{2}$). Here $y$ denotes any element of $L$ with $x.y \equiv 1 \bmod d$; it is unique modulo $M$, so those lattices do not depend on this choice of $y$. \ps
\item[(ii)] For $d$ odd, the lattice ${\rm N}_d(L ; x)$ only depends on $M$, hence on the $\Z/d$-line $\ell$ of $L/dL$ generated by $x$, so we may also denote it by ${\rm N}_d(L ; \ell)$. This fails for $d$ even. Indeed, if we set $x'  = x + d v$ with $v \in L$, a simple computation using {\it e.g.} the formula given in (i) shows
${\rm N}_d(L;\, x';\, \epsilon) = {\rm N}_d(L;\, x;\, \epsilon + v.v)$.  
% One way to see this is to observe that we have $z/d \equiv x/d \,+ \,(v.x)\, y \bmod M$,
% and thus $z/d - \frac{z.z}{2} y = x/d - (\frac{z.z}{2} - v.x)y \bmod M, 
% but z.z/2 = x.x /2+ d z.x + d^2/2 v.v
In the case $L$ is even, and only in this case, the lattice ${\rm N}_d(L ; x; \epsilon)$ only depends on $(L; \ell; \epsilon)$ (see Proposition \ref{parnei} (ii)).
\end{itemize}
}
\end{remark}

We denote by $\mathcal{N}_d(L)$ the set of all $d$-neighbors of $L$. 
As the sublattice $L\cap N$, for $N \in \mathcal{N}_d(L)$, plays an important role in this paper, we give a name to it.

\begin{definition} If $N$ is a $d$-neighbor of $L$, we call the lattice $M=L\cap N$ the {\rm visible part} of $N$.
By {\rm (Na)}, this is an integral lattice with ${\rm res}\, M \simeq {\rm H}(\Z/d)$.
\end{definition}

\subsection{The quadric ${\rm C}_L$}

Fix $L \in \mathcal{L}_n$, $d\geq 1$ and set again $e=1$ if $d$ is odd and $e=2$ otherwise.  Consider the finite quadric
\begin{equation}\label{defcl} {\rm C}_L(\Z/d) = \{ \ell \subset L \otimes \Z/d\, \, \, |\, \, \ell \simeq \Z/d \, \,\, \& \,\,\,\ell.\ell \equiv 0 \bmod ed \}.
\end{equation}
For $N$ in $\mathcal{N}_d(L)$, we denote by $l(N)$ the unique $\ell \in {\rm C}_L(\Z/d)$ satisfying ${\rm M}_d(\ell)=N \cap L$.
Alternatively, we have\footnote{Indeed, the image $\ell'$ of $dN$ in $L/dL$ is isomorphic to $N/(N\cap L) \simeq \Z/d$ and satisfies $dN.dN \equiv 0 \bmod ed\Z$, so we have $\ell' \in {\rm C}_L(\Z/d)$. But we have $dN.M \equiv 0 \bmod d\Z$ and thus $M={\rm M}(\ell) \subset {\rm M}(\ell')$, 
hence an equality and $\ell=\ell'$.}
$l(N)=(dN+dL)/dL$. The map
\begin{equation} \label{deflinemap} l : \mathcal{N}_d(L) \rightarrow {\rm C}_L(\Z/d), \, \, N \mapsto l(N) \end{equation}
will be called {\it the line map}. Proposition \ref{classdnei} (i) asserts:

\begin{cor}\label{cor:bijlinemap} The line map is $e:1$ and surjective.
\end{cor} \ps\ps

It would be easy to give a close formula for $|{\rm C}_L(\Z/d)|$, hence for $|\mathcal{N}_d(L)|$. 

\subsection{Parity of a $d$-neighbor} 
\label{subsect:paritynei}
We now discuss the {\it parity} of a neighbor, and the related notion of characteristic vectors. Recall that an integral lattice $L$ is called {\it even} if we have $v . v \in 2\Z$ for all $v \in L$, and {\it odd} otherwise. If $L$ is unimodular, the map $L \rightarrow \Z/2, v \mapsto v.v \bmod 2$ is $\Z$-linear, hence of the form $v \mapsto \xi \cdot v \bmod 2$ for some vector $\xi \in L$, uniquely determined modulo $2L$. Such vectors $\xi \in L$
	are called the {\it characteristic vectors} of $L$; they form a coset in $L/2L$ 
that we will denote by ${\rm Char}(L)$. %For $L= A \perp B$ we also have ${\rm Char}(L)=\{a +b \,\,|\,\, a \in {\rm Char}(A), b \in {\rm Char}(B)\}$.
As an example, we have 
\begin{equation}
\label{eq:charIn} 
{\rm Char}({\rm I}_n) = \{ (\xi_1,\dots,\xi_n) \in \Z^n\, \, |\, \, \xi_i \equiv 1 \bmod 2,\,\,\,1 \leq i \leq n\}.
\end{equation}
We have ${\rm Char}(L) = 2L$ if and only if $L$ is even, and if $L$ is odd and $\xi$ is in ${\rm Char}(L)$, then ${\rm M}_2(L; \xi)$ coincides with the largest even sublattice of $L$. It is clear that if $L$ and $L'$ are $d$-neighbors with $d$ odd, then $L$ is even if and only if $L'$ is even. The case $d$ even is more interesting:

\begin{prop}\label{parnei} Assume $L$ is an integral unimodular lattice, $d$ is even and $x \in L$ is $d$-primitive with $x.x \equiv 0 \mod 2d$. \begin{itemize}
\item[(i)] If $L$ is odd, the $d$-neighbor ${\rm N}_d(L; x; \epsilon)$ is even if, and only if, $x$ is a characteristic vector of $L$ and satisfies $ \frac{x. x}{2d} \equiv  (1+ \frac{d}{2})\epsilon \bmod 2$. \ps \item[(ii)]  If $L$ is even, then ${\rm N}_d(L; x; \epsilon)$ is even if and only if $\epsilon=0$.
\end{itemize}
\end{prop}

\begin{pf} 
By definition, ${\rm N}_d(L; x; \epsilon)\,=\,M\,+\,\Z \frac{\tilde{x}}{d}$ is even if and only if the lattice $M={\rm M}_d(L;x)$ and the integer $\tilde{x}.\tilde{x}/d^2$ are even, with $\tilde{x}$ as in Remark \ref{remdefline} (i). \par 
Assume first $L$ is odd.
As $L/M$ is cyclic of even order, there is a unique lattice $M \subset H \subset L$ with $L/H=\Z/2$, namely $H={\rm M}_2(L;x)$.  It follows that $M$ is even if and only if ${\rm M}_2(L ; x)$ is the largest even sublattice of $L$, i.e. if $x$ is a characteristic vector of $L$. A trivial computation, using $y . y \equiv x. y \equiv 1 \bmod 2$, then shows $\frac{\tilde{x}.\tilde{x}}{d^2} \equiv \frac{x.x}{2 d} + (1+\frac{d}{2}) \epsilon \bmod 2$, and concludes the proof of (i). \par 
If $L$ is even, a simple computation shows $\tilde{x}.\tilde{x}/d^2 \equiv \epsilon \bmod 2$, hence (ii).
\end{pf}

\subsection{Orbits} As emphasized in the introduction, it is a difficult question in general to understand the isometry classes of the $d$-neighbors of a given $L$. A standard observation though is that the isometry group ${\rm O}(L)$ of $L$ naturally acts on $\mathcal{N}_d(L)$, and so two neighbors in a same orbit are isometric. This group also acts on ${\rm C}_L(\Z/d)$ and we have the following obvious proposition.
\begin{prop} \label{propoequ} The line map \eqref{deflinemap} is ${\rm O}(L)$-equivariant. \end{prop}
Equivalently, we have $g({\rm M}_d(x))={\rm M}_d(g(x))$ for all $x \in L$ and $g \in {\rm O}(L)$.
For $d$ odd, the isometry class of a neighbor $N$ depends thus only on the ${\rm O}(L)$-orbit of its line $l(N)$. For $d$ even, the same holds up to the $\epsilon$ ambiguity; more precisely, by assertion (i) of Remark \ref{remdefline} we have for all $g \in {\rm O}(L)$, all $d$-primitive $x \in L$ with $x.x \equiv 0 \bmod 2d$, and all $\epsilon \in \Z/2$, the formula \begin{equation} \label{actOLNd} g({\rm N}_d(L; x; \epsilon)) = {\rm N}_d(L;g(x);\epsilon). \end{equation} 

% Let us finally define, for $d\geq 1$ and two unimodular lattices $L$ and $L'$, 
%the integer ${\rm n}_d(L,L')$ as the number of $d$-neighbors of $L$ which are isometric to $L'$. 
%The following relation due to Eichler follows from Remark (Na) above 
%\begin{equation} \label{symnpl} {\rm n}_d(L,L')|{\rm O}(L')|={\rm n}_d(L',L)|{\rm O}(L)| \end{equation}
%(see e.g. Scholium 3.1.7 of \cite{chlannes} for a ``one line'' proof).

%  Moreover, we have the equality ;
 % this sometimes imposes some relations between $|{\rm O}(L)|$ and $|{\rm O}(L')|$ when 
%  ${\rm N}_p(L,L')$ is small.   

\subsection{$d$-neighbors of ${\rm I}_n$}\label{subsect:dneiin}

%specific notations of I_n

${}^{}$\indent We finally specify the previous considerations to the standard odd unimodular lattice $L={\rm I}_n$, 
and relates the general definitions in this case to the notations already introduced in \S ~\ref{subsecintro:cyclicdnei}. Fix $d \geq 1$. The finite bilinear space ${\rm I}_n \otimes \Z/d$ is just the standard $\Z/d$-valued inner product on the space $(\Z/d)^n$,
and we set ${\rm C}_n(\Z/d) = {\rm C}_{{\rm I}_n}(\Z/d)$.  
The element $x\in \Z^n$ is $d$-primitive if, and only if, we have ${\rm gcd}(x_1,\dots,x_n,d)=1$.
The line ${\rm l}(x) \subset {\rm I}_n \otimes \Z/d$ it generates is in ${\rm C}_n(\Z/d)$ if and only if 
Formula~\eqref{defdiso} holds, {\it i.e.} $x$ is $d$-isotropic. We ease the notations by denoting 
by $${\rm M}_d(x), {\rm M}_d(\ell), {\rm N}_d(x), {\rm N}_d(\ell), {\rm N}_d(x;\epsilon)$$
the lattices 
${\rm M}_d({\rm I}_n; x), {\rm M}_d({\rm I}_n;\ell), {\rm N}_d({\rm I}_n; x), {\rm N}_d({\rm I}_n;\ell), {\rm N}_d({\rm I}_n; x;\epsilon)$. For $d$ even, we also denoted ${\rm N}_d(x)^{\pm}$ the lattices ${\rm N}_d(x; \epsilon)$ in the introduction.
A characteristic vector of ${\rm I}_n$ is $1^n$, so Proposition \ref{parnei} reads:

\begin{cor} The lattice ${\rm N}_d(x;\epsilon)$ is even if and only if $x_i$ is odd for each $i$ and we have
$\sum_i x_i^2 \equiv d(2+d)\epsilon \bmod 4d$ {\rm (}which forces $n \equiv 0 \bmod 8${\rm )}.
\end{cor}

%d=2k, n=2k(2+t2k)=4k(1+k)t mod 8 

% corollary N_d(x)^\epsilon is even iff d is even, x_i \equiv 1 \bmod 2, and sum_i x^i^2 \equiv (1+d/2)\epsilon \bmod 2.

%isometry group and orbits 

The isometry group ${\rm O}({\rm I}_n)$ is unusually large:
this is the group $\{ \pm 1\}^n \rtimes {\rm S}_n$  acting on $\Z^n$ by all possible permutations and sign changes of coordinates. The ${\rm O}({\rm I}_n)$-orbits on ${\rm C}_n(\Z/d)$, which are of great interest by Proposition \ref{propoequ},
will thus be in manageable quantity for small $d$ and $n$.
An element $x \in \Z^n$ will be called $d$-{\it ordered} if it satisfies \eqref{mainineq}, {\it i.e.}
$0 \leq x_1 \leq x_2 \leq \dots \leq x_d \leq d/2$.

\begin{fact} For any ${\rm O}({\rm I}_n)$-orbit $\Omega  \subset (\Z/d\Z)^n$, there is a unique $d$-ordered element $x \in \Z^n$ with $x \bmod d\, \, \in \,\Omega$.
\end{fact}

This obvious fact
%, together with the identity $g({\rm M}_d(x))={\rm M}_d(g(x))$ for  $g \in {\rm O}({\rm I}_n)$, 
explains why we always choose our $d$-isotropic elements $x$ to be $d$-ordered in our lists. 
%For a reasons to be discussed in Sect. \ref{sect:visible}
%we will usually assume furthermore $x_1=1$, and call $x$ {\it normalized} in this case.
Note however that two distinct $d$-isotropic and $d$-ordered elements of $\Z^n$ 
may generate the same line in $(\Z/d)^n$, hence give birth to the same $d$-neighbors
(see~\cite{allche} for a discussion about this).

\begin{remark} \label{warningdeven} {\rm Assume $x \in \Z^n$ is $d$-isotropic, $d$ is even and $g \in {\rm O}({\rm I}_n)$. We have $g({\rm N}_d(x;\epsilon))={\rm N}_d(g(x);\epsilon)$ by \eqref{actOLNd}. Beware however that, by assertion (ii) of Remark \ref{remdefline}, if we choose some $i$ and define $x' \in \Z^n$ by $x'_j=x_j$ for $j \neq i$, and $x'_i=x_i \pm d$, then $x'$ generates obviously the same line as $x$ in $(\Z/d)^n$, but we have ${\rm N}_d(x'; \epsilon)={\rm N}_d(x; \epsilon +1)$.}
\end{remark}

\begin{cor}\label{corxinull} Assume $x \in \Z^n$ is $d$-isotropic, with $d$ even and $x_i \equiv d/2 \bmod d$ for some $i \in \{1,\dots,n\}$. Then we have ${\rm N}_d(x; 0) \simeq {\rm N}_d(x; 1)$.
\end{cor}

\begin{pf} Set $x'=x - 2x_i \varepsilon_i$. Using Formula \eqref{actOLNd} for $g \in {\rm O}({\rm I}_n)$ defined by $g(\varepsilon_i)=-\varepsilon_i$ and $g(\varepsilon_j)=\varepsilon_j$ for $j \neq i$, we deduce ${\rm N}_d(x; 0) \simeq {\rm N}_d(x';0)$. By Remark~\ref{warningdeven} 
and $2 x_i \equiv d \bmod 2d$, we also have ${\rm N}_d(x';0)={\rm N}_d(x; 1)$.
\end{pf}

%wma x d-ordered with xn=d/2 and we have to show the 2 neighb are isom.
% xn -> d-xn

\section{Some invariants of lattices}\label{sect:inv}	
\label{sec:euclideansets}
\subsection{Configuration of vectors of given norm} 

Define an {\it Euclidean set} as a set $X$ equipped with 
an injection $X \overset{j}{\hookrightarrow} V$ into some Euclidean space $V$. 
We then denote by $\R X$ the Euclidean subspace of $V$ generated by $j(X)$.
Euclidean sets form a category if we define a morphism $(X,j) \rightarrow (X',j')$ as 
a map $f: X \rightarrow X'$ induced by a linear isometric embedding
$\widetilde{f} : \R X \rightarrow \R X'$, {\it i.e.} verifying\footnote{Note that given $f$, such an $\widetilde{f}$ is unique if it exists.} $\widetilde{f} \circ j = j' \circ f$.
 Note that it make sense to talk about the {\it scalar product} $x.y$ of two elements $x,y$ of an Euclidean set $(X,j)$
(namely $x.y=j(x).j(y)$), about the {\it rank} of $X$ ( {\it i.e.} $\dim {\rm Vect}_\R(j(X))$), about the {\it lattice generated by $X$} ({\it i.e.} $\Z \,j(X)$) ...\ps
	Let $L$ be an integral lattice in $V$. The {\it configuration of vectors of norm $i$} of $L$, already introduced in~\eqref{eq:fms},  is the Euclidean set 
\begin{equation}\label{defri} {\rm R}_i(L) = \{ v \in L\, \, |\, \, v.v=i\} \end{equation}	
with understood embedding ${\rm R}_i(L) \subset V$. Its isomorphism class is an invariant of the isometry class of $L$.  
Recall the notation ${\rm r}_i(L):=|{\rm R}_i(L)|$.
A natural variant of ${\rm R}_i(L)$ is the Euclidean set ${\rm R}_{\leq i}(L)$ 
defined by replacing $v.v=i$ with $v.v \leq i$ in \eqref{defri}.  
It is obvious that for two integral Euclidean lattices $L$ and $L'$, 
if we choose $i$ big enough so that ${\rm R}_{\leq i}(L)$ and ${\rm R}_{\leq i}(L')$ 
generate $L$ and $L'$, then $L$ is isometric to $L'$ if and only if 
${\rm R}_{\leq i}(L)$ isomorphic to ${\rm R}_{\leq i}(L')$. \ps

 For $i=1$, the isomorphism class of ${\rm R}_1(L)$ is obviously nothing more than the even integer ${\rm r}_1(L)$, since we have $w.v=0$ for $v \neq \pm w$ and $w,v \in {\rm R}_1(L)$.  For $i=2$ this is the same as the ${\bf ADE}$ root system of $L$, an important invariant that we review in \S \ref{subsect:rootsystems} below. 
 We are not aware of any general study or classification for the possible isomorphism classes of ${\rm R}_i(L)$ for $i \geq 3$. This question for $i=3$ is of great importance here, as experiments show that the unimodular lattices of dimension $n$ in our range are almost always generated over $\Z$ by their ${\rm R}_{\leq 3}$ : see~Table~\ref{tab:indexR3} for the {\it a posteriori} statistics, which explains much of Corollary~\ref{cori:invnorm3}.\ps
 
%As already explained in the introduction 
%(see Corollary~\ref{cori:invnorm3}),
%it will follow from the final classifications (but still has to be explained!) that unimodular lattices of rank up to $28$ are determined by their ${\rm R}_{\leq 3}$ {\color{blue}(even though they do not always span the lattices)}.

\tabcolsep=6pt
\begin{table}[H]
{\scriptsize 
\renewcommand{\arraystretch}{1.8} \medskip
\begin{center}
\begin{tabular}{c|c|c|c|c|c}
 $n \, \backslash \, d $ & $1$ & $2$ & $3$ & $4$ & $ \geq 5$  \\
\hline
 $26$ & $1857$ & $38$ & $2$  & $4$ & $0$  \\
 $27$ & $14425$ & $64$ & $1$ & $3$ & $0$  \\
\end{tabular} 
\end{center}
\caption{The number of isometry classes of rank $n$ unimodular lattices $L$ such that ${\rm R}_{\leq 3}(L)$ generates a sublattice of index $d$ in $L$ (including $d=\infty$).}
\label{tab:indexR3}
}
\end{table}

% L28=readvec("../data/unimodular_leq_28/nei_In/L28.gp");
% v = parapply(x->qfminim(nei(x),3)[1],L28);
% info_mimaa(v)=return( [vecmin(v)/2,vecmax(v)/2,vecsum(v)/(2*#v)*1.0]);
% 20 is for [94, [1; 3; 5; 7; 9; 11; 13; 15; 17; 19; 21; 23; 25; 27; 29; 31; 33; 35; 37; 39; 41; 43; 45; 47; 47; 47; 47; 47], 1, 1/162452817838080000, "D5", 1/42305421312000]
% gen3(v)={my(S,M,s); S=qfminim(nei(v),3)[3]; M=mathnf(S); s=matsize(M); if(s[1]<>s[2],return(0),return(abs(matdet(M))))};
% dim 26, [1 1857] [2   38] [3    2] [4    4]
% dim 27, [1 14425] [2    64] [3     1] [4     3]

\begin{remark} \label{finckepohst} {\rm (Computation of the set ${\rm R}_{\leq i}(L)$) If $L$ is an integral Euclidean lattice (given by a Gram matrix $\texttt{G}$), we use the Fincke-Pohst algorithm \cite{FP} to compute the sets ${\rm R}_{\leq i}(L)$ (function $\texttt{qfminim}(\texttt{G},i)$ in $\texttt{PARI/GP}$). As an indication, the average CPU time in ${\rm ms}$ on our machine to compute ${\rm R}_{\rm \leq i}(L)$, when $L$ is our list of $17059$ odd unimodular lattices of rank $27$, is about \texttt{2.3 ms} for $i=1$, \texttt{2.8 ms} for $i=2$ and \texttt{20 ms} for $i=3$. By comparison, computing a Gram matrix for such a lattice in neighbor form takes about \texttt{0.4 ms}. }
\end{remark}

\subsection{{\bf ADE} root systems and root lattices} \label{subsect:rootsystems}
${}^{}$
\indent By a {\it root} in an Euclidean space $V$, we mean an element $\alpha$ in $V$ with $\alpha.\alpha=2$.
The orthogonal symmetry about a root $\alpha$ is given by ${\rm s}_\alpha(v) = v - (v.\alpha) \alpha$.
An ${\bf ADE}$ {\it root system} is a finite Euclidean set $R$ consisting of 
roots such that for 
all $\alpha,\beta$ in $R$ we have $\alpha.\beta \in \Z$ and ${\rm s}_\alpha(R)=R$.
In other words, $R$ is a root system in $\R R$ in the sense of \cite[Ch. VI]{bourbaki}
satisfying $\alpha^\vee=\alpha$ for all $\alpha \in R$. 
Any such $R$ generates an even Euclidean lattice,
called the {\it associated root lattice}, and denoted ${\rm Q}(R)$ following Bourbaki.
We also set ${\rm res}\, R \,= \, {\rm res}\, {\rm Q}(R)$.
A {\it morphism} $R \rightarrow R'$ of root systems, also called an {\it embedding}, is a morphism
of Euclidean sets, or equivalently, a linear isometric embedding ${\rm Q}(R) \rightarrow {\rm Q}(R')$.
We talk about {\it sub root systems} for embeddings defined by an inclusion. \ps

For any integral Euclidean lattice $L$, then $R:={\rm R}_2(L)$ trivially is an ${\bf ADE}$ root system,
called the {\it root system of $L$}. We say that $L$ is a {\it root lattice} if we have ${\rm Q}(R)=L$.
The (non obvious but true) general equality $R = {\rm R}_2({\rm Q}(R))$ shows that the functors $R \mapsto {\rm Q}(R)$ and $L \mapsto {\rm R}_2(L)$ define inverse equivalences 
between the category of ${\bf ADE}$ root systems and that of root lattices (for linear isometries).
We trivially have ${\rm R}_2(L_1 \perp L_2) = {\rm R}_2(L_1) \coprod {\rm R}_2(L_2)$, where $\coprod$ denotes {\it orthogonal  disjoint union} of Euclidean sets, and ${\rm Q}(R_1 \coprod R_2) = {\rm Q}(R_1) \perp {\rm Q}(R_2)$. In particular, {\it irreducible} root systems correspond to {\it indecomposable} root lattices. \ps

Recall the standard lattices ${\rm A}_n$ ($n\geq 0$), ${\rm D}_n$ ($n\geq 0$) 
and ${\rm E}_n$ ($6 \leq n \leq 8$) form Sect. \ref{notations} (v).
All but ${\rm D}_1$ are root lattices, and all are indecomposable for $n\geq 1$ but 
${\rm D}_2$. We denote respectively by ${\bf A}_n$, ${\bf D}_n$ and ${\bf E}_n$ their root systems.  
By the ${\bf ADE}$ classification, any irreducible root system is isomorphic such a root system,
and the unique coincidences between them are ${\bf A}_0 = {\bf D}_0 = {\bf D}_1= \emptyset$,  
${\bf D}_2  \simeq {\bf A}_1 \coprod {\bf A}_1$ and ${\bf D}_3 \simeq {\bf A}_3$. \ps

\begin{remark} \label{computers} \label{rem:rscalc}
{\rm (An algorithm for computing root systems) 
If $L$ is an integral lattice in the Euclidean space $V$, the structure of its 
root system ${\rm R}_2(L)$ may be efficiently computed as follows. 
Determine first the set $R={\rm R}_2(L)$ as in Remark \ref{finckepohst}, 
choose a linear form 
$\varphi$ on $V$ with $0 \not \in\varphi(R)$, and set $R^+ = \{ \alpha \in R, \varphi(\alpha)>0\}$ ({\it positive roots}). 
(Actually \texttt{PARI}'s \texttt{qfminim}(\texttt{G},2) function directly returns such an $R^+$ rather than $R$).
Compute then the {\it Weyl vector} $\rho = \frac{1}{2}\sum_{\alpha \in R^+} \alpha$ 
and the {\it basis} $B = \{ \alpha \in R^+\, |\, \, \rho \cdot \alpha=1\}$ of $R$ 
associated to $R^+$. Compute the scalar products $b.b'$ for $b,b' \in B$ and
view $B$ as the vertices of the undirected graph 
with an edge between $b$ and $b'$ if and only if $b\neq b'$ and $b.b' \neq 0$ ({\it Dynkin diagram of $R$}).
The connected ({\it irreducible}) components of this union of trees are easily computed recursively, 
and identified as of type ${\bf A}_n$, ${\bf D}_n$ or ${\bf E}_n$ by simply looking at their unique vertice $x$ with valence $>2$ (if exists), and the sum of the valences of the $3$ neighbors of $x$. 
See \cite{cheweb} for our concrete implementation. As an indication, the average CPU time for computing the isomorphism class of the root system of a unimodular lattice of rank $27$ in our list is \texttt{5.9 ms}. 
As $|B| \leq \dim V$ is very small in practice, naive graphs algorithms are perfectly suitable for the last part above: 
$99.8$ \% of the CPU time is used for the computation of the sets $R^+$ and $B$.
}
\end{remark}

Let $L$ be an integral Euclidean lattice with {\rm ADE} root system $R:={\rm R}_2(L)$. 
For each $\alpha \in R$ the orthogonal reflexion ${\rm s}_\alpha(x) = x - (x.\alpha) \alpha$ 
lies in ${\rm O}(L)$. The {\it Weyl group} of $L$
is the subgroup ${\rm W}(L)$ of ${\rm O}(L)$ generated by those ${\rm s}_\alpha$ with $\alpha \in R$.
This is a normal subgroup of ${\rm O}(L)$ isomorphic to ${\rm W}(R):={\rm W}({\rm Q}(R))$. 
Moreover, if we choose a positive root system $R^+$ of $R$, and denote by $\rho$ the associated Weyl vector as in Remark~\ref{rem:rscalc}, we have another subgroup 
\begin{equation}
\label{def:OLrho} 
{\rm O}(L; \rho) := \{ \gamma \in {\rm O}(L)\, \, |\, \, \gamma(\rho)=\rho\}.
\end{equation}
As is well-known, ${\rm W}(L)$ acts simply transitively on the set of Weyl vectors of $L$ (and on the set of positive root systems). We have thus 
\begin{equation}
\label{eq:decompOLOLrho}
{\rm O}(L) = {\rm W}(L) \rtimes {\rm O}(L; \rho).
\end{equation}
Also, for $w \in {\rm W}(L)$ we obviously have ${\rm O}(L; w(\rho))=w {\rm O}(L; \rho)w^{-1}$.
It follows that the ${\rm W}(L)$-conjugacy class of the subgroup ${\rm O}(L; \rho)$ in ${\rm O}(L)$ is canonical, and does not depend on the choice of $\rho$. Moreover, each ${\rm O}(L; \rho)$ is naturally isomorphic to ${\rm O}(L)/{\rm W}(R)$.

\begin{definition} \label{def:redisogp}
Let $L$ be integral Euclidean lattice $L$ with root system $R$. 
Define the {\rm reduced isometry group} of $L$ as the group ${\rm O}(L)^{\rm red}={\rm O}(L)/{\rm W}(R)$,
and the {\rm reduced mass} of $L$ by the formula ${\rm rmass}(L):=\frac{1}{|{\rm O}(L)^{\rm red}|}$.
%We have ${\rm rmass}(L) \,=\, |{\rm W}(R)|\, {\rm mass}(L)$.
\end{definition}

%{\color{red} We also denote by ${\rm W}(L)^\pm$ the subgroup of ${\rm O}(L)$ generated by $\pm {\rm id}_L$ and ${\rm W}(L)$. }

A first important application of this notion is the following remark.
As the structure of ${\rm W}(R)$ is well-known, 
it follows that in order to determine generators of ${\rm O}(L)$, 
or simply its order, it is enough to do so for ${\rm O}(L;\rho)$.

\begin{remark} \label{pleskensouvignier} {\rm 
({\it Computation of the reduced isometry group})}
{\rm Let $L$ be an integral Euclidean lattice.
Choose a Weyl vector $\rho$ of $L$ as in Remark~\ref{computers}.
As was already observed and used in~\cite{chniemeier}, 
it turns out that the Plesken-Souvignier algorithm
\cite{pleskensouvignier}  directly allows to compute generators and the order of ${\rm O}(L;\rho)$.
Indeed, is enough to apply it to a pair consisting of a Gram matrix of $L$ 
and of a Gram matrix of the bilinear form $(x,y) \mapsto 4(\rho.x)(\rho.y)$, in a same basis. 
This actually returns order and generators for $\pm {\rm O}(L;\rho)$, but the similar information for 
${\rm O}(L;\rho)$ easily follows.\footnote{Souvignier's code is availble in \texttt{PARI/GP} as $\texttt{qfauto}(\texttt{G})$ with $\texttt{G}$ a Gram matrix of $L$. For all the computations of ${\rm O}(L; \rho)$ in this paper, we use the the LLL-algorithm (\texttt{PARI}'s $\texttt{qflllgram}$) to find suitable a $\texttt{G}$. See~\cite{allche} for more clever choices of bases.} This computation of ${\rm O}(L;\rho)$ is usually much faster than that of ${\rm O}(L)$
(the bigger ${\rm R}_2(L)$ is, the faster). As an indication, the average CPU time for the computation of ${\rm O}(L)$ for our $14493$ rank $27$ unimodular lattices with no norm $1$ vector using this method (and \texttt{PARI}'s $\texttt{qfauto}$) is  \texttt{8.9 s}.
}
\end{remark}

We end this section with a few more definitions, for a later use. \ps
%
%, this is a monoid for $\coprod$. \ps\ps
\ps
-- An embedding $f : R' \hookrightarrow R$ is called {\it saturated} if the subgroup $f({\rm Q}(R'))$ is saturated in ${\rm Q}(R)$ (see Sect. \ref{notations} (vi)). The saturated sub root systems of an ${\bf ADE}$ root system $R \subset V$ are those obtained by intersecting $R$ with a subspace of $V$; they are sometimes called {\it parabolic} and their Dynkin diagram are obtained from that of $R$ by removing a finite set of vertices. \ps
--  We denote by ${\bf RS}$ the set of isomorphism classes of ${\bf ADE}$ root systems. We often write 
$n_1 R_1\, n_2 R_2\,\dots n_k R_k$ for the orthogonal disjoint union of $n_i$ copies of $R_i$, for $i=1,\dots,k$.

\subsection{Vectors of norm $\leq 3$}\label{subsect:norm3}

${}^{}$
We start with some information on the number of vectors of norm $3$ of the lattices we are interested in.
The statement of the following proposition uses the notation ${\rm Exc}(L)$ for $L \in \mathcal{L}_n$, that will only be introduced in \S \ref{subsect:exclat} when discussing exceptional lattices. 

%Most lattices are {\it non exceptional}, {\it i.e.} satisfy $|{\rm Exc}(L)|=0$, and the possibilities for $|{\rm Exc}(L)|$ are given by~Proposition~\ref{prop:sizeExc}.

\begin{prop} 
\label{prop:thetaseries}
For all $L \in \mathcal{L}_n$ with $24 \leq n \leq 28$ and ${\rm r}_1(L)=0$ we have
${\rm r}_3(L)\,= \,\ \frac{4}{3}n (n^2 \,- \,69\,n \,+ 1208)\, + \,2(n-24)\, {\rm r}_2(L) \,- \,2^{36-n} \,|{\rm Exc}(L)|$.
\end{prop}
%${\rm r}_3(L)\,= \,\ \frac{1}{3} (4\,n^3 \,- \,276\,n^2 \,+ 4832\,n)\, + \,(2n-48)\, {\rm r}_2(L) \,- \,2^{36-n} \,|{\rm Exc}(L)|$.

\begin{comment}
%% GP CODE
TT2=2^4*q*(1+q^2+q^6+O(q^7))^4;
T3= 1 + 2*q + 2*q^4 + 2*q^9 + O(q^16);
T4= 1 - 2*q + 2*q^4 - 2*q^9 + O(q^16);
D=(1/16)*TT2*T4^4;
Dp=-(1/16)*TT2*T3^4;
Dpm=-(1/16)*T4^4*T3^4;
TL=T3^28+a*T3^20*D+b*T3^12*D^2+c*T3^4*D^3;
TLp=T4^28+a*T4^20*Dp+b*T4^12*Dp^2+c*T4^4*Dp^3;
TLpm=TT2^7+a*TT2^5*Dpm+b*TT2^3*Dpm^2+c*TT2*Dpm^3;
%%TMd(z) = (i/z)^14 2 TM(-1/z)
%%2TM(-1/z) = TL(-1/z)+TL(-1/z +1)
TMd= TL + TLpm;
%%Theta_M^# = 1 + (a + (-1/256*c + 56))*q + ...
%%Theta_L = 1 + (a + 56)*q + (32*a + (b + 1512))*q^2 + (468*a + (8*b + (c + 26208)))*q^3 + ...
r1=a+56;
r2=32*a+b+1512 = 32 r1 + b - 280;
r3=468*a+8*b+c+26208;
a=r1-56;
b=r2-32*a-1512=-32*r1 + (r2 + 280);
c=r3-(468*a+8*b+26208);
e=a + (-1/256*c + 56);
a=r1-56;
b=-32*r1 + (r2 + 280);
c=r3 + (-212*r1 + (-8*r2 - 2240));
%% e*256 = -r3 + 468*r1 + 8*r2 + 2240
r1=2*m;
r2=s2+2*m*(m-1);
r3=s3+2*m*s2+8*m*(m-1)*(m-2)/6;
s3=-8*m*(m-1)*(m-2)/6-2*m*s2+468*2*m+8*(s2+2*m*(m-1))+2240-256*2*m-256*E;
%-8*m*(m-1)*(m-2)/6+16*m^2 + (-2*s2 + 408)*m + (8*s2 + (-256*E + 2240))
%s3=-4/3*m^3 + 20*m^2 + (-2*s2 + 1216/3)*m + (8*s2 + (-256*E + 2240));
n=28-m;
% s3=4/3*n^3 - 92*n^2 + 4832/3*n + (8-2*(28-n))*s2 - 256*E;
\end{comment}

\begin{pf} Assume first we have $L \in \mathcal{L}_{28}$ (and possibly norm $1$ vectors) and let
$M$ denote the even part of $L$. A simple computation of the coefficient in $q$ of the theta series of $M^\sharp$,
following the arguments in \cite{bachervenkov} \S 4 (or \cite[Ch. 4 \S 4]{conwaysloane}), shows the relation
$ \,256 \,{\rm r}_1(M^\sharp) = - {\rm r}_3(L)\, +\, 8 \,{\rm r}_2(L) + \,468 \,{\rm r}_1(L) \,+ \,2240$.
Since $28 \equiv 4 \bmod 8$ we also have ${\rm r}_1(M^\sharp)={\rm r}_1(L) + |{\rm Exc}(L)|$ by Formula~\eqref{eq:bijexc4mod8}. Choose now $L_0 \in \mathcal{L}_n$ with $n\leq 28$ and ${\rm r}_1(L)=0$. 
We apply above relation to $L=L_0 \perp {\rm I}_{m}$ with $m=28-n$. We conclude the equality of the statement (for $L_0$) by the equalities ${\rm r}_1(L)=2m$, ${\rm r}_2(L)={\rm r}_2(L_0)+ 2^2 \binom{m}{2}$, ${\rm r}_3(L)\,=\,{\rm r}_3(L_0)\,+\, 2m \,{\rm r}_2(L_0) \,+\,2^3 \binom{m}{3}$, and $|{\rm Exc}(L)|=2^m|{\rm Exc}(L_0)|$ for $n\geq 24$. 
\end{pf}

\begin{example} 
\label{ex:r3dim2627}
{\rm 
We deduce ${\rm r}_3(L)\,=\,3120\,+\,4\,{\rm r}_2(L)\,-\,1024\, |{\rm Exc}(L)|$ in the case $n=26$, 
and ${\rm r}_3(L)\,=\,2664\,+\,6\,{\rm r}_2(L)\,-\, 512\, |{\rm Exc}(L)|$ for $n=27$.
In both cases, we have $|{\rm Exc}(L)|=2$ if $L$ is exceptional, 
and $|{\rm Exc}(L)|=0$ otherwise, by Proposition~\ref{prop:sizeExc}.
}
\end{example}

Our aim now is to discuss a few invariants of ${\rm R}_{\leq 3}$ that we have 
used during our simultaneous proofs of Theorem \ref{X2627} 
and Corollary~\ref{cori:invnorm3}. We sincerely apologize that what follows is
mostly empirical. We mostly relates facts 
(namely Propositions~\ref{prop:invdim26} and \ref{prop:invdim27}) that we 
observed during our search and only proved by case by case computations. 
In each case, it is an open problem to find conceptual explanations for our computations.\ps
\ps

A {\it component of size $s$} of a root system $R$ 
is a union of $s$ distinct irreducible components of $R$. 
%Denote by ${\bf RS}$ the set of isomorphism classes of ${\bf ADE}$ root systems,
We denote by $\underline{R} \in {\bf RS}$ the isomorphism class of $R$. 

\begin{definition} 
\label{def:invdeltas}
For $L$ an integral Euclidean lattice and $s\geq 0$ an integer, 
we denote by $\delta_s(L) \in \Z[{\rm RS} \times \N]$ 
the sum of $(\underline{C},m)$ of ${\rm RS} \times \Z_{\geq 0}$,
where $C$ runs among 
the components of size $s$ of ${\rm R}_2(L)$ and $m=|{\rm R}_3(C^\perp \cap L)|$.
\end{definition}

The invariant $\delta_0(L)$ is just ${\rm r}_3(L)$, a weak information by Example~\ref{ex:r3dim2627}.
The invariants $\delta_k(L)$ have already been used for instance by 
Megarban\'e in his study of the rank $26$ even lattices of determinant $3$ \cite{megarbane} 
(see Example~\ref{ex:G26G27}). As an example, let us consider again the $7$ isometry classes of unimodular lattices of rank $26$ with root system $10 {\bf A}_1$, listed in Table~\ref{tab:10A1dim26}.  The invariants $\delta_1$ and $\delta_2$ of these lattices (ordered as in that table) are given by Table~\ref{tab:10A1dim26inv}:

\tabcolsep=2pt
\begin{table}[H]
{\scriptsize 
\renewcommand{\arraystretch}{1.8} \medskip
\begin{center}
\begin{tabular}{c|c}
 $\delta_1$ & $\delta_2$ \\
\hline
 $10\,\, ({\bf A}_1,2578)$ & $4\,\,(2{\bf A}_1, 1968) \,+ \, 6\, (2{\bf A}_1, 2000) \,+\, 20 \,(2{\bf A}, 2032)\}+\,15\,(2{\bf A}_1,2064)\,$ \\
 $10\,\, ({\bf A}_1,2578)$ & $13 \,\,(2{\bf A}_1, 2000) \,+ \, 16\, (2{\bf A}_1, 2032) \,+\, 16 \,(2{\bf A}, 2064)\}$ \\
 $10\,\, ({\bf A}_1,2578)$ & $11 \,\,(2{\bf A}_1, 2000) \,+ \, 24\, (2{\bf A}_1, 2032) \,+\, 10 \,(2{\bf A}, 2064)\}$ \\
 $10\,\, ({\bf A}_1,2578)$ & $20 \,\,(2{\bf A}_1, 2000) \,+ \, 25\, (2{\bf A}_1, 2064) $ \\
 $10\,\, ({\bf A}_1,2578)$ & $4 \,\,(2{\bf A}_1, 1936) \,+ \, 16\, (2{\bf A}_1, 2000) \,+\, 25 \,(2{\bf A}, 2064)\}$ \\
 $10\,\, ({\bf A}_1,2578)$ & $20 \,\,(2{\bf A}_1, 1936) \,+ \, 25\, (2{\bf A}_1, 2064) $ \\
 $({\bf A}_1, 18) \,+\, 9\,({\bf A}_1,1042) $ & $9 \,\,(2{\bf A}_1, 16) \,+ \, 36\, (2{\bf A}_1, 912) $ 
\end{tabular} 
\end{center}
\caption{The invariants $\delta_1$ and $\delta_2$ of the $7$ lattices of Table~\ref{tab:10A1dim26}.}
\label{tab:10A1dim26inv}
}
\end{table}

In particular, all those lattices are distinguished by $\delta_2$ (but not by $\delta_1$).
It follows from our computations that this is a general fact in rank $26$.

\begin{prop} 
\label{prop:invdim26}
Two unimodular lattices of rank $26$ are isometric if, and only if, they have the same root system and the same invariants $\delta_1$ and $\delta_2$.
\end{prop}

The invariants $\delta_s$ are however not strong enough in dimension $27$.
%In the general, they are typically weak for lattices with empty (or small) ${\rm R}_2$. 
We now discuss a second invariant of ${\rm R}_{\leq 3}$. 

\begin{definition}
\label{def:graphGL}
For an Euclidean integral lattice $L$, we define ${\rm G}(L)$ as the undirected graph with vertices the nonzero pairs $\{\pm x\}$
with $x \in {\rm R}_{\leq 3}(L)$, and with $|x.y|$ arrows between $\{\pm x\}$ and $\{ \pm y\}$. 
\end{definition}

Of course, we have $|x.y| \leq 3$ for all $x,y \in {\rm R}_{\leq 3}(L)$.  The isomorphism class of this graph ${\rm G}(L)$ only depends on that of ${\rm R}_{\leq 3}(L)$. 
In Table~\ref{tab:sizeR3} below, we give the {\it a posteriori} information for the number of vertices of ${\rm G}(L)$ 
for $L$ in ${\mathcal L}_n$ with ${\rm r}_1(L)=0$ for $n=26,27$ (agreeing Example~\ref{ex:r3dim2627}).

\tabcolsep=6pt
\begin{table}[H]
{\scriptsize 
\renewcommand{\arraystretch}{1.8} \medskip
\begin{center}
\begin{tabular}{c|c|c|c}
 $n$ & ${\rm min}$ & ${\rm max}$ & ${\rm average}$ \\
\hline
 $26$ & $556$ & $2850$ & $\simeq 1776$ \\
 $27$ & $820$ & $3277$ & $\simeq 1573$ \\
 % $28$ & $20$ & $3388$ & $\simeq 1318$
\end{tabular} 
\end{center}
\caption{The minimum, maximum, and average, number of vertices of the graph ${\rm G}(L)$ for $L \in \mathcal{L}'_n$.}
\label{tab:sizeR3}
}
\end{table}

Any of invariant of graphs can be applied to study ${\rm G}(L)$. 
From a computational point of view, an especially simple one
 that we can consider is the rank ${\rm h}_p(G)$, 
of the adjacency matrix ${\rm mod}\, p$ of a graph $G$. 
Here $p$ is any given prime. 
So for $s\geq 0$ and a prime $p$, we can define 
a variant $$\delta_{s,p}(L) \in \Z[{\rm RS} \times \N \times \N]$$ of $\delta_s$ 
by replacing each $(\underline{C},m)$ in Definition~\ref{def:invdeltas} 
by $(\underline{C},m,r)$ with $r={\rm h}_p(\mathrm{G}(C^\perp \cap L))$. 
It turns out that these $\delta_{s,p}$ do suffice to distinguish 
all unimodular lattices in rank $27$. More precisely, our computations show:

\begin{prop} 
\label{prop:invdim27}
Two unimodular lattices of rank $27$ are isometric if, and only if, 
they have the same root system 
and the same invariants $\delta_{s,p}$ with $s\leq 3$ and $p=5$, 
unless their root system is in the following list : 
{\small
$$3{\bf A}_1, 6{\bf A}_1, 7{\bf A}_1, 3{\bf A}_1{\bf A}_2, 5{\bf A}_1{\bf A}_2, 7{\bf A}_1{\bf A}_2,
4{\bf A}_12{\bf A}_2, 6{\bf A}_12{\bf A}_2, 8{\bf A}_12{\bf A}_2, 5{\bf A}_13{\bf A}_2.$$
}
For such root systems, we need the invariants $\delta_{s,p}$ for $s\leq 7$ and $p=5,7$.
\end{prop}

We will not insist much on these invariants $\delta_{s,p}$ here, 
but rather refer to the companion work \cite{allche} in 
which we will define a another invariant of ${\rm G}(L)$ 
(inspired by the work of Bacher and Venkov \cite{bachervenkov}) 
which will turn out to be fine enough to distinguish all unimodular lattices of rank $\leq 28$.

\section{The visible root system of a $d$-neighbor of ${\rm I}_n$}
\label{sect:visible}
\label{sec:visiblers}

\subsection{The root system}
\label{subsect:visrs} Let $x \in \Z^n$ be $d$-isotropic, and let $N={\rm N}_d(x')$ be an associated $d$-neighbor 
of ${\rm I}_n$. Recall ${\rm I}_n \cap N = {\rm M}_d(x)$. In this section, we study the visible root system 
of $N$ in the sense of Definition~\ref{def:visrs}, namely the subroot system
${\rm R}_2({\rm M}_d(x))$ of ${\rm R}_2(N)$. The term {\it visible} reflects the fact that this root system is immediately seen on the shape of $x$. The receipe is as follows.\ps

To any $x \in \Z^n$ and $d\geq 1$ we have an equivalence relation $\sim$ on $\{1,\dots,n\}$ defined by
\begin{equation} \label{eqrelx} i \sim j \Leftrightarrow x_i \equiv \pm x_j \mod d .\end{equation}
There are two {\it distinguished} subsets $D$ and $D'$, defined respectively as the subset of $i \in \{1,\dots,n\}$ with $x_i \equiv d/2 \bmod d$ or $x_i \equiv 0 \bmod d$ (so $D=\emptyset$ for $d$ odd). We set ${\rm m}(x)=|D|$, ${\rm m}'(x)=|D'|$ and we denote by ${\rm a}(x)$ 
 the integer partition ${\rm a}_1(x) \geq {\rm a}_2(x) \geq \dots$ of $n-{\rm m}(x)-{\rm m}'(x)$ defined by the sizes of the equivalence classes of $\sim$  different from $D$ and $D'$. All of $\sim$, ${\rm a}(x)$, the ${\rm a}_i(x)$, ${\rm m}(x)$ and ${\rm m}'(x)$ only depend on the line $\ell={\rm l}(x) \subset (\Z/d)^n$ generated by $x$, and we also use the notations $\sim_\ell$, ${\rm a}(\ell)$, ${\rm a}_i(\ell)$, ${\rm m}(\ell)$ and ${\rm m}'(\ell)$ for them.
 
\begin{prop}\label{visrs} For all $x$ in $\Z^n$ and $d\geq 1$ we have an isomorphism 
$${\rm R}_2({\rm M}_d(x)) \simeq {\bf D}_{{\rm m}(x)} \, {\bf D}_{{\rm m}'(x)} \, {\bf A}_{{\rm a}_1(x)-1} \, {\bf A}_{{\rm a}_2(x)-1} \, \dots.$$
\end{prop}
\begin{pf}  Recall ${\rm R}_2({\rm I}_n) = \{ \pm \varepsilon_i \pm \varepsilon_j \, \, |\,\, 1 \leq i < j \leq n\}$ from 
\eqref{eq:r1r2in}. 
By definition of ${\rm M}_d(x)$, 
the root $\varepsilon_i-\varepsilon_j$ (resp. $\varepsilon_i+\varepsilon_j$) of ${\rm I}_n$ 
belongs to ${\rm M}_d(x)$ if and only $x_i \equiv x_j \bmod d$ (resp. $x_i \equiv - x_j \bmod d$). 
Up to applying an element of ${\rm O}({\rm I}_n)$ if necessary, 
we may assume $x$ is $d$-ordered.
In this case, the equivalence classes of the relation $\sim$ clearly are intervals. 
Also, the visible roots are the $\varepsilon_i-\varepsilon_j$ whenever $x_i=x_j$,  
and for $d$ even the $\pm(\varepsilon_i+\varepsilon_j)$ whenever $x_i \equiv x_j \equiv 0 \bmod d/2$ (with $i\neq j$ in both cases):
we recognize the root system of the statement.
\end{pf}

\subsection{Saturation properties} 
\label{subsect:satprop}
As noticed in \S\ref{subsect:visiblepart}, we will usually assume $x_i \neq 0 \bmod d$ for all $i$, or equivalently ${\rm m}'(x)=0$, in order to have ${\rm r}_1({\rm M}_d(x))=0$.
Our aim until the end of this \S \ref{sect:visible} is to discuss the relations between 
the visible root system $R^{\rm v}$ and the actual root system of a $d$-neighbor $N$ of ${\rm I}_n$, 
as well as the constraints on the embedding $R^{\rm v} \hookrightarrow N$.
%A first important property of the visible root system is the following.

\begin{lemma}\label{satrv1} Let $x \in \Z^n$ and $d\geq 1$ with $x_i \not \equiv 0 \bmod d$ for each $i$.
Set $M={\rm M}_d(x)$ and $R^{\rm v}={\rm R}_2(M)$. Then ${\rm Q}(R^{\rm v})$ is saturated is $M$.
\end{lemma}

\begin{pf} Applying an element of ${\rm O}({\rm I}_n)$ if necessary, we may assume $x$ is $d$-ordered.
As observed in Proposition \ref{visrs}, the classes $I$ of the equivalence relation $\sim$ on $\{1,\dots,n\}$ 
associated to $x$ and $d$ are intervals in this case, and simply determine $R^{\rm v}$: 
each $I \neq D$ gives rise to an ${\bf A}_{|I|-1}$ component (note $x_i \not \equiv 0 \bmod d$ for $i \in I$), 
and $I=D$ to a ${\bf D}_{|I|}$ component.
For each class $I$, define an abelian group $\Delta_I$ as follows: 
set $\Delta_I=\Z$, unless $I=D$, $|D|\geq 2$ and $\Delta_I=\Z/2$.
Define also $\varphi_I : \Z^n \rightarrow \Delta_I$ by  $\varphi_I(v) \equiv  \sum_{i \in I} v_i$. 
We have a surjective linear map
%\begin{equation} \label{mapphi} 
$\varphi=\prod_I \varphi_I : \Z^n \longrightarrow \prod_{I} \Delta_I$, %\end{equation}
where $I$ runs among the equivalence classes of $\sim$.
We clearly have
$$\varphi(M) = \{ (w_I) \in \prod_{I} \Delta_I \, \, |\, \, \sum_I x_I w_I \equiv 0 \bmod d\},$$
where $x_I \in \Z/d$ denotes the common class of $x_i$ for $i \in I$.
(Note that if $D \neq \emptyset$ the product $x_D w_D$ is well-defined in $\Z/d$ as we have $d$ even and $x_D \equiv d/2$.)
The description of $R^{\rm v}$ recalled above shows ${\rm ker}\, \varphi \,= \,{\rm Q}(R^{\rm v})$, hence
$\varphi(M) \simeq M/{\rm Q}(R^{\rm v})$. The group $\prod_I \Delta_I$ is torsion-free, hence so is its subgroup
$\varphi(M)$, unless $d$ is even and $|D|\geq 2$. In this case, the unique torsion element $(w_I)$ is defined by $w_I\equiv 0$ for $I \neq D$, 
and $w_D \equiv 1$. But this element does not belong to $\varphi(M)$ as $x_D \equiv d/2 \neq 0 \bmod d$.
\end{pf}

\begin{definition}\label{dkerdef} Let $L$ be an integral Euclidean lattice, $R={\rm R}_2(L)$ and $R'$ a sub root system of $R$. We say that $R'$ is a $d$-{\rm kernel} of $L$ if there is a surjective linear map $\varphi : L \rightarrow \Z/d$ with  ${\rm R}_2({\rm ker}\, \varphi)=R'$. 
In the case $L={\rm Q}(R)$ we also say that $R'$ is a $d$-{\rm kernel} of $R$.
\end{definition}

In other words, the $d$-kernels of $L$ are the $R'={\rm R}_2(M)$ for $M$ a sublattice of $L$ with $L/M \simeq \Z/d$. In such a situation, and if we set $R={\rm R}_2(L)$, 
we also have ${\rm Im}\, ({\rm Q}(R) \rightarrow L/M) \simeq \Z/d'$ for some $d'\, |\,d$, hence $R'$ is a $d'$-kernel of $R$.

\begin{cor} 
\label{satrv2}
Assume $x \in \Z^n$ is $d$-isotropic with $x_i \not \equiv 0 \bmod d$ for each $i$.
Set $M={\rm M}_d(x)$, $R^{\rm v}={\rm R}_2(M)$ 
and let $N$ be a $d$-neighbor of ${\rm I}_n$ with line $l(x)$. 
Then $R^{\rm v}$ is a $d$-kernel of $N$. 
If furthermore $d$ is prime to ${\rm a}_i(x)$ for all $i$,
and odd in the case ${\rm m}(x)>1$, then ${\rm Q}(R^{\rm v})$ is saturated in $N$.
\end{cor}

\begin{pf} The first assertion is clear. 
As ${\rm Q}(R^{\rm v})$ is saturated in $M$ by Lemma~\ref{satrv1}, 
the saturation $S$ of ${\rm Q}(R^{\rm v})$ in $N$ satisfies $S \cap M = {\rm Q}(R^{\rm v})$.
The isotropic subgroup $S/{\rm Q}(R^{\rm v})$ of ${\rm res}\, R^{\rm v}$ embeds
thus in $\Z/d$. So its order $s$ satisfies
\footnote{We even have $s^2\,\,|\,\,|{\rm res}\, R^{\rm v}|$. 
Indeed, for any isotropic subgroup $I$ of a finite bilinear abelian group $A$ 
we have $|A|=|A/I^\perp|\cdot |I^\perp/I|\,\cdot  |I|$ and 
$A/I^\perp \simeq {\rm Hom}(I,\Q/\Z) \simeq I$.}
$s\,\, |\,\, d$  and $s\,\,|\,\,|{\rm res}\, R^{\rm v}|$.
We conclude as by Proposition~\ref{visrs} and Table~\ref{tab:venkov} 
we have $|{\rm res}\, R^{\rm v}|=f \prod_{i} {\rm a}_i(x)$ 
with $f=4$ (case ${\rm m}(x)>1$) or $f=1$ (otherwise). 
\end{pf}

\begin{remark} \label{remaconclpropsat} 
  {\it
%  Proposition \ref{satrv} shows that the visible root system $R^{\rm v}$ satisfies 
 % strong saturation constraints in the $d$-neighbor $N$.  
  As a consequence, if $d$ is a sufficiently big prime then ${\rm Q}(R^{\rm v})$ is saturated in $N$,
  which is a strong constraint. Theorem~\ref{thmi:stat2} shows that 
  generically this is actually the only constraint. 
  %Moreover, the analysis above explains why the cases where $d$ is not a prime allows 
  %more interesting lattice constructions. 
  %From this perspective, the general (non necessarily cyclic) neighbors of ${\rm I}_n$ 
  %would be even more permissive. 
  }
\end{remark}

 Of course, there are many examples of unimodular lattices $N$ with root system $R$ such that 
 if $S$ denotes the saturation of ${\rm Q}(R)$ in $L$ then $S/{\rm Q}(R)$ is not a cyclic group.
 We will thus often have $R^{\rm v} \subsetneq {\rm R}_2(N)$ in practice. 
 By the remark following Definition \ref{dkerdef}, and still in the notations of Corollary \ref{satrv2},
 the visible root system $R^{\rm v}$ is also a $d'$-kernel of ${\rm R}_2(N)$ for some $d'|d$.
 We postpone to \S \ref{subsect:dkerdesc} the description of all the $d$-kernels of a given root system.
 We content ourselves here with the following simple observation (see Remark \ref{classdkernelrem} for a more precise statement).

\begin{prop} 
Let $R'$ be a saturated sub root system of the {\bf ADE} root system $R$ with ${\rm rk}\, R' < {\rm rk}\, R$.
Then $R'$ is a $d$-kernel of $R$ for all $d$ big enough.
\end{prop}

\begin{pf} 
We have an abelian group decomposition ${\rm Q}(R) = {\rm Q}(R') \oplus P$ 
with $P \simeq \Z^r$ and $r={\rm rk} R - {\rm rk} R'>0$. 
As $R$ is finite, we may find a surjective linear map $\varphi : {\rm Q}(R) \rightarrow \Z$
with $\varphi(R')=0$ and $\varphi(r) \neq 0$ for $r \not\in R'$. 
We conclude by using $\varphi \otimes \Z/d$ with $d> \varphi(r)$ for all $r \in R$. 
\end{pf}

\subsection{An example : the {\it safe} case} 
\label{subsect:safe}
For the purpose of unimodular hunting, we are led to the following definition.

\begin{definition} 
\label{def:safe}
Let $R$ and $S$ be {\bf ADE} root systems. 
We say that $(R,S)$ is {\rm safe} if for any integral lattice $L$ with ${\rm R}_{\leq 2}(L)=R$, there is an isometric embedding ${\rm Q}(S) \rightarrow L$ whose image is saturated in $L$.
\end{definition}

By definitions, if $(R,S)$ is safe then $S$ is isometric to a saturated sub root system of $R$.
Moreover, $(R,S')$ is also safe 
for any saturated sub root system $S'$ of $S$. 
Note also that if $L$ is an integral lattice with root system $R$, 
and if $V$ is the Euclidean space generated by $R$, 
then $L \cap V$ is saturated in $L$, so that we may actually assume $L \subset V$ in the definition above. 
In this case, $L \subset {\rm Q}(R)^\sharp$ is uniquely determined by 
the isotropic subspace $I = L/{\rm Q}(R)$ of ${\rm res}\,R$, by Sect. \ref{notations} (iii). 
In terms of the Venkov map recalled {\it loc. cit.} (iv), the assumption ${\rm R}_{\leq 2}(L) \simeq R$ 
is equivalent to ${\rm qm}(x) \neq \frac{1}{2}, 1$ for all $x \in I$.
  
\begin{definition} 
\label{def:detecting}
We call an ${\bf ADE}$ root system $S$ {\rm detecting} if for all integral lattices 
${\rm Q}(S) \subsetneq L \subset {\rm Q}(S)^\sharp$ we have $S \subsetneq {\rm R}_{\leq 2}(L)$.
Equivalently, $S$ is detecting if for all $x$ in ${\rm res}\,S$ with ${\rm qm}(x) \in \frac{1}{2}\Z$ we have ${\rm qm}(x)\leq 1$.
\end{definition}

The interest of this notion for us is the following proposition.

\begin{prop} 
\label{prop:safecriterion}
Assume $R$ is the orthogonal disjoint union of its sub root systems $S$ and $T$. If $S$ is detecting then $(R,S)$ is safe.
\end{prop}

\begin{pf} Write $R= S \perp T$ and let $L$ be an integral lattice with ${\rm R}_{\leq 2}(L)=R$.
The saturation $L'$ of ${\rm Q}(S)$ in $L$ is orthogonal to $T$, so we must have ${\rm R}_2(L')=S$.
We also have ${\rm R}_1(L')\subset {\rm R}_1(L)=\emptyset$.  As $S$ is detecting, we deduce $L'=Q(S)$.
\end{pf} 

Here are a few examples of detecting root systems.
For instance, it follows from Table~\ref{tab:venkov} that ${\bf A}_m$ is detecting if, and only if, there is no integer $1 \leq i \leq m$ such that $m+1\,\,|\,\,i^2$ and $i(m+1-i)>2(m+1)$. This holds in particular if $m+1$ is square free or for all $m\leq 10$. Similarly, ${\bf D}_m$ is detecting unless we have $m \equiv 4 \bmod 8$ and $m \neq 4$, 
and ${\bf E}_m$ is detecting for $6 \leq m \leq 8$. Here is another example.
% Cas Am. soit u = gcd(m+1,i). 
\begin{example} 
\label{ex:detecting}  $m{\bf A}_1\, n{\bf A}_2$ is detecting if, and only if, we have $0 \leq m,n \leq 5$ and
$(m,n) \neq (2,3), (4,3)$.
{\rm Indeed, the only nonzero value of ${\rm qm}$ on ${\rm res}\, {\rm A}_1$ {\rm (}resp. ${\rm res}\, {\rm A}_2${\rm )} 
is $\frac{1}{4}$ {\rm (}resp. $\frac{1}{3}${\rm )}}.
\end{example}

For a given $R$, we will often have to find a sub system $S \subset R$ as large as possible such that $(R,S)$ is safe,
which usually reduces to a problem in coding theory. For instance for $R \simeq n {\bf A}_1$, it amounts to ask
for the maximal integer $m$ such that for any {\it even linear binary code} $I$ in 
$(\Z/2)^n$ with minimal distance $\geq 6$, there is a partition $\{1,\dots,n\}=S \coprod T$ with $|S|=m$ and such that the natural projection $(\Z/2)^n \rightarrow (\Z/2)^T$ is injective on $I$. We leave as an exercise to the reader 
to check the following assertion.

\begin{example} 
\label{ex:safena1}
Assume either $n\leq 8$ and $m \leq n-1$, or $n \leq 10$ and $m \leq n-2$, then $(n {\bf A}_1, \,m {\bf A}_1)$ is safe. %(Proof left to the reader)
\end{example}

 \subsection{Classification of the $d$-kernels of ${\bf ADE}$ root systems}\label{subsect:dkerdesc}
% We now explain how to find all the $d$-kernels of a given root system. 
We first reduce to the irreducible case.

\begin{lemma} Let $R$ be an ${\bf ADE}$ root system, $R= \coprod_i R_i$ its irreducible decomposition, $d \geq 1$ an integer, $S$ a sub root system of $R$, and $S_i:=S \cap R_i$. Then $S$ is a $d$-kernel of $R$ if, and only if, there are divisors $d_i$ of $d$ with ${\rm lcm} \{ d_i\}_{i \in I}=d$ and such that $S_i$ is a $d_i$-kernel of $R_i$ for each $i$. 
\end{lemma}

\begin{pf} 
Just use $R=\coprod_i R_i$ and ${\rm Q}(R) = \bigoplus_i {\rm Q}(R_i)$.
\end{pf}

Lemma~\ref{visrs} can be viewed as a classification of the $d$-kernels of ${\rm I}_n$.
It would be easy to classify the $d$-kernels of ${\bf A}_n$ and ${\bf D}_n$ 
using a similar same method: see Remark~\ref{dkerAD} for the result. 
We follow a different approach, which works in all cases including type ${\bf E}$, 
inspired by classical works of Borel-de-Siebenthal and Dynkin. \ps

%Note that any linear map ${\rm Q}(R) \rightarrow \Z/d$ is the reduction modulo $d$
%of some linear map ${\rm Q}(R) \rightarrow \Z$, which in turn has the form $x \mapsto x \cdot \xi$ for some
%$\xi \in {\rm Q}(R)^\sharp$. \ps\ps
%({\rm Case} ${\bf A}_n$) Any linear map $\varphi : {\rm A}_n \rightarrow \Z/d$ has the form 
%$\varphi(v) = \sum_{i=1}^{n+1} x_i \xi_i \bmod $ for some $0 \leq \xi_i <d$. 
%Up to ${\rm S}_{n+1} \subset {\rm O}({\rm A}_n)$, we may assume $\xi_1 \leq \xi_2 \leq  \dots \leq \xi_n$.
%By looking at ${\rm R}_2({\rm A}_n)$, we see that there is an integer partition $n+1=a_1+a_2+\cdots$ with ${\rm R}_2({\rm ker}\, \varphi) \simeq {\bf A}_{a_1-1} \, {\bf A}_{a_2-1} \, \dots$. Any integer partition of $n+1$ with length $\leq d$ appears this way.\ps\ps
%({\rm Case} ${\bf D}_n$) A linear map $\varphi : {\rm D}_n \rightarrow \Z/d$ has the form $x \mapsto 
%x.\xi \bmod d$ with either $\xi \in {\rm I}_n$, or $2 \xi \in {\rm I}_n$ with odd coordinates. 
%In the first case, we have seen that there is an integer partition $n=m+m'+a_1+a_2+\cdots$ with ${\rm R}_2({\rm ker}\, \varphi) \simeq {\bf D}_m \, {\bf D}_{m'} \, {\bf A}_{a_1-1} \, {\bf A}_{a_2-1} \, \dots$.
%Any integers $m$ and $m'$ with $m+m' \leq n$ and $m'=0$ for $d$ odd, and any
%integer partition of $n-m-m'$ with length $< d/2-1$, appear this way.
%In the second case, up to ${\rm O}({\rm I}_n) \subset {\rm O}({\rm D}_n)$ 
%we may assume $1 \leq 2\xi_1 \leq 2\xi_2 \leq  \dots \leq 2\xi_n \leq d$, 
%and we reach the same conclusion 
%with the extra constraint $m=0$.\ps\ps

Fix an irreducible ${\bf ADE}$ root system $R$ of rank $n\geq 1$.
Choose a positive root system $R^+ \subset R$, 
with associated basis $\{\alpha_i\}_{1 \leq i \leq n}$ and dual basis $\{ \varpi_i \}_{1\leq i \leq n}$ in 
the {\it weight lattice} ${\rm Q}(R)^\sharp$ ({\it fundamental weights}).
In all examples below, we choose the same numbering of simple roots as \cite{bourbaki} to fix ideas.
For each $1 \leq i \leq n$, set $h_i \,=\, {\rm max} \{  \varpi_i.\alpha\, \, |\, \, \alpha \in R\}$.
Recall that $\widetilde{\alpha}=\sum_{i=1}^n \,h_i\, \alpha_i$ is in $R^+$ ({\it highest root}).
We also set $\alpha_0=-\widetilde{\alpha}$, $h_0=1$, $\varpi_0=0$ and $I=\{0,1,\dots,n\}$.  
For each subset $J \subset I$, we define a sub root system $R_J \subset R$ by the formulas 
\begin{equation} 
\label{rsrtJ}
Q_J \,=\, \sum_{i \in I-J} \,\Z\, \alpha_i\, \, \, \, \, {\rm and}\, \, \, \, \, R_J = {\rm R}_2(Q_J).
\end{equation} 
As $Q_J$ is a subgroup of ${\rm Q}(R)$ generated by roots we have $Q_J= {\rm Q}(R_J)$.

%In the case $0 \in J$, we recognize the parabolic root system with basis $\{\alpha_i}_{i \in I-J}$,
%and whose Dynkin diagram is obtained from that of $R$ by removing the vertices $\alpha_i$ for all nonzero 
%$i \in J$. 

\begin{lemma}\label{propRJ}  Assume we have $J \subset I$ with $J \neq \emptyset$.
Then  $\{ \alpha_i\, |\, i \in I-J\}$ is a basis of the root system $R_J$, and we have 
\begin{equation}\label{descRJ} R_J=\{ \alpha \in R \, \,|\,\, \exists n \in \{-1,0,1\},\, \, \forall i \in J,\,\, \varpi_i.\alpha = n h_i, \}. \end{equation}
\end{lemma}

\begin{pf} 
The sum in \eqref{rsrtJ} is direct since $J \neq \emptyset$.
Note that for $\alpha \in R$ and $n \in \Z$, the condition $\varpi_i.\alpha = n h_i$ for all $i \in J$ implies
$n \in \{ -1, 0, 1\}$, and even $n=0$ in the case $0 \in J$, by definition of the $h_i$.
All assertions are then simple consequences of the fact that for any $\alpha \in R^+$, both $\alpha$ and
$\widetilde{\alpha}-\alpha$ are a finite sum of $\alpha_i$ for $i \neq 0$. 
\end{pf}
The Dynkin diagram of $R_J$ is thus obtained from the affine Dynkin diagram of $R$ by removing each $\alpha_j$ for $j \in J$. Familiar cases include the {\it parabolic} case $0 \in J$, and the case
$|J|=1$ (see \cite[Ch. VI, \S 4, Ex. 4]{bourbaki}).

\begin{prop}\label{classdkernel} 
For $d\geq 1$, the $d$-kernels of $R$ are the $w(R_J)$ with $w \in {\rm W}(R)$ and $J \subset I$ 
such that there exists integers $x_j>0$ for all $j \in J$ with $\sum_{j \in J} x_j h_j =d$ and 
${\rm gcd}(d,\,\,\{ x_j\,\,|\,\,j \in J\}\, \, )=1$.
\end{prop}

\begin{pf} 
Any linear form $\varphi : {\rm Q}(R) \rightarrow \Z/d$ has the form $\varphi(x)=\xi.x \bmod d$ 
for some $\xi \in {\rm Q}(R)^\sharp$. We also have $\varphi(w(x))=(w^{-1}\xi).x \bmod d$ for $w \in {\rm W}(R)$. Applying the {\it affine Weyl group} ${\rm Q}(R) \rtimes {\rm W}(R)$ to $\xi/d$, 
we may assume that $\xi$ is in the {\it closed alcove} defined by $0 \leq \xi . \alpha \leq d$ for all $\alpha \in R^+$. In particular,
\begin{equation} \label{r2kerphi} R' := {\rm R}_2({\rm ker}\, \varphi)= \{ \alpha \in R\,\,|\,\, \xi . \alpha \in \{-d,0,d\} \}.\end{equation} 
On the other hand, we may write $\xi = \sum_{i \in I} x_i \varpi_i$ in a unique way such that $\sum_{i \in I} x_i h_i = d$, namely $x_i=\xi.\alpha_i \in \Z$ for $i\neq 0$ and $x_0 = d - \sum_{i=1}^n x_i h_i$
({\it Kac coordinates}). 
The alcove inequalities are equivalent to $x_i \geq 0$ for all $i \in I$. 
Set $J=\{i \in I\, \,|\, \, x_i \neq 0\}$, so that $\sum_{i \in J} \,x_i \,h_i \,=\, d$ and 
$x_i \geq 1$ for all $i \in J$. For $\alpha \in R^+$, the element 
$\xi.\alpha \,=\, \sum_{i \in J}\, x_i \,\varpi_i .\alpha$ is $0$ (resp. $d$)
if and only if we have $\varpi_i . \alpha=0$ (resp. $\varpi_i . \alpha=h_i$) for all $i \in J$.
Formulas \eqref{descRJ} and \eqref{r2kerphi} show then $R_J=R'$. 
Also, $\varphi$ is surjective if and only if the set of $x_j$ with $j \in J$ is coprime with $d$.
Conversely, any $R_J$ as in the statement is a $d$-kernel by this same analysis, setting 
$\xi = \sum_{i \in J} x_i \varpi_i$ and $\varphi(x) = \xi.x \bmod d$.
\end{pf}

\begin{remark} \label{classdkernelrem}
{\rm  
Assume $J \subset I$ contains some element $i$ with $h_i=1$ ({\it e.g.} $0 \in J$) and $|J|>1$.
The existence of $x_j$ as in Proposition~\ref{classdkernel} is then equivalent to 
$\sum_{i \in J} h_i \leq d$, setting $x_j=1$ for all $j \neq i$. 
}
\end{remark}

\begin{remark} 
{\rm If we omit the condition that the set of $x_j$ is coprime to $d$, 
we obtain a description of all the $d'$-kernels of $R$ for some $d'$ dividing $d$.}
\end{remark}

\begin{example}
\label{dkerAD}
{\rm 
\begin{itemize}
\item[(${\bf A}_n$)] 
The affine diagram of ${\bf A}_n$ ($n\geq 1$)
is a ``circle'' with $h_i=1$ for all $i$. By removing $s \geq 2$ vertices, we see that the $d$-kernels of ${\bf A}_n$ with $d\geq 2$ are isomorphic to ${\bf A}_{a_1}\,{\bf A}_{a_2}\,\dots {\bf A}_{a_s}$ with $2 \leq s \leq d$, $a_i \geq 0$ and $n+1 = s + \sum_{i=1}^s a_i$.
\item[(${\bf D}_n$)] The affine diagram of ${\bf D}_n$ ($n\geq 4$) is a ``bone'', 
with boundary heights $h_0=h_1=h_n=h_{n-1}=1$ and $h_i=2$ otherwise. 
Removing for instance $s \geq 2$ inner vertices, we see that for $d$ even, the root systems of the form
 ${\bf A}_{a_1}\,{\bf A}_{a_2}\,\dots {\bf A}_{a_{s-1}} \, {\bf D}_m \,{\bf D}_{m'}$ with 
 $2 \leq s \leq d/2$, $a_i \geq 0$, $m,m'\geq 2$ and $n +1 =m+m'+ s + \sum_{1 \leq i \leq s-1} a_i$, 
 are $d$-kernels of ${\bf D}_n$.
 % We leave as an exercise to the reader to check that $d$-kernels of ${\bf D}_n$ with $d\geq 3$ 
%are the root systems of the form
 %${\bf A}_{a_1}\,{\bf A}_{a_2}\,\dots {\bf A}_{a_s} \, {\bf D}_m \,{\bf D}_{m'}$ with 
% $m,m'\geq 0$, $n=m+m'+ \sum_{i=1}^{s} (a_i+1)$, $a_i \geq 0$ and: 
%either $1 \leq s \leq (d-1)/2$, with furthermore $m=0$ for $d$ odd, {\color{red} or  $2 \leq s \leq d/2$ and $m'=0$, 
% with furthermore $m=0$ for $d$ even, and $s \equiv 0 \bmod 4$ in the case $s=[d/2]$.}
\end{itemize}
}
\end{example}

\begin{example} 
{\rm 
({\it $2$-kernels}) The $2$-kernels of ${\bf A}_n$ with $n\geq 1$ are the ${\bf A}_p\, {\bf A}_{n-1-p}$ with $0 \leq p \leq n-1$. Those of ${\bf D}_n$ with $n\geq 4$ are ${\bf A}_{n-1}$ and the 
${\bf D}_p\,{\bf D}_{n-p}$ for $1 \leq p \leq n-1$. The $2$-kernels of ${\bf E}_n$ are ${\bf A}_1\,{\bf A}_5$ and ${\bf D}_5$ for $n=6$,  ${\bf A}_1\,{\bf D}_6$, ${\bf A}_7$ and ${\bf E}_6$ for $n=7$, ${\bf A}_1\, {\bf E}_7$ and ${\bf D}_8$ for $n=8$.
}
\end{example}

%Examples with d=3 !

% Many properties of the unimodular lattice 
% are influenced by the same property for the much simpler lattice ${\rm M}_d(x)$. 

\section{King's mass formulas}\label{sect:mass}

Our aim in this section is to discuss various {\it mass formulas}.
It will be convenient to recall first a few elementary but useful concepts about {\it groupoids}.

\subsection{Preliminaries on groupoids} 
\label{subsec:gpds}

A groupoid is a category $X$ all of whose arrows are isomorphisms. 
%The motivating case here is the groupoid with objects the integral lattices,\footnote{
% it is often convenient not to assume that lattices are embedded in 
% a fixed Euclidean space, and we will tend to do so in this subsection.}  
% and with morphisms the isometries.
We say that such an $X$ is {\it finite} if it has finitely many objects up to isomorphism,
and if ${\rm Aut}(x)$ is a finite group for all objects $x$ of $X$. 
We then denote by $\underline{X}$ the finite set of isomorphism classes of objects in $X$.
The {\it class number} of a finite groupoid $X$ is ${\rm h}(X)=|\underline{X}|$, and its 
{\it mass} is the rational number 
\begin{equation}\label{defmassX} {\rm mass}(X) = \sum_x 1/|{\rm Aut}(x)|,\end{equation}
where $x$ runs among representatives of $\underline{X}$.
Two {\it equivalent} finite groupoids (in the sense of categories) have the same class numbers and masses.
Assume now $F : X \rightarrow Y$ is any {\it morphism} of groupoids, which simply means that $F$ is a functor.
If $y$ is an object of $Y$, we define the (naive) {\it fiber} of $F$ at $y$ as the full subcategory 
$F^{-1}y$ of $X$ whose objects $x$ satisfy $F(x) \simeq y$. 
This is a groupoid which only depends on the isomorphism class of $y$.
If $X$ is finite, then so are the fibers of $F$, and we trivially have
\begin{equation} \label{slicemassX} {\rm mass}(X) = \sum_y {\rm mass}(F^{-1}y), \end{equation}
where $y$ runs among representatives of $\underline{Y}$. \ps 

A typical finite groupoid is obtained as follows. 
Let $G$ be a finite group acting on a finite set $S$. 
We denote by $[S/G]$ the finite groupoid with set of objets $S$, 
and with morphisms $s \rightarrow s'$ the set of elements $g \in G$ with $gs=s'$
(with obvious compositions). 
Then ${\rm h}([S/G])$ is the number of $G$-orbits in $S$ and we have 
\begin{equation} \label{massSsurG} {\rm mass}([S/G]) = |S|/|G| \end{equation} 
by the {\it orbit-stabilizer} formula. In the proofs below, we will usually compute the mass of a given $X$ fiberwise, 
using \eqref{slicemassX} for a suitable morphism $F : X \rightarrow Y$ and \eqref{massSsurG} by identifying each fiber of $F$
with a concrete $[S/G]$.

\subsection{Minkowski-Siegel-Smith mass formulas}
\label{subsect:minkosiegsmith}

Recall that the {\it mass} of an integral lattice $L$ is defined by 
$${\rm mass}(L) = \frac{1}{|{\rm O}(L)|}.$$
More generally, if $\mathcal{L}$ is any collection of integral lattices (embedded in some fixed Euclidean space or not), 
but consisting of finitely many isometry classes $L_1,...,L_h$, the {\it mass} of $\mathcal{L}$ 
may be defined and ${\rm mass}(\mathcal{L}) = \sum_i \frac{1}{|{\rm O}(L_i)|}$. Alternatively, this is the mass of the 
finite groupoid (that we will still denote by $\mathcal{L}$) with objects $\mathcal{L}$ and with morphisms given by isometries. \ps

When $\mathcal{L}$ is a genus of integral lattices, 
the famous {\it Minkowski-Siegel-Smith mass formula} gives 
a closed formula for ${\rm mass}(\mathcal{L})$. 
Recall that $\mathcal{L}_n$ is the disjoint union of the two genera
$\mathcal{L}_n^{\rm odd}$ and $\mathcal{L}_n^{\rm even}$ consisting of odd or even lattices: 
 see {\it e.g.} \cite{conwaysloanemass} for concrete formulas for their respective masses. 
\ps

\subsection{King's mass formulas}
\label{subsect:kingsmassform}

For any groupoid $\mathcal{L}$ of integral Euclidean lattices $L$, 
and any ${\bf ADE}$ root system $R$, 
we denote by $\mathcal{L}(R)$ the groupoid of objects $L$ of $\mathcal{L}$
whose root system ${\rm R}_2(L)$ is isometric to $R$.
%(of course, $\mathcal{L}(R)$ only depends on the isomophism class of $R$). 
King explains in \cite{king} an algorithm to compute, for each root system $R$, the mass 
\begin{equation} \label{mnevenR} {\rm mass}(\mathcal{L}_n^{\rm even}(R)).\end{equation} 
He uses for this the expression for the ``mass-weighted'' sum of the Siegel theta series of rank $n$ even unimodular lattices as a Siegel Eisenstein series (Siegel-Weil formula), Katsurada's formula for its Fourier coefficients, 
as well as elementary properties of root lattices. 
Using the computer, he was then able to determine the quantity \eqref{mnevenR} for all $R$ up to $n=32$.
\ps

As explained in \S 4 of \cite{king}, using ideas of Conway and Sloane
in \cite{conwaysloane} Chap. 16 one can deduce from
this computation the mass of  $\mathcal{L}_n^{\rm odd}(R)$ for all $R$ and $n \leq 30$.
As the argument is only sketched {\it loc. cit.}, and with tables not freely available, 
we now give more details about how this computation can be made. 
Actually, we follow a different method and give more general formulas. \ps

Assume $m>1$ is an integer, and $R$ and $R'$ are ${\bf ADE}$ root systems. Denote by ${\rm N}_m(R,R')$
the number of root systems $S \subset R$ which are isomorphic to ${\bf D}_m$, 
saturated\footnote{This condition is empty for $m \neq 8$, and for $m=8$ it means that the irreducible component of $R$ containing $S$ is not of type ${\bf E}_8$.} in $R$, and with $S^\perp \cap R$ isometric to $R'$ (see \S \ref{subsect:rootsystems}).
This integer only depends on the isomorphism class of $R$ and $R'$. 
Its concrete determination is easily deduced from Table 4 in \cite[\S 8]{king}.
Recall the root lattice ${\rm D}_n$ for $n>1$ (see Sect. \ref{notations} (v)). 
The integer $|{\rm O}({\rm D}_n)|$ is $|{\rm O}({\rm I}_n)|=2^n\, n!$ for $n\neq 4$, 
and $3 \cdot 2^4 \,4!$ for $n=4$ ({\it triality}). 

%For $n\geq 1$ an integer, set $\delta_n=1$ for $n\neq 4$, and $\delta_4=3$.
%$f_n=1$ for $n \not \equiv 4 \bmod 8$, and $f_n=3$ otherwise. We also set $f'_n=f_n$ for $n \neq 4$, 
%and $f'_4=1$.

\begin{prop} 
\label{kingform} 
Fix integers $n,m>1$ with $m+n \equiv 0 \bmod 8$. We have
$${\rm mass}(\mathcal{L}_n^{\rm odd}(R)) \,=\, \frac{|{\rm O}({\rm D}_m)|}{2} \,\sum_{R'} {\rm N}_m(R,R')\, {\rm mass}(\mathcal{L}_{n+m}^{\rm even}(R')),$$
where $R$ is any {\bf ADE} root system of rank $\leq n$, and $R'$ runs among all the isomorphism classes of ${\bf ADE}$ root systems of rank $\leq n+m$.
\end{prop}

Our aim now is to prove this proposition. 
Denote by $\mathcal{D}_n$ the groupoid of (abstract) even lattices $D$ of rank $n$ such that 
the finite quadratic space ${\rm res}\, D$ is isomorphic to ${\rm res}\, {\rm D}_n$.
%The structure of ${\rm q}: {\rm res}\, {\rm D}_n \rightarrow \Q/\Z$ is simple and well-known:
%for $n$ odd it is isomorphic to $ \Z/4 \rightarrow \Q/\Z, \,x \mapsto \frac{n}{8} x^2$, 
%and for $n=2k$ even it is isomorphic to $\Z/2 \times \Z/2 \rightarrow \Q/\Z,\, (x,y) \mapsto \frac{k}{4}(x^2+y^2)$ (case $k$ odd), $\frac{1}{2} xy$ (case $k \equiv 0 \bmod 4$) or $\frac{1}{2}(x^2+xy+y^2)$ (case $k \equiv 2 \bmod 4$).
Set $f_n=1$ for $n \not \equiv 4 \bmod 8$, and $f_n=3$ otherwise.\ps

\begin{lemma}  
\label{ma1} 
For any $n>1$ and any ${\bf ADE}$ root system $R$ we hve
$${\rm mass}(\mathcal{L}_n^{\rm odd}(R)) \,=\, f_n\, {\rm mass}(\mathcal{D}_n(R)).$$
\end{lemma}

\begin{pf} If $L$ is an integral lattice, it has a largest even sublattice denoted $L^{\rm even}$, and which is of index $2$ if $L$ is odd.
As ${\rm D}_n=({\rm I}_n)^{\rm even}$, and $\mathcal{L}_n^{\rm odd}$ is the genus of ${\rm I}_n$, it follows that $L \mapsto L^{\rm even}$ defines a morphism
of groupoids ${\rm ev}_n: \mathcal{L}_n^{\rm odd} \rightarrow \mathcal{D}_n$.  As ${\rm R}_2(L)={\rm R}_2(L^{\rm even})$ for all integral lattices $L$, and by Formula \eqref{slicemassX}, it is enough to show that for all
$D$ in $\mathcal{D}_n$ we have 
\begin{equation} \label{fibevn} {\rm mass}({\rm ev}_n^{-1} D)\,=\,f_n\, {\rm mass}(D). \end{equation}
	Fix $D$ in $\mathcal{D}_n$ and let $S$ be the set of odd unimodular lattices $L$ in $V=D \otimes \R$
with $L^{\rm even}=D$. Any object in ${\rm ev}_n^{-1}D$ is isomorphic to some element of $S$. Better,
any isometry $L \rightarrow L'$ with $L,L'$ in $S$ is induced by an isometry of $D=L^{\rm even}=(L')^{\rm even}$.
The functor $L \mapsto L$, $[S/{\rm O}(D)] \rightarrow {\rm ev}_n^{-1}D$, is thus an equivalence. By Formula~\eqref{massSsurG} 
it only remains to show $|S|=f_n$. By Sect. \ref{notations} (iii) and (iv), the map $\beta_D : L \mapsto L/D$ defines a bijection between $S$ and the 
set of order $2$ subgroups of ${\rm res}\, D$ generated by an element $v$ with $v.v \equiv 0$ but ${\rm q}(v) \not \equiv 0$.
An inspection of ${\rm res}\, D \simeq {\rm res}\, {\rm D}_n$ (Table \ref{tab:venkov}) concludes.
\end{pf}

Proposition \ref{kingform} is a special case of the following more general proposition.
It will be convenient to make first a definition.

\begin{definition} 
\label{def:MandMprime}
Let $R$ be an ${\bf ADE}$ root system and set $Q={\rm Q}(R)$. We say that $R$ {\rm satisfies 
${\rm (M)}$} if the the finite quadratic space ${\rm res}\, Q$ is anisotropic, or equivalently, 
if $Q$ is saturated in any even lattice containing it. We say that $R$ {\rm  satisfies 
${\rm (M')}$} if we have ${\rm qm}(x)=1$ for all isotropic $x \in {\rm res}\, {\rm Q}(R)$, or equivalently, 
if all the even integral lattices between $Q$ and $Q^\sharp$ are root lattices. 
\end{definition}

\begin{example} 
\label{ex:MMprime}
{\rm It is clear that {\rm (M)} implies {\rm (M')}. 
The irreducible root systems satisfying {\rm (M)} are the following : 
${\bf A}_m$ in the case $m+1$ is either squarefree or of the form $4q$ with $q$ odd and squarefree, ${\rm D}_m$ for $m \not \equiv 0 \bmod 8$, ${\bf E}_6$, ${\bf E}_7$ and ${\bf E}_8$. 
The root systems ${\bf D}_8$, ${\bf A}_m$ for $m \leq 14$, and $m{\bf A}_1$ for $m\leq 7$, satisfy {\rm (M')}. 
}
\end{example}
\ps

%
% when does there is some 1<=i<=n such that  2(n+1) divides i(n+1-i) ?
% set k = gcd(i,n+1), i = ki' and n+1=ku, with gcd(i',u)=1.
% we have 1<=i'<u.
% we have 2ku | ki'(ku-ki'), so 2u | i'k(u-i'), and then u | k.
% write k=uk', so n+1=u^2 k'. So n+1 not squarefree.
% remains 2 | i'k'(u-i'), so either i', k', u-i' is even.
% always possible if u odd or k' even.
% if u even and k' odd, i' even, but then u>i' implies u>2.
% reciproques given by i'=1.

For {\bf ADE} root systems $R_1,R_2,R_3$ we denote by ${\rm N}(R_1,R_2,R_3)$ 
the number of saturated sub root systems $T$ of $R_2$ satisfying $T \simeq R_1$ and $T^\perp \cap R_2 \simeq R_3$. When $R_1$ satisfies {\rm (M')}, then ${\rm N}(R_1,R_2,R_3)$ can be deduced\footnote{This is especially simple if $R_1$ satisfies (M), since this number is just the number of sub root systems $T$ of $R_2$ satisfying $T \simeq R_1$ and $T^\perp \cap R_2 \simeq R_3$.} from Table 4 in 
{\rm \cite[\S 8]{king}}. \ps

\begin{prop}
\label{kingform2} 
Let $R_0$ be a rank $m$
%, residually anisotropic, 
${\bf ADE}$ root system satisfying {\rm (M')}, let $n >1$ an integer with $m+n \equiv 0 \bmod 8$, and $\mathcal{G}$ the groupoid\footnote{We shall not use it, but $\mathcal{G}$ is actually a single genus, by \cite[Cor. 1.9.4]{nikulin}. } of rank $n$ even Euclidean lattices $L$ with ${\rm res}\, L \simeq - {\rm res}\, {\rm Q}(R_0)$ (as quadratic spaces). Then for any ${\bf ADE}$ root system $R$ we have
$${\rm mass}(\mathcal{G}(R))\, =\, \frac{|{\rm O}({\rm Q}(R_0))|}{|{\rm O}({\rm res}\,R_0)|}\, \,\sum_{R'} \,{\rm N}(R_0,R',R) \,\mathrm{mass}(\mathcal{L}_{m+n}^{\rm even}(R'))$$
where $R'$ runs over all isomorphism classes of ${\bf ADE}$ root systems of rank $\leq m+n$.
\end{prop}

In the statement above, ${\rm O}({\rm res}\,R_0)$ denotes the isometry group of the finite {\it quadratic} space 
${\rm res}\, R_0 \,= \,{\rm res}\, {\rm Q}({\rm R}_0)$ (equipped with the form ${\rm q}$ of \S \ref{notations} (iv)).

\begin{pf} 
(Proposition \ref{kingform2} $\Rightarrow$ Proposition \ref{kingform}).
%Choose $m,n>1$ with $m+n \equiv 0 \bmod 8$ as in Proposition \ref{kingform}.
Observe that we have ${\rm res}\, {\rm D}_n\, \simeq - \,{\rm res}\, {\rm D}_m$, so we 
are in the situation of Proposition \ref{kingform2} with $R_0 \simeq {\bf D}_m$, 
${\rm Q}(R_0) \simeq {\rm D}_m$,  and $\mathcal{G}= \mathcal{D}_n$. As the isometry group of 
${\rm res}\, {\rm D}_m$ is isomorphic to ${\rm S}_3$ 
for $m \equiv 4 \bmod 8$, and to $\Z/2$ otherwise, we have 
$|{\rm O}({\rm res \,D}_m)|=2 f_m$. We conclude by Lemma~\ref{ma1}, $f_n=f_m$ and ${\rm N}({\bf D}_m,R',R)={\rm N}_m(R',R)$. 
\end{pf}

\begin{pf}
We now prove Proposition \ref{kingform2}. Denote by 
$\mathcal{U}_{n,R_0}$ the groupoid with objects the 
pairs $(U,Q)$ where $U$ is an even unimodular lattice of rank $m+n$ 
and $Q$ is a saturated subgroup of $U$ isomorphic to the root lattice ${\rm Q}(R_0)$, 
and with morphisms $(U,Q) \rightarrow (U',Q')$ the isometries $U \rightarrow U'$ sending $Q$ onto $Q'$.
For $(U,Q)$ in $\mathcal{U}_{n,R_0}$ the orthogonal $Q^\perp$ of $Q$ in $U$ is an object of $\mathcal{G}$
 since we have isomorphisms of finite quadratic spaces 
 ${\rm res} \,Q^\perp\, \simeq - {\rm res} \,Q$ 
by Prop. 1.6.1 in \cite{nikulin}. The map $U \mapsto Q^\perp$ induces 
 a morphism of groupoids ${\rm orth} : \mathcal{U}_{n,R_0} \rightarrow \mathcal{G}$.
 
 \begin{lemma} 
 \label{fiborthnm} 
 For all $D$ in $\mathcal{G}$ we have 
 ${\rm mass}({\rm orth}^{-1}D)= \frac{|{\rm O}({\rm res}\,R_0)|}{|{\rm O}({\rm Q}(R_0))|} {\rm mass}(D)$.
\end{lemma}  

\begin{pf} 
Fix $D$ in $\mathcal{G}$. Set $M = {\rm Q}(R_0) \perp D$, $V=M \otimes \R$, and
let $S$ be the set of even unimodular lattices $U$ of $V$ containing $M$ and with ${\rm Q}(R_0)$ saturated in $U$.
This set has a natural action of the group
$G={\rm O}({\rm Q}(R_0)) \times {\rm O}(D)$, and the map
$L \mapsto (L,{\rm Q}(R_0))$ trivially induces an equivalence 
$[S/G] \isomo {\rm orth}^{-1}D$. 
We have thus ${\rm mass}({\rm orth}^{-1}D)|=|S|/|G|$.
By Sect. \ref{notations} (iii) and (iv), the map $\beta_M : U \mapsto U/M$ induces a bijection between $S$ and the 
set of quadratic Lagrangians in ${\rm res}\,M \,= \,{\rm res}\, {\rm Q}(R_0) \perp {\rm res}\, D$
which are transversal to ${\rm res}\, {\rm Q}(R_0)$.
But in a finite quadratic space $A \perp B$ with $B \simeq -A$, there are 
$|{\rm O}(A)|$ Lagrangians which are transversal to $A$, namely the 
$I_\varphi=\{ x+\varphi(x), x \in A\}$, with $\varphi$ any isomorphism $A \isomo -B$, 
thus $|S|=|{\rm O}({\rm res}\,R_0)|$.
\end{pf}

%Any 
%$(U,Q)$ in $\mathcal{U}_{n,m}$ is isomorphic to $(U,{\rm D}_m)$ for some $U$ in $S$.
%Also,  for $L,L'$ in $S$, any isometry $L \rightarrow L'$ preserving ${\rm D}_m$ sends 
%$D={\rm D}_m^\perp \cap L = {\rm D}_m^\perp \cap L'$ into itself, so $L \mapsto (L,{\rm D}_m)$
%defines an equivalence $S/\!\!/{\rm O}({\rm D}_m) \isomo {\rm orth}_{n,m}^{-1}

Denote by $\mathcal{U}_{n,R_0,R}$ the full subcategory of $\mathcal{U}_{n,R_0}$ whose objects
$(U,Q)$ satisfy ${\rm R}_2(Q^\perp) \simeq R$, i.e. $Q^\perp$ is an object of $\mathcal{G}(R)$.
This is the disjoint union of all the fibers of ${\rm orth}$ over $\mathcal{G}(R)$, so
Lemma \ref{fiborthnm} and Formula \eqref{slicemassX} imply
\begin{equation}\label{ma2} {\rm mass}(\mathcal{U}_{n,R_0,R}) \,=\, \frac{|{\rm O}({\rm res}\,R_0)|}{|{\rm O}({\rm Q}(R_0))|}\, {\rm mass}(\mathcal{G}(R)).\end{equation}
Consider now the forgetful functor ${\rm forg} : \mathcal{U}_{n,R_0,R} \rightarrow \mathcal{L}_{m+n}^{\rm even}$,
$(U,Q) \mapsto U$.
Fix $U$ in $\mathcal{L}_{m+n}^{\rm even}$, and consider the set ${\rm S}(R_0,U,R)$ of saturated subgroups $Q$ 
of $U$ satisfying $Q \simeq {\rm Q}(R_0)$ and ${\rm R}_2(Q^\perp) \simeq R$.
Then ${\rm S}(R_0,U,R)$ has a natural ${\rm O}(U)$-action, and the fiber ${\rm forg}^{-1}U$ 
is trivially equivalent to the groupoid $[{\rm S}(R_0,U,R)/{\rm O}(U)]$.

\begin{lemma} 
\label{lem:satroot}
We have $|{\rm S}(R_0,U,R)|={\rm N}(R_0,R',R)$.
 \end{lemma}

\begin{pf}
Set $R'={\rm R}_2(U)$ and let $Q= {\rm Q}(S)$ be a sub root-lattice of $U$, hence of ${\rm Q}(R')$.
If $Q$ is saturated in $U$, then $Q$ is obviously saturated in ${\rm Q}(R')$.
Assume conversely that $Q$ is saturated in ${\rm Q}(R')$. 
Let $Q'$ be the saturation of $Q$ in $U$.
As $Q'$ is even, we have $Q' = {\rm Q'}(S')$ for $S'={\rm R}_2(Q')$ since $R_0$ satisfies (M'). 
But then we have $S' \subset {\rm R}_2(U)=R'$ and thus $Q' \subset {\rm Q}(R')$, and then $Q=Q'$ as 
$Q$ is saturated in ${\rm Q}(R')$, so $Q$ is saturated in $U$. 
%This proves the lemma.
\end{pf}

We have thus ${\rm mass}({\rm forg}^{-1}U)\, =\,{\rm N}(R_0,R',R)\, {\rm mass}(U)$ by Formula~\eqref{massSsurG} and Lemma~\ref{lem:satroot}.
Using Formula \eqref{slicemassX}, we deduce
$${\rm mass}(\mathcal{U}_{n,R_0,R}) \,=\, \sum_{R'}\, {\rm N}(R_0,R',R)\, {\rm mass}(\mathcal{L}_{m+n}^{\rm even}(R')),$$
the sum being over all isomorphism classes of {\bf ADE} root systems $R'$, and conclude by \eqref{ma2}.\end{pf}

%\begin{remark}\label{rem:compNfromking} {\rm ({\it On the computation of ${\rm N}(R_0,R',R)$})
%Assume that the only quadratic isotropic subspace of ${\rm res}\, Q(R_0)$ is $\{0\}$, or equivalently, that $R_0$ is saturated in any even lattice containing it. This holds for instance for $R_0 \simeq {\bf D}_m$ when $m \not \equiv 0 \bmod 8$. Then ${\rm N}(R_0,R',R)$ is the number of sub root systems of $R'$ isometric to $R_0$ and with orthogonal complement isometric to $R$, hence it can be deduced from Table 4 in 
%{\rm \cite[\S 8]{king}}. More generally, one can deduce ${\rm N}(R_0,R',R)$ from this table if 
%all the even integral lattices between ${\rm Q}(R_0)$ and ${\rm Q}(R_0)^\sharp$ are root lattices, or equivalently, if we have ${\rm qm}(x)=1$ for all isotropic $x \in {\rm res}\, {\rm Q}(R_0)$.}
%\end{remark} 

As a consequence of King's results and Proposition~\ref{kingform2}, 
we obtain new lower bounds for the number of isometry classes in the genera of even lattices of rank $32-r$ and residue isometric to $- {\rm res}\, {\rm Q}(R_0)$
with $R_0$ an ${\bf ADE}$ root system of rank $r \geq 1$ satisfying (M'). This will be a useful ingredient in forthcoming work with O. Ta\"ibi. As a very simple example, let us consider the genus $\mathcal{G}$ of even lattices of rank $31$ and determinant $2$.
We are in the case $R_0 \simeq {\bf A}_1$, $m=1$ and $n=31$. The numerical application of the Proposition above shows that there are exactly $18\,437$ root systems in this genus and using the same argument as in~\cite[\S 6]{king} we obtain :

\begin{cor} 
There are at least $6\,678\,411\,375$ even lattices of determinant $2$ in dimension $31$.
Moreover, the mass of those lattices having an empty root system is $11882632915662079/591224832 \simeq 
20098331.92$.
\end{cor}

\begin{comment}
\begin{table}[H]
\tabcolsep=5pt
{\scriptsize \renewcommand{\arraystretch}{1.8} \medskip
\begin{center}
\begin{tabular}{c|c|c|c|c|c|c|c}
$n$ & $1$ & $2$ & $3$ & $4$ & $5$ & $6$ & $7$ \\
\hline
 $#$root systems & $ $18437$  \\
 lower bound &  
 $6,6 \cdot 10^9$ & $5,5 \cdot 10^8$ & $3,7 \cdot 10^7$ & $2,5 \cdot 10^6$ & $253739$ & $39663$ 
\end{tabular} 
\end{center}
} 
{\small \caption{ For the genus of rank $32-n$ lattices with residue $- {\rm res}\, {\rm A}_n$, 
number of different root systems and lower bound on the number of lattices}}
\label{tab:appkingAn}
\end{table}
\end{comment}

\subsection{Odd unimodular lattices without norm 1 vector}

We denote by $\mathcal{L}'_n$ the groupoid of all $L$ in $\mathcal{L}_n$ with ${\rm r}_1(L)=0$.
Note that in the canonical decomposition $L= A \perp B$ recalled in Formula \eqref{decnorm1} of the introduction, we obviously have ${\rm mass}(L) = \frac{{\rm mass}(B) }{2^{m}m!}$.
Using the convention ${\rm mass}(\mathcal{L}_0)={\rm mass}(\mathcal{L}'_0)=1$ (that is, ${\rm Aut}(0)=\{1\}$ !)
we deduce
\begin{equation}
\label{powersermno} \sum_{n\geq 0} \,{\rm mass}(\mathcal{L}'_n) \,x^n\, = \,e^{-x/2}\, \sum_{n\geq 0} \,{\rm mass}(\mathcal{L}_n)\, x^n.\end{equation}
The relevant numerical values are given in Table \ref{tab:mnu} below.\ps

\begin{table}[htp]
{\tiny \renewcommand{\arraystretch}{1.4} \medskip
%\hspace{-2 cm}
\begin{tabular}{c|c|c|c}
$n$ & ${\rm mass}(\mathcal{L}'_n) $ & $n$ & ${\rm mass}(\mathcal{L}'_n) $ \\
\hline 
$0$ & $\scalebox{.6}{1}$ & $20$ & $\scalebox{.7}{4060488226549/11479871952566228090880000}$  \\
$8$ & $\scalebox{.7}{1/696729600}$ & $21$ & $\scalebox{.7}{138813595637/54497004983156736000000}$ \\
$12$ & $\scalebox{.7}{1/980995276800}$ & $22$ & $\scalebox{.7}{1475568922019/45471119389159682211840}$ \\
$14$ & $\scalebox{.7}{1/16855282483200}$ & $23$ & $\scalebox{.7}{21569773276937492389/28590262351867673365708800000}$  \\
$15$ & $\scalebox{.7}{1/41845579776000}$ & $24$ & 
$\scalebox{.7}{4261904533831299496396870055017/129477933340026851560636148613120000000}$ \\
$16$ & $\scalebox{.7}{5213041/277667181515243520000}$ & $25$ & 
$\scalebox{.7}{103079509578355844357599/37291646545914356563968000000}$ 
\\
$17$ & $\scalebox{.7}{1/49662885888000}$ & $26$ & 
$\scalebox{.7}{15661211867944570315962162816169/34253421518525622105988399104000000}$

 \\
$18$ & $\scalebox{.7}{1073351/32780153373327360000}$ & $27$ & 
$\scalebox{.7}{18471746857358122138056975582390629/121385562506275173338096389324800000}$ \\
$19$ & $\scalebox{.7}{37813/450541700775936000}$ & $28$ & 
$\scalebox{.7}{1722914776839913679032185321786744287148737/16573175467523436999761427022479360000000}$
\end{tabular} 
}
\caption{ {\small The nonzero ${\rm mass}(\mathcal{L}'_n) $ for $0 \leq n\leq 28$}}\label{tab:mnu}
\end{table}

Note also that still in the canonical decomposition $L= A \perp B$ above, we have ${\rm r}_1(B)=0$ 
and thus ${\rm R}_2(L) = {\rm R}_2(A) \coprod {\rm R}_2(B)$ and ${\rm R}_2(A) \simeq {\bf D}_m$.
If $\mathcal{L}_n^{\rm odd'}$ denotes the groupoid of rank $n$ odd unimodular lattices with
no norm $1$ vector, and if $R$ is an {\bf ADE} root system, we have thus
\begin{equation} \label{nonorm1rootmass}{\rm mass}(\mathcal{L}_n^{\rm odd}(R)) = 
\,\sum_{(m,S)} \, \frac{1}{2^m \,m!}\, {\rm mass}(\mathcal{L}_{n-m}^{\rm odd'}(S)) \end{equation}
where $(m,S)$ runs among the pairs with $m$ an integer $\geq 0$ and $S$ an isomorphism class of root systems
such that ${\bf D}_m \coprod S \simeq R$.
Of course, the term ${\rm mass}(\mathcal{L}_{n}^{\rm odd'}(R))$ appears in the sum \eqref{nonorm1rootmass} for $m=0$, so we obtain an expression for it
in terms of King's ${\rm mass}(\mathcal{L}_{32}^{\rm even}(R'))$ for $n \leq 30$:
see {\rm \cite{cheweb}} for tables and Table \ref{tab:nbrs} in the introduction for a few information.

\section{Visible isometries}
\label{sect:visisometry}
\subsection{The visible isometry group of a $d$-neighbor}

Fix a $L \in \mathcal{L}_n$, $d\geq 1$ an integer, an isotropic line $\ell \in {\rm C}_L(\Z/d)$. The group ${\rm O}(L)$ naturally acts on ${\rm C}_L(\Z/d)$ and we denote by 
 ${\rm O}(L; \ell)$ the stabilizer of $\ell$. As it is equivalent to stabilize $\ell$ and its orthogonal ${\rm M}_d(\ell)/dL$, we also have 
 \begin{equation}
 \label{eq:olmdl}
 {\rm O}(L; \ell) = {\rm O}(L) \cap {\rm O}( {\rm M}_d(\ell)).
 \end{equation}
Observe that the action of any $g \in {\rm O}(L;\ell)$ on the line $\ell \simeq \Z/d$ is $\Z/d$-linear, hence given by the multiplication by some element 
$\nu(g) \in (\Z/d)^\times$. Then \begin{equation}
\label{eq:morphnu}
\nu : {\rm O}(L; \ell) \rightarrow (\Z/d)^\times, g \mapsto \nu(g)
\end{equation}
is a group homomorphism. For instance we have $\nu(- {\rm id}_L)=-1$. 
Our aim now is to determine ${\rm O}(L; \ell)$
when $L$ is the lattice ${\rm I}_n$. For this we choose $x \in \Z^n$ such that $x \bmod d$ generates $\ell$, denote by ${\rm X}(\ell)$ the multi set $\{\!\{ \,\,\pm x_1,\,\dots,\,\pm x_n\,\}\!\}$ of elements of $(\Z/d)/\{ \pm 1\}$, and set\footnote{For a multiset $X$ over $(\Z/d)/\{\pm 1\}$ and $\lambda \in (\Z/d)^\times$ we set $\lambda X = \{\!\{ \lambda x \,|\, x \in X\}\!\}$.} 
\begin{equation}
\label{eq:defHl}
{\rm H}(\ell) = \{ \lambda \in (\Z/d)^\times\, \, |\, \,\lambda \,{\rm X}(\ell)\, =\, {\rm X}(\ell)\,\,\}.
\end{equation}
Then ${\rm H}(\ell)$ is a subgroup of $(\Z/d)^\times$ not depending on the choice of $x$.
Recall also from Sect.~\ref{sect:visible} that attached to $\ell$ is a natural partition of the integer $n$ $$n={\rm m}(\ell)+{\rm m'}(\ell)+\sum_{i \in I} {\rm a}_i(\ell).$$

\begin{prop} 
\label{prop:visautOLl}
For $\ell \in {\rm C}_n(\Z/d)$ and $\nu$ as in Formula~\eqref{eq:morphnu} we have 
$${\rm Im}\, \nu \,=\, {\rm H}(\ell)\,\, {\rm and} \,\,{\rm Ker}\,\, \nu \,\,\simeq \,\,{\rm O}({\rm I}_{{\rm m}(\ell)}) \times {\rm O}({\rm I}_{{\rm m}'(\ell)}) \times \prod_{i \in I} {\rm S}_{{\rm a}_i(\ell)}.$$
\end{prop}

\begin{pf} 
Recall that ${\rm O}({\rm I}_n) \simeq \{\pm 1\}^n \rtimes {\rm S}_n$ acts on ${\rm I}_n=\Z^n$ as arbitrary signed permutations of the coordinates. We have thus $\nu(g) \in {\rm H}(\ell)$ for all $g \in {\rm O}(L; \ell)$. Conversely, assume we have $\lambda \in {\rm H}(\ell)$. By definition, there are $\sigma \in {\rm S}_n$ and $e \in \{\pm 1\}^n$ such that  for all $ 1\leq i \leq n$, $\lambda x_i = e_i x_{\sigma(i)}$. In other words, there is $g \in {\rm O}({\rm I}_n)$ with $g (x) = \lambda x \bmod d$, {\it i.e.} with $\nu(g)=\lambda$. We have proved ${\rm Im}\, \nu = {\rm H}(\ell)$. \par 
	An element $\sigma \in {\rm O}(L; \ell)$ is in $\ker \nu$ if, and only if, we have $\sigma(x) \equiv x \bmod d$. The natural action of the signed permutation $\sigma$ on $\{1,\dots,n\}$ preserves thus the partition of $\{1,\dots,n\}$ defined in \S \ref{sect:visible}, with sign changes only allowed on the subsets $D$ and $D'$. 
Indeed, for $i \in \Z/d$  we have $-i \equiv i \bmod d$ if, and only if, $d$ is even and $i \equiv 0$ or $i/2 \bmod d$. This proves the assertion about $\ker \nu$. 
\end{pf}

We now go back to the case of a general $L \in \mathcal{L}_n$. Set $e=1$ for $d$ odd, $e=2$ otherwise.
Then ${\rm O}(L; \ell)$ preserves the lattice ${\rm M}_d(\ell)$, and thus permutes the $e$ cyclic $d$-neighbors of $L$ with line $\ell$, which defines a natural group morphism 
\begin{equation}
\label{eq:visOLl} 
{\rm O}(L; \ell) \rightarrow {\rm S}_e.
\end{equation}

\begin{definition} 
We denote by ${\rm O}(L; \ell)^0$ the kernel of the morphism \eqref{eq:visOLl}.
If $N$ is any of the $1$ or $2$ cyclic $d$-neighbors of $L$ with line $\ell$, then 
${\rm O}(L; \ell)^0$ co\"incides with the subgroup ${\rm O}(N)^{\rm v}:={\rm O}(N) \cap {\rm O}(L)$ of ${\rm O}(N)$. We call it the {\rm visible isometry group} of the $d$-neighbor $N$ of $L$.
\end{definition}

The idea behind these definitions is that visible isometries can be concretely determined from an inspection of the line $\ell$.
Of course, in the case $d$ odd we always have ${\rm O}(L; \ell) = {\rm O}(L; \ell)^0= {\rm O}(N)^{\rm v}$. 
This also holds for $d$ even unless ${\rm O}(L; \ell)$ permutes the $2$ cyclic $d$-neighbors with line $\ell$.
An element of ${\rm O}(N)^{\rm v}$ will also be called a {\it visible isometry} of $N$. \ps

\begin{remark} \label{rem:othervisdef} For $d$ odd and $M=L \cap N$, 
we have ${\rm O}(N)\cap {\rm O}(M) \subset {\rm O}(L)$ since $L$ is the unique cyclic $d$-neighbor of $N$ with visible part $M$
by Proposition~\ref{classdnei}. This shows the equality ${\rm O}(N)^{\rm v}={\rm O}(N) \cap {\rm O}(M)$ for $d$ odd. 

\end{remark}

For the purpose of unimodular hunting, the most basic application of visible isometries is that if we fix $\sigma \in {\rm O}(L)$ and
restrict to isotropic lines $\ell \in {\rm C}_L(\Z/d)$ satisfying $\sigma(\ell)=\ell$, then the corresponding $d$-neighbors of $L$
necessarily contain $\sigma$ as a visible isometry (or perhaps $\sigma^2$, if $d$ is even and $\sigma$ permutes the $2$ neighbors with line $\ell$).  In what follows we explain some efficient construction of such lines for $L={\rm I}_n$.

\subsection{Some stable isotropic lines for ${\rm I}_n$}
\label{sec:stabisoBV}

Let us focus on the case $L={\rm I}_n$. 
Fix two integers $q$ and $k$ with $qk \leq n$, and consider an element $\sigma \in {\rm S}_n \subset {\rm O}({\rm I}_n)$ which is a product of $k$ cycles of length $q$ with disjoint supports $C_1,\dots,C_k \subset \{1,\dots,n\}$.
%To fix ideas we may assume $C_i = \{ j+ai \, |\, 1 \leq j \leq a-1\}$ and $\sigma(j+ai)=j+ai+1$ for $1 \leq j <a$ and $0 \leq i <b$. 
Set $C= \coprod_{j=1}^k C_j$.
Choose any odd prime $p \equiv 1 \bmod q$, and for such a $p$, choose an element $\omega \in (\Z/p)^\times$ of order $q$.
Choose any $x \in (\Z/p\Z)^n$ such that 
$$x_i \neq 0 \iff i \in C,\, \, \, \, \, {\rm and}\, \,\,\,\,x_{\sigma^{-1}(i)} = \omega x_i \, \, \, {\rm for\, \, all}\, \, i \in C.$$ There are exactly 
$(p-1)^k$ such elements $x$. Note that line $\ell = \Z/p \,x$ generated by $x$ is automatically in ${\rm C}_n(\Z/p)$ for $q \neq 2$, 
as we have $\sum_{s=0}^{q-1} \omega^{2s}=0$. In the case 
$q=2$ (so $\omega=-1$) we add the assumption that $\ell$ is isotropic. 
By construction, we have $\sigma(x) = \omega x$ an thus
$$\sigma \in {\rm O}(L; \ell)\, \, {\rm and}\, \, \nu(\sigma) = \omega.$$
So $\sigma$ is a visible isometry of the $p$-neighbor ${\rm N}_p(x)$. 
This construction is promising, but too restrictive in practice. Indeed, as we (must) have $x_i=0$ for 
$i \notin C$ the visible root system of ${\rm N}_p(x)$ necessarily contains ${\bf A}_{n-qk}$, which is too restrictive to find the most interesting lattices. 
We can circumvent this problem using the following trick. Consider an extra integer $d$ prime to $p$, and choose any $d$-isotropic vector $y \in {\rm C}_n(\Z/d)$ such that 
$$y_{\sigma(i)}=y_i \, \, \, \forall i \in C.$$
We clearly have $\sigma(y) = y$, so $\sigma \in {\rm O}(L; \Z/d\,y)$, and $\nu(\sigma) = 1$. 
By the chinese remainder theorem, there is a unique $z \in (\Z/pd)^n$ with $z \equiv x \bmod p$ and $z \equiv y \bmod d$.
Then $z$ is $pd$-isotropic and generates a line $\ell'$ such that
$$\sigma \in {\rm O}(L; \ell'),\, \, \, \nu(\sigma) \equiv \omega \bmod p\, \, \, {\rm and}\, \, \nu(\sigma)=1 \bmod d.$$
Note that for $q$ odd, $\sigma$ necessarily belongs to the subgroup $\sigma \in {\rm O}(L; \ell')^0$.\ps
Using such lines, we may even produce lattices with empty visible root systems ! Indeed, it is enough to ensure that we have
$y_i \not \equiv \pm y_j \bmod d$ for $i$ and $j$ not in the same $C_k$. 
\ps\ps
\subsection{An example : The Bacher-Venkov unimodular lattices with no root in dimension $28$}
\label{sec:examplebvdim28}
${}^{}$\indent
Let us consider the problem of finding representatives of ${\rm X}_{28}^\emptyset$, which is the most
difficult computation in \cite{bachervenkov}. For this computation, we will use the variant ${\rm BV}$ of an invariant used by Bacher-Venkov, which is defined in the companion paper~\cite{allche}.
Using either the arguments {\it loc. cit.} or King's results, 
we know that the mass of ${\rm X}_{28}^\emptyset$ is \scalebox{.8}{$17924389897/26202009600$}. The denominator of this mass factors as $2^{12} \cdot 3^{9} \cdot 5^2 \cdot13$. \ps

(a) This suggests to look first for lattices with a visible automorphism of order $13$. 
By enumerating isotropic lines of the form above for $q=13$, the prime $p=53$, $k=2$ and odd $d\leq 17$, we find for $d=17$ (so $pd=901$)
two unimodular lattices with no element of norm $\leq 2$ and a visible isometry of order $13$, after running only over $12$ isotropic lines ! Their masses turn out to be respectively \scalebox{.8}{$1/18341406720$} and \scalebox{.8}{$1/116480$}. For instance the first one is 
${\rm N}_{901}(x)$ with $x \in (\Z/901)^{28}$ defined by the following formula, with $t \equiv 16 \bmod 53$ (of order $13$):
\begin{equation*}
{\tiny
\begin{array}{ccc} 
x \bmod 53 \, \, & = & \, \, (1, t, t^2, t^3, t^4, t^5, t^6, t^7, t^8, t^9, t^{10}, t^{11}, t^{12}, 1, t, t^2, t^3, t^4, t^5, t^6, t^7, t^8, t^9, t^{10}, t^{11}, t^{12}, 0, 0) \\ 
x \bmod 17 \, \, & = & \, \, (1, 1, 1, 1, 1, 1, 1, 1, 1, 1, 1, 1, 1, 4, 4, 4, 4, 4, 4, 4, 4, 4, 4, 4, 4, 4, 6, 7)
\end{array}
}
\end{equation*}
This is one of the last lattices found in \cite{bachervenkov}, whereas this method finds it instantly.
The remaining mass is then {\small $17021999/24883200$}, whose denominator is $2^{12} \cdot 3^5 \cdot 5^2$. \ps

(b) Now we use $q=5$, with the prime $p=11$. For $k=5$, we find $6$ new lattices in ${\rm X}_{28}^\emptyset$ in the first thousand of lines, namely 
$3$ for $d=35$, with masses \scalebox{.8}{$1/400$}, \scalebox{.8}{$1/7680$} and \scalebox{.8}{$1/15360$}, and $3$ others for $d=37$, with masses \scalebox{.8}{$1/696729600$}, \scalebox{.8}{$1/7680$} and \scalebox{.8}{$1/3317760$}.
The remaining mass is \scalebox{.8}{$474647137/696729600$}, but other tries do not seem to find new lattices. On the other hand,  we have $696729600 = 2^{14} \cdot 3^5 \cdot 5^2 \cdot 7$, and the prime $7$ now appears.\footnote{Actually, $7$ already appeared actually in the two masses found in (a), but disappeared in the remaining mass there by an unlucky cancelation.} \ps

(c) So we try $q=7$, with the prime $p=43$. For $k=4$ and $d\leq 15$, we find $4$ new lattices in the first two hundreds of lines for $d=15$ (hence $pd=645$), with respective mass \scalebox{.8}{$1/344064$},  \scalebox{.8}{$1/112$} and \scalebox{.8}{$1/96768$} and \scalebox{.8}{$1/9676800$}.
The remaining mass is then \scalebox{.8}{$836459/1244160$}, with denominator $2^{10}\cdot 3^5 \cdot 5$ (no more $7$). \ps

(d) Trying $q=3$, $p=7$ and $k=8$ we immediately find $7$ lattices for $d=45$ with masses \scalebox{.8}{$1/96$}, \scalebox{.8}{$1/1728$}, \scalebox{.8}{$1/96$}, \scalebox{.8}{$1/15552$}, \scalebox{.8}{$1/55296$}, \scalebox{.8}{$1/6144$} and \scalebox{.8}{$1/192$}. The remaining mass is then \scalebox{.8}{$4957/7680$}, with denominator $2^9 \cdot 3 \cdot 5$ (luckily the exponent of $7$ is now $0$ and that of $5$ is $1$). For $k=9$ we find two more lattices for $d=41$, with mass \scalebox{.8}{$1/24$} and \scalebox{.8}{$1/192$}, hence the remaining mass \scalebox{.8}{$4597/7680$}, with denominator $2^9\cdot 3\cdot 5$. \ps

(e) At this point, we have actually found all the lattices with smallest mass, so the most efficient (and lazy) method is to run the general algorithm described in the introduction, biased with an empty visible root system. This way we do quickly find the $17$ remaining lattices as $d$-neighbors of ${\rm I}_{28}$ with $61 \leq d \leq 70$. 

%Trying $q=2$, $p=7$ and $k=12$ we quickly find $11$ lattices for $d=61$ and $65$. Nevertheless, at this point we have actually found all the lattices with smallest mass, so the most efficient is to run the general algorithm with empty visible root system, which does find all the remaining lattices as $d$-neighbors of ${\rm I}_{28}$ with $61 \leq d \leq 65$.

\subsection{A constraint : the {\it type} of a prime order isometry of a unimodular lattice }
\label{subsect:constr}
Let $L \in \mathcal{L}_n$ and $\gamma \in {\rm O}(L)$ of odd prime order $q$.
The characteristic polynomial of $\gamma$ has the form\footnote{Here $\Phi_d$ 
denotes the $d$-th cyclotomic polynomial.}
 $$\Phi_q^k \Phi_1^{l},\,\,\,\,{\rm  with}\,\,\,\,  n = l+k(q-1).$$ 
 As explained in~\cite[\S 4]{nebeaut}, the rank $l$ sublattice $L_1 = \ker (\gamma-1) \subset L$
satisfies $L_1^\sharp/L_1 \simeq (\Z/q)^s$ with $s \leq k,l$ and $s \equiv k \bmod 2$.
The pair $(k,s)$ is called the {\it type} of $\gamma$; it only depends on the $\Z_q[\gamma]$-module $L \otimes \Z_q$. \ps
In the case $L={\rm I}_n$, then $\sigma$ is necessarily a product of $k$ cycles of length $q$, 
and the orthogonal of $L_1$ in $L$ is clearly isometric to ${\rm A}_{q-1}^k$. So $\sigma$ is of type $(k,k)$. 
As a consequence, for any stable isotropic line $\ell \in {\rm C}_n(\Z/d)$
with $d$ prime to $q$, the visible isometry $\sigma$ of the associated $d$-neighbors will also have type $(k,k)$. This is an important restriction in the above method, although types $(k,k)$ 
seems to be the most common ones in practice. Another restriction, more obvious, is the fact that 
we must have $kq \leq n$, or equivalently $l \geq k$, instead of the most general case $k(q-1) \leq n$.
This excludes for instance, in the case $n \equiv 0 \bmod q-1$, the unimodular lattices defined by {\it Hermitian} $\Z[\zeta_q]$-lattices of rank $\frac{n}{q-1}$. 

\begin{example} {\rm
Going back to the Example of \S~\ref{sec:examplebvdim28}, there are $11$ elements in ${\rm X}_{28}^\emptyset$ having an isometry of order $5$.
Two of them, with masses \scalebox{.8}{$1/320$} and \scalebox{.8}{$1/160$}, have not been found in step (b) (nor in steps (a), (c) and (d) as we have $160=2^5 \cdot 5$). Indeed, using \cite{GAP} we can check that these two lattices have a single conjugacy class of order $5$ isometry, and whose type is $(6,4)$. 
But there is no isometry of ${\rm I}_{28}$ with characteristic polynomial $\Phi_5^6 \Phi_1^4$.
Actually, a third lattice also has a single conjugacy class of order $5$ isometry, which is of type is $(6,4)$. 
But this lattice has an order $13$ isometry of type $(2,2)$, and we found it in step (a).
}
\end{example}

\subsection{``{\it The neighbors of a lattice with small mass likely have a small mass}''}
\label{subsect:neismallmasshassmallmass}
${}^{}$ Here is an alternative method that we used in many instances to find lattices of small mass, and whose slogan is the title of this subsection. The ideas is simply that if we have $L_0 \in \mathcal{L}_n$ with a large isometry group, then ${\rm C}_{L_0}(\Z/d)$ will usually contain many points with a non-trivial stabilizer in ${\rm O}(L_0)$, hence leading to $d$-neighbors of $L_0$ with non-trivial visible isometry groups. 
 \par
In practice we often take $d=2$, and assume $L_0$ given as a $d_0$-neighbor of ${\rm I}_n$ with $d_0$ odd, say $L_0={\rm N}_{d_0}(x_0)$. So we expect to find lattices with non-trivial isometry groups of the form ${\rm N}_{2d_0}(x)$ with $x \equiv x_0 \bmod d_0$ by Lemma \ref{lem:voisdevois}. In practice we often combine this idea with that of the visible root system, by imposing that ${\rm N}_{d_0}(x)$ and ${\rm N}_{2d_0}(x)$ have the same visible root systems (we talk about "strict $2$-neighbors" of ${\rm N}_{d_0}(x)$). 
%Note that the constraints given in \S~\ref{subsect:constr} do not apply to $\gamma_0$, as $\gamma_0$ is {\it not} an isometry of ${\rm I}_n$ in general.

\begin{example} \label{ex:12A1X27} 
Consider for instance the problem of finding the elements $L \in {\rm X}_{27}$ 
with ${\rm r}_1(L)=0$ and ${\rm R}_2(L) \simeq 12\,{\bf A}_1$. 
From King's results, the total reduced mass of those lattices is \scalebox{.8}{$368401/138240$}. 
Using the visible 
root system $9 {\bf A}_1$, and an already lengthy enumeration of the corresponding $d$-neighbors for odd $d$ from $d=37$ to $d=45$, we find $26$ such lattices, with remaining reduced mass \scalebox{.8}{$731/276480$}. The one with smallest mass is ${\rm N}_{45}(x)$ with
$$x=(\scalebox{.7}{{\color{magenta} 1, 1}, 2, {\color{magenta} 4, 4}, {\color{magenta} 6, 6}, 7, {\color{magenta} 9, 9}, 10, 11, 12, {\color{magenta}13, 13}, 14, {\color{magenta}16, 16}, {\color{magenta}17, 17}, {\color{magenta}18, 18}, 19, 20, {\color{magenta}21, 21}, 22}) \in \Z^{27}$$
whose reduced mass is \scalebox{.8}{$1/2048$}. By considering ``solely'' the $\simeq 2^{16}$ strict $2$-neighbors of this lattice we quickly find the three remaining ones,
with respective reduced mass \scalebox{.8}{$1/384$}, \scalebox{.8}{$1/46080$} and \scalebox{.8}{$1/55296$}:
 see Table~\ref{tab:12A1dim26}. 
\end{example}

\tabcolsep=4pt
\begin{table}[H]
{\scriptsize 
\renewcommand{\arraystretch}{1.8} \medskip
\begin{center}
\begin{tabular}{c|c|c|c|c|c|c|c|c|c|c|c|c}
$\texttt{reduced mass}$ & \scalebox{.8}{$1/4$} & \scalebox{.8}{$1/8$} & \scalebox{.8}{$1/16$} & \scalebox{.8}{$1/32$} & \scalebox{.8}{$1/48$} & \scalebox{.8}{$1/64$} & \scalebox{.8}{$1/384$} & \scalebox{.8}{$1/640$}  & \scalebox{.8}{$1/720$} & \scalebox{.8}{$1/2048$} & \scalebox{.8}{$1/46080$} & \scalebox{.8}{$1/55296$} \\ \hline
$\sharp {\texttt{lattices}}$ & $8$ & $2$ & $4$ &  $3$ & $2$ & $1$ & $4$  & $1$  & $1$ & $1$ & $1$ & $1$
\end{tabular} 
\end{center}
}
\caption{{\small The $29$ lattices in ${\rm X}_{27}$ with no norm $1$ elements and root system $12 {\bf A}_1$.}}
\label{tab:12A1dim26}
\end{table}

\section{An example : the lattice ${\rm N}_{2n+1}(1,2,\dots,n)$}
\label{sec:ex12n}

In this section we give an example of study of the non-visible part of a neighbor.
More precisely, we fix an integer $n\geq 1$ with $n \not \equiv 1 \bmod 3$ and
set $x=(1,2,\dots,n)$. As already explained in the introduction Formula~\eqref{eq:oddshortborchleech}, 
we have a unimodular lattice 
$$ {\rm N}_n:={\rm N}_{2n+1}(x) \in \mathcal{L}_n$$
defined by ${\rm N}_n={\rm M}_n + \Z \frac{1}{2n+1} x'$ with ${\rm M}_n:={\rm M}_{2n+1}(x)$ 
and $x'=x+\frac{(2n+1)n^2(n+1)}{6} \epsilon_1$.  As we have $\pm i, \pm i \pm j \not \equiv 0 \bmod 2n+1$ for $1 \leq i < j \leq n$,
we have no visible element of norm $\leq 2$, {\it i.e.} ${\rm R}_2({\rm M}_n)=\emptyset$.
Better :

\begin{prop} 
\label{prop:R2Nn} We have ${\rm r}_2({\rm N}_n)= 0$ for all $n\geq 23$. 
\end{prop}

\begin{pf} 
Fix  $z\in {\rm N}_n \smallsetminus {\rm M}_n$. 
For some divisor $b$ of $2n+1$, we may write $z=m + \frac{k}{b} x'$ with $m \in M$ and $1 \leq k<b$ 
coprime with $b$.  Write $$b=2s+1\, \, \, {\rm and}\, \, \, 2n+1=(2t+1)(2s+1).$$ We have $n=(2s+1)t+s$, so
the coordinates of $x$ are $\equiv$ modulo $b$ to $\pm 1, \pm 2, \dots, \pm s, 0$ ($t$ times) and then to $1,2,\dots,s$, modulo $b$. Observe that if $S$ is a subset of $X:=\Z/(2s+1)\smallsetminus \{0\}$ satisfying
$X = S \coprod -S$, then the same holds for $kS$ for all $k \in (\Z/(2s+1))^\times$.
It follows that the coordinates mod $b$ of the element $bz = bm+ k x'$, also are $2t+1$ times $\pm 1, \pm 2, \dots, \pm s$, and $t$ times $0$. A trivial coordinatewise inequality shows then
%Applying coordinatewise the trivial inequality of Lemma~\ref{lem:ineqelem}, 
\begin{equation}
\label{ineq:zz} z.z\,\,\, \geq \,\,\,\frac{1}{b^2} (2t+1)\sum_{i=1}^s i^2 \,\,\,=\,\,\,  \frac{2n+1}{24}\left(1 - \frac{1}{(2s+1)^2}\right).
\end{equation}
As $2n+1 \not \equiv 0 \bmod 3$ we have $s\geq 2$, and thus $z.z  \geq \frac{53}{24}\frac{24}{25}>2$ for $n \geq 26$. 
The two remaining cases $n=23$ and $24$ (Short Leech and Odd Leech lattices) could be further analyzed in this style (or checked with a computer), but they are classical so we omit them.\end{pf}

%could be analyzed further 
%or easily checked with the computer, or we can argue as follows.
%For $n=23$, we have $2n+1=47$ prime, hence $t=0$ and $s=n=23$, and $z.z \geq \frac{92}{47}=2-2/47$. 
%Assuming $z.z=2$ the inequality~\eqref{ineq:zz} is strict, but this contradicts the inequality case of Lemma~\ref{lem:ineqelem} below.
%{\color{blue} For $n=24$ we have $2n+1=49$, and $z.z=2$ implies $s=3=t$. We are in the equality case of~\eqref{ineq:zz},  
%so the coordinates of $bz$ are, up to sign, each of $1, 2, 3$ with multiplicity $7$, and $3$ with multiplicity $0$.}

%We have used the following elementary lemma.

%\begin{lemma}
%\label{lem:ineqelem}
%Let $x \in \Z$, $s\in \mathbb{N}$ and $k$ the unique integer in $[-s,s]$ with $x \equiv k \bmod 2s+1$. 
%We have $x^2 \geq k^2$.
%, and if the inequality is strict then we have 
%$x^2 \geq k^2+(2s+1)(2s-2|k|+1)$ with equality if $s=|k|$.
%\end{lemma}
% 
%We now study ${\rm R}_3({\rm N}_n)$.

\begin{prop} 
\label{prop:R3Nn}
For $n\geq 5$ the lattice ${\rm M}_n$ is generated by ${\rm R}_3({\rm M}_n)$, and for $n\geq 36$ 
we have ${\rm R}_3({\rm M}_n)={\rm R}_3({\rm N}_n)$.
\end{prop}

\begin{pf} Let ${\rm L}_n$ be the orthogonal of $x=(1,2,\dots,n)$ in ${\rm I}_n$. 
We have ${\rm L}_n\subset {\rm M}_n$.
A simple computation shows that ${\rm L}_5$ is generated by its norm $3$ vectors (there are $8$ such vectors).
Using ${\rm L}_{n} = {\rm L}_{n-1} \times \{0\} + \Z e$ with $e=\scalebox{0.9}{(-1,0,\dots,0,-1,1)}$, the same holds for ${\rm L}_n$ with $n\geq 5$. The vector $e'=\scalebox{.9}{(0,1,0,\dots,0,1,1)}$ satisfies $e'.e'=3$ and $e'.x=2n+1$, 
so we have ${\rm M}_n={\rm L}_n+ \Z e'$ and the first assertion holds.\par
A non-visible vector $z$ in ${\rm N}_{2n+1}(1,2,\dots,n)$ satisfies the inequality \eqref{ineq:zz}, 
for some integer $s\geq 2$ such that $2s+1$ divides $2n+1$.
For $n> 37$, the right hand side of \eqref{ineq:zz} is $\geq \frac{77}{24}\frac{24}{25}>3$. 
For $n=36$, then $2n+1=73$ is prime, so we have $s=n$ and $z.z \geq \frac{222}{73} >3$.
\end{pf}

We finally determine the isometry group of ${\rm N}_n$.
Recall from Sect.~\ref{sect:visisometry} the visible subgroup ${\rm O}({\rm N}_n)^{\rm v}={\rm O}({\rm N}_n) \cap {\rm O}({\rm I}_n)$ 
as well as the morphism $$\nu : {\rm O}({\rm N}_n)^{\rm v} \rightarrow (\Z/(2n+1))^\times.$$ 

\begin{prop} 
\label{prop:ONn}
For all $n\geq 1$, the morphism $\nu$ defines an isomorphism ${\rm O}({\rm N}_n)^{\rm v} \isomo (\Z/(2n+1)\Z)^\times$.
Moreover, for $n\geq 32$ we have $ {\rm O}({\rm N}_n)^{\rm v}={\rm O}({\rm N}_n)$.
\end{prop}

\begin{pf} 
The first assertion is a consequence of Prop.~\ref{prop:visautOLl}.
Indeed, in the notations of this proposition, and given the shape of $x$, we have ${\rm m}(\ell)={\rm m}'(\ell)=0$ and ${\rm a}_i(\ell)=1$ for all $i \in I$, so $\nu$ is injective, as well as 
${\rm X}(\ell)=(\Z/(2n+1)\smallsetminus \{0\})/\{\pm 1\}$, so we have ${\rm H}(\ell)=(\Z/(2n+1))^\times$ and $\nu$ is surjective.\ps
For $n\geq 36$, we have ${\rm O}({\rm N}_n) \subset {\rm O}({\rm M}_n)$ by Proposition~\ref{prop:R3Nn}, hence the last assertion by Remark~\ref{rem:othervisdef}.
In the remaining cases $n=32, 33, 35$, it is enough to check $|{\rm O}({\rm N}_n)|=\varphi(2n+1)$, which
follows from a computation using the Plesken-Souvignier algorithm 
(actually, we have ${\rm R}_3({\rm M}_n) \subsetneq {\rm R}_3({\rm N}_n)$ for those $n$).
\end{pf}

\begin{remark} For $n=23, 24, 26, 27, 29$ and $30$, then $|{\rm O}(N_n)|/\varphi(2n+1)$  is respectively $1\,839\,366\,144\,000$, $23\,876\,075\,520$, $360\,000$, $192$, $4$ and $2$.
\end{remark}

\section{Exceptional lattices and visible characteristic vectors}
\label{sec:excsec}
%move in introduction ?

	As already observed by Bacher and Venkov in their study of lattices with no root \cite{bachervenkov}, certain unimodular lattices that they term {\it exceptional} tend to be harder to find. For this same reason they played an important role in our search, which explains this section.

\subsection{Exceptional lattices}
\label{subsect:exclat}
 Let $L$ be a unimodular lattice of rank $n$. Recall that ${\rm Char}(L)$ denotes the set of characteristic vectors of $L$ (see \S~\ref{subsect:paritynei}). We have
\begin{equation}
\label{eq:normcharvec}
\forall \xi \in {\rm Char}(L), \, \, \xi.\xi \equiv n \bmod 8.
\end{equation}
Indeed, if $L$ is even this holds as $\xi \in 2L$, and if $L$ is odd it also holds as 
$L$ and ${\rm I}_n$ are isometric over $\Z_2$. The following definition is a generalization of a terminology of Bacher and Venkov~\cite[\S 3]{bachervenkov}.

%Alternatively, for $v \in L$ we have $v \in {\rm Char}(L)$ if, and only if, the lattice $\{ x \in L\, |\, x \cdot v \equiv 0 \bmod 2\}$ is even. 

\begin{definition} 
A unimodular lattice $L$ of rank $n$ is called {\rm exceptional} if there exists $\xi \in {\rm Char}(L)$ with $\xi \cdot \xi <8$.
\end{definition}

As an example, it follows from Formula \eqref{eq:charIn} that the lattice ${\rm I}_n$ is exceptional only for $n<8$, in which case it has exactly $2^n$ characteristic vectors of norm $<8$. Observe that for $L= A \perp B$, we have ${\rm Char}(L)=\{a +b \,\,|\,\, a \in {\rm Char}(A), b \in {\rm Char}(B)\}$. We deduce the following fact, where ${\rm \rho}(n)$ denotes the unique integer $r$ satisfying $n \equiv r \bmod 8$.

\begin{prop} 
\label{prop:excwithI1}
Let $m,n\geq 0$ and $L$ be a unimodular lattice of rank $n$. Then ${\rm I}_m \perp L$ is exceptional if, and only if, 
$L$ is exceptional and $m+{\rm \rho}(n)<8$.
\end{prop}

The exceptional unimodular lattices of rank $n\equiv 0 \bmod 8$ are just the even unimodular lattices. Also, there is clearly no exceptional unimodular lattice of rank $n$ with ${\rm \rho}(n)=1$ and ${\rm r}_1(L)=0$. 
We will say more about the cases ${\rm \rho}(n)>1$ below. 
For $L \in \mathcal{L}_n$ we denote by ${\rm Exc}(L)$ the set of $e \in {\rm Char}(L)$ with $e.e={\rm \rho}(n)$.
Proposition~\ref{prop:sizeExc} below gives information on ${\rm Exc}(L)$ when $2 \leq \rho(n) \leq 4$. 
In order to prove it we recall a classical specificity of the case $n \equiv 4 \bmod 8$. 

\begin{prop}
\label{prop:defcompa}
Assume $L \in \mathcal{L}_n$ with $n\equiv 4 \bmod 8$, and let.
$M$ denote the largest even lattice in $L$ (the "even part" of $L$). 
%As is well-known, and already seen during the proof of Proposition~\ref{prop:sizeExc}, 
there are exactly two other unimodular lattices in $\mathcal{L}_n$ with same even part $M$. 
\end{prop}

Indeed, this immediately follows from ${\rm res}\, M \simeq \,{\rm res}\, {\rm D}_4$ and Table~\ref{tab:venkov}. The two other lattices of the statement will be called the {\it companions} of $L$.

%This is a finite set equipped with a natural action of ${\rm O}(L)$. 

\begin{prop} 
\label{prop:sizeExc} Assume $L \in \mathcal{L}_n$ is exceptional with ${\rm r}_1(L)=0$. 
\begin{itemize}
\item[(i)] For ${\rm \rho}(n) \in \{2,3\}$, or ${\rm \rho}(n)=4$ and ${\rm r}_2(L)=0$, we have $|{\rm Exc}(L)|=2$.
\item[(ii)] For ${\rm \rho}(n)=4$ we have $|{\rm Exc}(L)|\leq 2n$.
\end{itemize}
\end{prop}

The ideas in the proof below are inspired from the proof of Prop. 4.1 in \cite{bachervenkov}, 
which contains the special case ${\rm r}_2(L)=0$.

\begin{pf} Assume $L \in \mathcal{L}_n$ with ${\rm \rho}(n)=4$ 
(but do not assume ${\rm r}_1(L)=0$ for the moment).
Define $M$ as even part of $L$, namely $M={\rm M}_2(L; \xi)$ for any $\xi \in {\rm Char}(L)$,
and denote by $L'$ and $L''$ the two companions of $L$.
Then we have $M^\sharp = L \cup L' \cup L''$, and the map $M^\sharp \rightarrow L, v \mapsto 2v,$ defines bijections 
\begin{equation}
\label{eq:bijexc4mod8}
{\rm R}_1(M^\sharp) \isomo {\rm R}_1(L) \coprod {\rm Exc}(L)\, \, \, {\rm and}\, \, {\rm R}_1(L') \coprod {\rm R}_1(L'') \isomo {\rm Exc}(L).
\end{equation}
Assume furthermore $L$ exceptional. Up to renaming $L'$ and $L''$, we may thus assume ${\rm r}_1(L') \neq 0$. But if any of $L, L'$ and $L''$ has a norm $1$ vector, then the orthogonal symmetry about this vector defines an isometry between the two other lattices (see~\cite[Prop. 2.3]{bachervenkov}). This shows $L \simeq L''$, and also that we have $L' \simeq L''$ in the case $L$ has a norm $1$ vector. \par
	Assume first $L = {\rm I}_r \perp L_0$ with $r \in \{1,2\}$ and ${\rm r}_1(L_0) = 0$. Then we have $L \simeq L' \simeq L''$ and thus $|{\rm Exc}(L)|=2r+2r=4r$. But the discussion before Prop.~\ref{prop:excwithI1} shows 
	$|{\rm Exc}(L)|=2^r \cdot |{\rm Exc}(L_0)|$. This proves $|{\rm Exc}(L_0)|=4r/2^r=2$, hence part (i) (for the lattice $L_0$) in the case $2 \leq {\rm \rho}(n) \leq 3$. \par
	Assume finally ${\rm r}_1(L)=0$. Then we have ${\rm r}_1(L'')=0$ and $L' \simeq {\rm I}_k \perp N$ 
for some $1 \leq k \leq n$ and ${\rm r}_1(N) = 0$. We have $|{\rm Exc}(L)|=2k$. Note that the norm $2$ vectors of $L$, $L'$ and $M$ are the same. This shows $k=1$ for ${\rm r}_2(L)=0$, and proves assertions (ii) and (i) for ${\rm \rho}(n)=4$.
\end{pf}

\begin{remark}
\label{rem:casen4}
 {\rm (The case ${\rm \rho}(n)=4$)} Let $n\geq 1$ be an integer $\equiv 4 \bmod 8$. 
Let $\mathcal{A}_n$ be the groupoid of exceptional lattices $A \in \mathcal{L}_n$ with ${\rm r}_1(A)=0$, 
and let $\mathcal{B}_n$ be the groupoid of non-exceptional lattices $B \in \mathcal{L}_n$ with ${\rm r}_1(B) \neq 0$.
The last paragraph in the proof above shows that each $L$ in $\mathcal{A}_n$ has a unique companion $L' \in \mathcal{L}_n$ with ${\rm r}_1(L') \neq 0$, necessarily non-exceptional. Better, $L \mapsto L'$ defines a natural functor $\mathcal{A}_n \rightarrow \mathcal{B}_n$ which induces a bijection on isometry classes on both sides, and satisfies the following properties : 
$|{\rm Exc}(L)|={\rm r}_1(L')$, ${\rm R}_2(L) = {\rm R}_2(L')$, and ${\rm O}(L)$ is an index $2$ subgroup of ${\rm O}(L')$.
\end{remark}

\subsection{Mass formulae for $2 \leq {\rm \rho}(n) \leq 7$.}

In the cases $2 \leq {\rm \rho}(n) \leq 7$, the exceptional unimodular lattices are related to the genus\footnote{These lattices form a single genus by \cite[Chap. 15 p. 387]{conwaysloane}. An example is given by the root lattice ${\rm A}_{r-1} \perp {\rm E}_8^{\frac{n-r}{8}}$, with $r={\rm \rho}(n)$. For all such lattices $L$, we have in particular ${\rm res}\, L \simeq {\rm res}\, {\rm A}_{r-1}$ (a cyclic group of order $r$).}
$\mathcal{G}_{n}$ of even Euclidean lattices of rank $n-1$ and determinant ${\rm \rho}(n)$. 
This is presumably quite classical, but we briefly recall how in this section, and use this to derive a few interesting mass formulae. \ps

  For $L \in \mathcal{L}_n$ and $e \in {\rm Exc}(L)$ we denote by 
$L(e)$ the orthogonal of $e$ in $L$. This is an even lattice as $e \in {\rm Char}(L)$, with determinant $e.e={\rm \rho}(n)$ as $L$ is unimodular, so we have $L(e) \in \mathcal{G}_n$.
Also, $\Z e$ is saturated in $L$.
Set $\delta(n)=1$ if ${\rm \rho}(n)=2$, and $2$ otherwise.

\begin{lemma} 
\label{lem:glueexc}
Assume $2\leq {\rm \rho}(n) \leq 7$, let $N \in \mathcal{G}_n$ and set $M = N \perp \Z e$ with $e.e={\rm \rho}(n)$. There are $\delta(n)$ integral unimodular lattices $L$ containing $M$ and in which $\Z e$ is satured, and they are permuted transitively by $1 \times {\rm O}(\Z e) \simeq \Z/2$.
\end{lemma}

\begin{pf} The finite bilinear space $V:={\rm Res}\, \Z e$ is isometric to $\Z/r$ equipped with the $\Q/\Z$-valued bilinear form $(i,j) \mapsto \frac{ij}{r} \bmod \Z$, and is isometric to $-{\rm res}\, N$. By Sect.~\ref{notations} (iii), the unimodular lattices $L$ of the statement naturally correspond to the bilinear Lagrangians in ${\rm res}\, M \simeq -V \perp V$ which are transversal to $0 \perp V$. We conclude as $V$ is cyclic of order $r={\rm \rho}(n)$ and the only solutions to $x^2 = 1$ in $\Z/r$ are the $\delta(n)$ elements $x=\pm 1$ for $2 \leq r \leq 7$.
\end{pf} 

For $L \in \mathcal{L}_n$ the group ${\rm O}(L)$ naturally acts on ${\rm Exc}(L)$, as well as on the quotient
set ${\rm Exc}(L)^{\pm }:={\rm Exc}(L)/\{ \pm {\rm id}_L\}$.
Let $\mathcal{E}_n$ be the natural groupoid whose objects are the pairs $(L, \pm e)$ with $L \in \mathcal{L}_n$ and $\pm e \in {\rm Exc}(L)^{\pm}$. (For ${\rm \rho}(n)=2,3$, Proposition~\ref{prop:sizeExc} shows that the datum of $\pm e$ is unique, hence superfluous). 
Lemma~\ref{lem:glueexc} shows:

\begin{prop}
\label{prop:foncEG} The natural functor $F: \mathcal{E}_n \rightarrow \mathcal{G}_{n}, \,(L, \pm e) \mapsto L(e),$ is essentially surjective, and for all $(L, \pm e)$ in $\mathcal{E}_n$ the natural morphism ${\rm O}(L, \pm e) \rightarrow {\rm O}(L(e))$ is surjective and $\delta(n)$ to $1$. 
\end{prop}

For an (isomorphism class of) ADE root system $R$, let us denote by $\mathcal{E}^R_n$ the subgroupoid of pairs $(L,e)$ in $\mathcal{E}_n$ with ${\rm R}_2(L(e)) \simeq R$.

\begin{cor} 
\label{cor:massformexc}
Assume $2 \leq {\rm \rho}(n) \leq 7$. Let $L_1,\dots,L_h$ be representatives for the isometry classes of exceptional unimodular lattices of rank $n$. 
\begin{itemize} 
\item[(i)] We have $\frac{1}{\delta(n)}\, {\rm mass}\,\mathcal{G}_{n} \,\,={\rm mass}\, \mathcal{E}_n\,\,=\,\, \sum_{i=1}^h \,\, |{\rm Exc}(L_i)^{\pm}|\,{\rm mass}(L_i)$.\ps
\item[(ii)] For each {\rm ADE} root system $R$, we have $\frac{1}{\delta(n) }{\rm mass}\,\mathcal{G}_{n}(R) \,\,= \,\,{\rm mass} \,\mathcal{E}^R_n$.
\end{itemize}
\end{cor}

\begin{pf} The first equality in (i), as well as (ii), follow from Proposition~\ref{prop:foncEG}. The second equality in (i) follows from the obvious equivalence of groupoids $\mathcal{E}_n \,\,\simeq \,\,\coprod_{i=1}^h\, [{\rm Exc}(L_i)^{\pm}/{\rm O}(L_i)]$ (see~Sect.~\ref{sect:mass} for several similar arguments).
\end{pf}

\begin{remark} 
\label{rem:rootmassesGnR}
{\rm 
(Determination of the terms ${\rm mass}\, \mathcal{G}_n(R)$)
{\it
Assume ${\rm \rho}(n)=2, 3, 4, 5$ or $6$, and set
$R_0={\bf E}_7, {\bf E}_6, {\bf D}_5, {\bf A}_4$ or ${\bf A}_1 {\bf A}_2$ 
respectively. Then we have ${\rm res}\, {\rm A}_{{\rm \rho}(n)-1} \simeq - {\rm res}\, R_0$. In particular, $\mathcal{G}_n$ coincides with the genus $\mathcal{G}$ in the statement of Proposition~\ref{kingform2}, with $m:=9-{\rm \rho}(n)$. 
Assuming furthermore $n\leq 30$, we may thus deduce ${\rm mass}\, \mathcal{G}_n(R)$ from that proposition and {\rm \cite{king}}.
}
}
\end{remark}

\begin{example} 
\label{ex:G26G27}
{\it {\rm (The cases $n=26$ and $27$)} For such an $n$ and an exceptional lattice $L \in \mathcal{L}_n$ with ${\rm r}_1(L)=0$,
we have $|{\rm Exc}(L)|=2$ by Prop.~\ref{prop:sizeExc}.
\begin{itemize}
\item[(i)] The genus $\mathcal{G}_{26}$ of even lattices of rank $25$ and determinant $2$ has been determined by Borcherds in {\rm \cite{borcherds}}. There are $121$ such lattices: $24$ of the form ${\rm A}_1 \perp N$ with $N$ a Niemeier lattice, and $97$ with a dual lattice with minimum $>1/2$. Accordingly, the corresponding $26$ dimensional exceptional lattice $L$ is either ${\rm I}_1 \perp N$, or satisfies 
${\rm r}_1(L)=0$ and ${\rm R}_2(L(e)) \coprod \{\pm e\}={\rm R}_2(L)$. 
\par
\item[(ii)] The genus $\mathcal{G}_{27}$ of even lattices of rank $26$ and determinant $3$
has also been determined by Borcherds {\rm \cite{borcherdsthesis}}, up to a few indeterminacies that were settled by M\'egarban\'e in {\rm \cite{megarbane}}. They are $678$ such lattices. Also, for $(L,\pm e)$ in $\mathcal{E}_{26}$, we easily check the equivalences between\footnote{For $v \in L(e)^\sharp$ with $v.v=2/3$, then $\pm v \pm e/3$ has norm $1$ and $\pm v \pm 2e/3$ has norm $2$.}: {\rm (a)} ${\rm r}_1(L)= 0$, 
{\rm (b)} there is no $v \in L(e)^\sharp$ with $v.v=2/3$,  {\rm (c)} ${\rm R}_2(L(e)) = {\rm R}_2(L)$.
\par
\item[(iii)] The genus $\mathcal{G}_{28}$ of even lattices of rank $27$ and determinant $4$
is also denoted $\mathcal{D}_{28}$ in Sect.~\ref{subsect:kingsmassform} and thus is easily deduced from $\mathcal{L}_{27}$. See also Remark~\ref{rem:casen4} for yet another approach. \end{itemize}
}
\end{example}

%\footnote{Recall also the following result due to Elkies \cite{elkies}:  if $L \in \mathcal{L}_n$ is such that $\xi.\xi \geq n$ for all
%$\xi \in {\rm Char}(L)$, then $L \simeq {\rm I}_n$.}

%As an example, consider the case $R_0 \simeq {\bf A}_m$ for $m\geq 1$. Then 
%${\rm res} {\rm A}_m$ is isomorphic to $\Z/(m+1)$ with quadratic form $x \mapsto \frac{-m}{2(m+1)} x^2$

\subsection{Visible exceptional characteristic vectors}
\label{subsec:vischarvec}
	Our main aim now is to explain how to produce $d$-neighbors of ${\rm I}_n$ which are exceptional. 
		
\begin{prop}
\label{prop:viscar}
 {\rm (Visible characteristic vectors)}
Let $L$ be an odd unimodular lattice of rank $n$,
$d$ be an even integer, $x \in L$ be a $d$-isotropic vector, $\epsilon \in \{0,1\}$, and $N:={\rm N}_d(L; x; \epsilon)$ the associated $d$-neighbor of $L$. 
An element $e \in L$ is in ${\rm Char}(N)$ if, and only if, the following properties are satisfied:
\begin{itemize}
\item[(i)] $x .e \equiv 0 \bmod d$,\ps
\item[(ii)] either $e \in {\rm Char}(L)$, or $\exists \xi \in {\rm Char}(L)$ such that $x\,\, \equiv \,\,\xi \,-\, e\,\, \bmod 2L$,\ps
\item[(iii)] $\frac{x .e}{d} \,\,\equiv \,\, \frac{x.x}{2d} \,+\,  \epsilon (1+\frac{d}{2} ) \bmod 2$.
\end{itemize}
\end{prop}	

% check : this (ii) exaclty the condition checked in algo par_seek_except_refined
	
\begin{pf} Set $M=L \cap N$.
By definition, we have $e \in M$ if and only if (i) holds.
We first check that for $e \in M$, condition (ii) is equivalent to $m.m \equiv e.m \bmod 2$ for all $m \in M$. 
Fix $\xi \in {\rm Char}(L)$. For $m \in M$ we have $m.m \equiv \xi.m \bmod 2$.
We may thus assume $e \notin {\rm Char}(L)$. The condition $(\xi-e).m \equiv 0 \bmod 2$
for all $m \in M$ amounts to ask that the kernel $M={\rm M}_d(x)$ of the linear form $L \rightarrow \Z/d, v \mapsto v \cdot x \bmod d$, is included in the kernel $H$ of the nonzero linear form $L \rightarrow \Z/2, v \mapsto  v . (\xi - e) \bmod 2$. 
As $L/M$ is cyclic, the unique index $2$ lattice of $L$ containing $M$ is ${\rm M}_2(x)$, so we have $H={\rm M}_2(x)$, so that condition is equivalent to $x \equiv \xi -e \bmod 2$. \par
By definition of $N$, we have $N=M+ \Z\frac{\tilde{x}}{d}$ with $\tilde{x}=x+\,dr\,y$, $r = -\frac{x.x}{2d}+ \epsilon d/2 \in \Z$ and $y \in L$ with $y.x \equiv 1 \bmod d$ (see Remark~\ref{remdefline}). It only remains to check that condition (iii) is equivalent to $e.\frac{\tilde{x}}{d} \equiv \frac{\tilde{x}}{d}.\frac{\tilde{x}}{d} \bmod 2$. We have
$$\tilde{x}.\tilde{x}\,\, \equiv\,\, x.x \,+\,2\,dr\,x.y \,+\,d^2 r^2 \,y.y\,\, \equiv\,\,\epsilon d^2 \,+\, d^2 r \,y.y \,\,\bmod \,2d^2,$$
$$d \,\tilde{x}. e  \,\,\equiv  \,\,d\, x.e\, +\, d^2 r\, y.e \,\,\bmod 2d^2.$$ 
Using $y.y - y.e \equiv (\xi - e).y \equiv x.y \equiv 1 \bmod 2$ we obtain
$$\tilde{x}.\tilde{x}\,\,-\,\, d \,\tilde{x}.e \,\,\equiv \,\, \epsilon d^2\,+\,d^2r\,+\, d\,x.e\,\,\equiv \,\,\epsilon d^2 \,(1+d/2) \,- \, d\, \frac{x.x}{2}\, +\, d \,x.e\,\, \bmod 2d^2,$$
and conclude the proof.
\end{pf}	
		
\begin{remark} Part (iii) implies $2 x.e \equiv x.x \bmod 8$, and then $(x-2e).x \equiv 0 \bmod 8$. If we write $x= e + \xi$ with $\xi \in {\rm Char}(L)$, we deduce the congruence $e.e \equiv \xi . \xi \bmod 8$, in agreement with Formula~\eqref{eq:normcharvec}. 
\end{remark}

\subsection{The case $L={\rm I}_n$}
\label{subsec:vischarvecIn}
	
We now discuss the special case $L={\rm I}_n$.
Fix an integer $1 \leq r < n$ with $n \equiv r \bmod 8$ and consider the element $e \in \Z^n$ defined by
% $e_i=0$ for $i \leq n-r$, and $e_i=1$ for $n-r< i \leq n$. We have 
$$e = (\underbrace{0,\dots,0}_{n-r},\underbrace{1,\dots,1}_{r}), \, \, {\rm with}\, \, \, e.e=r.$$
Fix an even integer $d$, a $d$-isotropic $x \in \Z^n$ and $\epsilon \in \{0,1\}$. 
Then conditions (i), (ii) and (iii) on $(x,\epsilon)$ in Proposition~\ref{prop:viscar} take the following forms :\begin{itemize}{\it
\item[{\rm (i)}] $\sum_{n-r<i \leq n} x_i \equiv 0 \bmod d$, \ps
\item[{\rm (ii)}] $x_i$ is odd for $i \leq n-r$, and $x_i$ is even for $i>n-r$,  \ps
\item[{\rm (iii)}] $2 \sum_{n-r<i \leq n} x_i \,\equiv\, \sum_{i=1}^n x_i^2 + d \epsilon (2+d) \bmod 4d$.}
\end{itemize}
(For (ii), just use that $\xi=(1,1,\dots,1)$ is in ${\rm Char}({\rm I}_n)$).
These equations do have many solutions $(x,\epsilon)$ for sufficiently large $d$, and for such an $(x,\epsilon)$ the associated $d$-neighbor ${\rm N}_d(x; \epsilon)$ has the 
concrete vector $e$ as a characteristic vector of norm $r$. The special case $r<8$ do lead to constructions of exceptional lattices.
%, and the pair $({\rm N}_d(x; \epsilon), \pm e)$ belongs to the groupoid $\mathcal{E}_n$ above. 
Note that these choices of $x$ also allow to prescribe the visible root system, as well as visible isometries, to some extent. If we furthermore impose that the $x_i$ are distinct mod $d$ for $n-r<i$, 
the visible root system of ${\rm N}_d(x; \epsilon)$ is a subroot system of $e^\perp$. \ps

\begin{example}
\label{ex:exc10A1dim26}
{\it
In the introduction, we considered the problem of finding all unimodular lattices of rank $26$ with root system $10 \,{\bf A}_1$. Among the $7$ such lattices, a single one is exceptional, namely the last, and most complicated to find there!, of Table~\ref{tab:10A1dim26}. This lattice $L$ is straightforward to find by the method above. For instance, using the visible root system $2{\bf A}_1\,{\bf D}_2 \,\simeq\, 4\,{\bf A}_1$ we immediately find 
$L \simeq {\rm N}_{92}(x)$ with $x = $ {\small $({\color{magenta} 1, 1},  3, 5, 7, 9, 11, 13, 15, 17, 19, 21, {\color{magenta} 23, 23}, \\ 25, 27, 29, 31, 33, 35, 37, 39, 41, 43, {\color{blue} 46, 46})$}.
}
\end{example}

%[100, [1; 1; 3; 7; 9; 11; 13; 15; 17; 19; 21; 23; 25; 27; 29; 31; 33; 35; 37; 39; 41; 43; 45; 47; 50; 50], 0, 1/92897280, "10A1", ["fHBV", 12421112404259529000]]
%[100, [1; 3; 7; 9; 11; 13; 15; 17; 19; 21; 23; 25; 27; 29; 31; 33; 35; 37; 39; 41; 43; 45; 47; 49; 50; 50], 0, 1/92897280, "10A1", ["fHBV", 12421112404259529000]]
%[92, [1; 1; 3; 5; 7; 9; 11; 13; 15; 17; 19; 21; 23; 23; 25; 27; 29; 31; 33; 35; 37; 39; 41; 43; 46; 46], 0, 1/92897280, "10A1", ["fHBV", 12421112404259529000]]
%[92, [1; 1; 5; 7; 9; 11; 13; 15; 17; 19; 21; 23; 23; 25; 27; 29; 31; 33; 35; 37; 39; 41; 43; 43; 46; 46], 0, 1/92897280, "10A1", ["fHBV", 12421112404259529000]]
%[92, [1; 1; 5; 9; 11; 13; 15; 17; 19; 21; 23; 23; 25; 27; 29; 31; 33; 35; 37; 39; 39; 41; 43; 43; 46; 46], 0, 1/92897280, "10A1", ["fHBV", 12421112404259529000]]
%[92, [1; 1; 7; 9; 11; 15; 17; 19; 21; 23; 23; 25; 27; 29; 31; 33; 33; 35; 37; 39; 41; 41; 43; 43; 46; 46], 0, 1/92897280, "10A1", ["fHBV", 12421112404259529000]]

\subsection{Application : the exceptional unimodular lattices of rank $\leq 27$ }
\label{subsect:exc2627}

	Let ${\rm E}_n$ denote the set of isometry classes of exceptional unimodular lattices $L$ of rank $n$ with ${\rm r}_1(L) = 0$. As explained after Prop.~\ref{prop:excwithI1}, we have $|{\rm E}_{24}|=24$ (the Niemeier lattices) and ${\rm E}_{25}=\emptyset$. Moreover, by Example~\ref{ex:G26G27} we know:
	
\begin{prop} 
\label{prop:exc2627}
We have $|{\rm E}_{26}|=97$ and $|{\rm E}_{27}|=557$.
%Moreover, for all these $97+557$ lattices $L$ we have $|{\rm Exc}(L)|=2$. 
%as well as  ${\rm O}(L) \simeq {\rm O}(L') where L'=e^\perp
\end{prop}

Most of these exceptional lattices are found as $d$-neighbors of ${\rm I}_n$ using our general algorithms without any specific efforts, 
but for a few of them it is much more efficient to use the isotropic lines described in \S\ref{subsec:vischarvec}, and prescribing suitable visible root systems. We already gave an example of such a situation in dimension $26$ (Example~\ref{ex:exc10A1dim26}). In Table~\ref{tab:exdim27} below, we give a few examples of exceptional unimodular lattices of dimension $27$ with no norm $1$ vectors obtained by this method:

\tabcolsep=3pt
\begin{table}[H]
{\scriptsize 
\renewcommand{\arraystretch}{1.8} \medskip
\begin{center}
\begin{tabular}{c|c|c|c|c}
${\rm R}_2$ & $d$ & $x \in \Z^{27}$ & $\,\,\epsilon\,\,$ & ${\rm reduced \,mass}$ \\

\hline

$6{\bf A}_1$ & $70$ & \scalebox{.65}{({\color{magenta} 1, 1}, 3, 5, 7, 9, 11, 13, 15, 17, {\color{magenta} 19, 19}, 21, {\color{magenta} 23, 23}, {\color{magenta} 25, 25}, {\color{magenta} 27, 27}, 29, 31, {\color{magenta} 33, 33}, 35, {\color{blue} 4, 32, 34})} & 1 &  $1/23040$ \\

$3{\bf A}_1\,{\bf A}_2$ & $74$ & \scalebox{.65}{({\color{magenta}1, 1, 1}, 3, 5, 7, 9, 11, 13, 15, {\color{magenta}17, 17}, {\color{magenta} 19, 19}, 21, 23, 25, 27, {\color{magenta}29, 29}, 31, 33, 35, 37, {\color{blue} 8, 32, 34})} &  0 & $1/483840$\\

$3{\bf A}_1$ & $82$ & \scalebox{.65}{$({\color{magenta} 1, 1}, 3, 5, 7, 9, 11, 13, 15, 17, 19, 21, 23, 25, 27, {\color{magenta} 29, 29}, {\color{magenta} 31, 31}, 33, 35, 37, 39, 41,{\color{blue} 20, 24, 38})$ } & 0 & $1/1512000$ \\

${\bf A}_2$ & $100$ & \scalebox{.65}{$({\color{magenta} 1, 1, 1}, 3, 5, 7, 9, 11, 13, 15, 17, 19, 21, 23, 25, 27, 29, 31, 33, 37, 39, 41, 43, 45, {\color{blue}2, 48, 50})$} & 0  & $1/277136640$ \\

${\bf A}_3$ & $94$ & \scalebox{.65}{$({\color{magenta} 1, 1, 1, 1}, 5, 7, 9, 11, 13, 15, 17, 19, 21, 23, 25, 27, 29, 31, 33, 35, 37, 39, 41, 47, {\color{blue} 4, 44,46})$} & $0$ & $1/489646080$ \\

$\emptyset$ & $96$ & \scalebox{.65}{$(1, 3, 5, 7, 9, 11, 13, 15, 17, 19, 21, 23, 25, 27, 29, 31, 33, 35, 37, 39, 41, 43, 45, 47, {\color{blue} 30, 32,34})$ } & 1 & $1/1268047872$ \\
\end{tabular} 
\end{center}
}
\caption{A few rank $27$ exceptional unimodular lattices of the form ${\rm N}_d(x; \epsilon)$.}
\label{tab:exdim27}
\end{table}

\section{A few more examples in dimension $26$}
\label{sect:exdim26}

In this section, we illustrate our methods by studying a few more examples.
We consider first, in dimension $26$, the root systems in Table~\ref{tab:rsX26} below. 
In the end, they will turn out to be exactly those with at least $7$ isometry classes 
of lattices. 
\tabcolsep=2pt
\begin{table}[H]
{\tiny \renewcommand{\arraystretch}{1.8} \medskip
\begin{center}
\begin{tabular}{c|c|c|c|c|c}
 $R$ &  
$2{\bf A}_1 2{\bf A}_2 2{\bf A}_3 2{\bf A}_4$ & 
$4{\bf A}_1 4{\bf A}_2 2{\bf A}_3$ & 
$3{\bf A}_1 3{\bf A}_2 2{\bf A}_3 {\bf A}_4$ & 
$5{\bf A}_1 3{\bf A}_2 2{\bf A}_3$ & 
$7{\bf A}_1 3{\bf A}_2 {\bf A}_3$ \\ 
$\texttt{reduced mass}$ & $6$ & $77/16$ & $6$ & $4$ & $8/3$ \\
$\sharp\texttt{lattices}$ & $16$ & $15$ & $15$ & $13$ & $12$ \\ \hline
$R$ &  
$6{\bf A}_1 4{\bf A}_2 {\bf A}_3$ & 
$2{\bf A}_1 2{\bf A}_2 3{\bf A}_3 {\bf A}_4$ & 
$6{\bf A}_1 4{\bf A}_2$ & 
$4{\bf A}_1 2{\bf A}_2 4{\bf A}_3$ & 
$3{\bf A}_1 3{\bf A}_2 3{\bf A}_3$ \\ 
$\texttt{reduced mass}$ & $161/48$ & $5$ & $545/576$ & $31/16$ & $15/4$ \\ 
$\sharp\texttt{lattices}$ & $11$ & $11$ & $10$ & $10$ & $10$ \\ \hline
$R$ &  
${\bf A}_1 2{\bf A}_2 2{\bf A}_3 {\bf A}_4 {\bf A}_5$ & 
$8{\bf A}_1 2{\bf A}_2 2{\bf A}_3$ & 
$6{\bf A}_1 2{\bf A}_2 3{\bf A}_3$ & 
$5{\bf A}_1 5{\bf A}_2 {\bf A}_3$ & 
$5{\bf A}_1 3{\bf A}_2 {\bf A}_3 {\bf A}_4$ \\ 
$\texttt{reduced mass}$ & $17/4$ & $85/64$ & $23/12$ & $5/2$ & $13/4$ \\ 
$\sharp\texttt{lattices}$ & $9$ & $9$ & $9$ & $9$ & $9$ \\ \hline
$R$ &  
$4{\bf A}_1 2{\bf A}_2 2{\bf A}_3 {\bf A}_4$ & 
$8{\bf A}_1 2{\bf A}_2$ & 
$4{\bf A}_1 4{\bf A}_2 {\bf A}_3 {\bf A}_4$ & 
$3{\bf A}_1 5{\bf A}_2 2{\bf A}_3$ & 
$3{\bf A}_1 2{\bf A}_2 3{\bf A}_3 {\bf A}_5$ \\ 
$\texttt{reduced mass}$ & $4$ & $491/1344$ & $7/2$ & $9/4$ & $7/4$ \\ 
$\sharp\texttt{lattices}$  & $9$ & $8$ & $8$ & $8$ & $7$ \\ \hline

$R$ &  
$8{\bf A}_1 4{\bf A}_2$ & 
$6{\bf A}_1 2{\bf A}_2 2{\bf A}_3$ & 
$5{\bf A}_1 2{\bf A}_2 2{\bf A}_3 {\bf A}_5$ & 
$2{\bf A}_1 {\bf A}_2 {\bf A}_3 2{\bf A}_4 {\bf A}_5$ & 
$2{\bf A}_1 4{\bf A}_2 2{\bf A}_3 {\bf A}_4$ \\ 
$\texttt{reduced mass}$ & $41/32$ & $15/8$ & $7/4$ & $3$ & $3$ \\  
$\sharp\texttt{lattices}$ & $8$ & $7$ & $7$ & $7$ & $7$ \\ \hline
$R$ &  $2{\bf A}_1 3{\bf A}_2 {\bf A}_3 {\bf A}_4 {\bf A}_5$ &  $10{\bf A}_1$ & ${\bf A}_1 2{\bf A}_2 {\bf A}_3 2{\bf A}_4 {\bf A}_5$ & $9{\bf A}_1 3{\bf A}_2$ & $6{\bf A}_1 2{\bf A}_2 2{\bf A}_3 {\bf A}_4$ \\
$\texttt{reduced mass}$ & $13/4$ & $4424507/116121600$  & $11/4$ & $15/16$ & $7/4$  \\
$\sharp\texttt{lattices}$  & $7$ & $7$ & $6$ & $6$ & $6$
\end{tabular} 
\end{center}
} 
{\small \caption{The root systems $R$ such that there are at least $7$ isometry classes in ${\rm X}_{26}$ with root system $R$ an no norm $1$ vectors (and a few others).}}
\label{tab:rsX26}
\end{table}

As a first example, consider the root system $R:=2{\bf A}_1 \,2{\bf A}_2 \,2{\bf A}_3 \,2{\bf A}_4$.
In this case the reduced mass is $6$, 
hence a priori $R$ is the root system of at least $12$ lattices in ${\rm X}_{26}$. By searching for isotropic lines with visible root system ${\bf A}_1 \,2{\bf A}_2 \,2{\bf A}_3 \,2{\bf A}_4$ (one ${\bf A}_1$ less!), we do quickly find $16$ lattices for $16 \leq d \leq 25$ : $8$ with reduced mass \scalebox{.8}{$1/2$} and $8$ with reduced mass \scalebox{.8}{$1/4$}, the last one being\footnote{That is, 
$x:=(1, 1, 1, 1, 1, 2, 2, 2, 4, 4, 4, 4, 5, 5, 6, 6, 6, 6, 6, 11, 11, 11, 12, 12, 12, 12)$.}
$${\rm N}_{25}(x)\,\,{\rm with}\,\, x=(1^5, 2^3, 4^4, 5^2, 6^5, 11^3, 12^4).$$
\noindent As already said, all these $16$ lattices are distinguished by their invariant $\delta_2$. 
However, they are not distinguished by $\delta_1$. In fact, if we rather choose the invariant $\delta_1$ in the search above, we stop finding new lattices for $d\geq 26$ despite running several thousands of isotropic lines whose associated neighbor has root system $R$. At this point, the remaining mass is \scalebox{.8}{$3/4$}, which represents $(3/4)/6 =12.5 \%$ of the trotal mass, so if we believe in Theorems~\ref{thmi:stat1} \&~\ref{thmi:stat2} (say, ignoring the unknown terms ${\rm emb}(R,-)$) we should have found instead hundreds of times the missing lattices. This is a very strong argument that the chosen invariant is wrong and has to be refined. We have used this strategy several times during our search before we discover fine enough invariants !\ps\ps

Most of the isometry classes in ${\rm X}_{26}$ with root system as in Table~\ref{tab:rsX26} are not especially hard to find. The most complicated is the case $R= 10 {\bf A}_1$, which is why we discussed it at length in the introduction. One lattice with root system $8 {\bf A}_1\,2{\bf A}_2$ has reduced mass \scalebox{.8}{$1/1344$} and is not immediate to find. We immediately find it using a visible isometry of order $7$ and visible root system $7{\bf A}_1$, and more precisely the $29 \cdot 27$-isotropic vector $x \in \Z^{26}$ with 
{\small
$$x \equiv (1^{14},2^7,3,5,9,10,11) \bmod 27$$ 
$$x \equiv(1,t,t^2,t^3,t^4,t^5,t^6,1,t,t^2,t^3,t^4,t^5,t^6,1,t,t^2,t^3,t^4,t^5,t^6,0,0,0,0,0) \bmod 29,$$
}
for $t \equiv 16 \bmod 29$ (an element of order $7$).\ps

We now give another example, namely the case of the root system $R= 22{\bf A}_1\,{\bf D}_4$.
The reduced mass is \scalebox{.8}{$53/60480$}. This is a typical case where we cannot choose anything very close to $R$ as a visible root system. Nevertheless, for even $d$ then $11 {\bf A}_1 {\bf D}_4$ is a possible visible root system a priori. It amounts to choose $d$-isotropic vectors $x \in \Z^{26}$ having $4$ coordinates equal to $d/2$, as well as $11$ other pairs of equal coordinates. This forces $d\geq 24$. And indeed, we quickly find the two following lattices, with respectives masses \scalebox{.8}{$1/1152$} and \scalebox{.8}{$1/120960$}, and conclude :
$${\rm N}_{30}(1^2, 2^2, 3^2, 4^2, 7^2, 8^2, 9^2, 10^2, 11^2, 13^2, 14^2, 15^4)$$
$${\rm N}_{46}(1^2, 5^2, 7^2, 9^2, 11^2, 13^2, 15^2, 17^2, 19^2, 20^2, 21^2, 23^4).$$
Both lattices can also be understood with the help of ${\rm X}_{22}$. For instance,
the orthogonal of the ${\rm D}_4$ in the latter is actually the even sublattice of the unique lattice in ${\rm X}_{22}$
with root system $22{\bf A}_1$ (and whose even part is also the orthogonal of some $2{\bf A}_1$ inside the Niemeier lattice with root system $24 {\bf A}_1$ !).

\section{A few constructions of lattices in neighbor form}
\label{sect:lattconstr}

First, it is straightforward to ``add ${\rm I}_m$'' to unimodular lattices in neighbor forms. \ps

\begin{lemma} 
\label{lem:ajoutin0}
Let $x \in \Z^n$ be $d$-isotropic and define $y \in \Z^{n+m}$ by $y_i=x_i$ for $i\leq n$, 
$y_i=0$ otherwise. Then $y$ is $d$-isotropic, and we have ${\rm N}_d(y)={\rm N}_d(x) \perp {\rm I}_m$ if $d$ is odd, 
and ${\rm N}_d(y; \epsilon) = {\rm N}_d(x; \epsilon) \perp {\rm I}_m$ if $d$ is even and $\epsilon \in \{0,1\}$.
\end{lemma}

\begin{pf} We have ${\rm M}_d(y) = {\rm M}_d(x) \perp {\rm I}_m$. Consider a $d$-neighbor ${\rm N}_d(x')$ associated with $x$, and define $y' \in \Z^{n+m}$ by $y'_i=x'_i$ for $i \leq n$, $y'_i=0$ otherwise.  As $x$ (resp. $x'$) coincides with 
$y$ (resp. $y'$) inside $\Z^{n+m}$, we have ${\rm N}_d(y')={\rm N}_d(x') \perp {\rm I}_m$. 
\end{pf}

\begin{lemma} 
\label{lem:voisdevois}
Let $L,N \in \mathcal{L}_n$.
Assume that $N$ is a $d$-neighbor of ${\rm I}_n$ associated to the $d$-isotropic vector $x \in \Z^n$,
and that $L$ is a $d'$-neighbor of $N$, with ${\rm gcd}(d,d')=1$. Then $L$ is a $dd'$-neighbor of ${\rm I}_n$
associated to a $dd'$-isotropic vector $y \in \Z^n$ satisfying $y_i \equiv x_i \bmod d$ for all $i=1,\dots, n$.
\end{lemma}

\begin{pf} 
Using ${\rm gcd}(d,d')=1$, and localizing at primes dividing $dd'$, one readily observes $L \cap {\rm I}_n \subset N$, 
as well as ${\rm I}_n / (L \cap {\rm I}_n)\, \simeq\, \Z/dd'$. Let $\ell \subset {\rm I}_n \otimes \Z/dd'$ be the $dd'$-isotropic line satisfying $L \cap {\rm I}_n = {\rm M}_{dd'}(\ell)$. By the inclusion 
${\rm M}_{dd'}(\ell) \subset N \cap {\rm I}_n = {\rm M}_{d}(x)$, the reduction mod $d$ of $\ell$ has to be $l(x)$, and the statement follows.
\end{pf}

Our aim now is to address the following more interesting :\ps\ps\ps

\noindent {\bf Problem :} {\it {\rm (Addition of ${\bf D}_m$)}  Assume we know all (isometry classes of) unimodular lattices 
$L \in \mathcal{L}_n$ with given root system $R$, and in neighbor form. Choose $m\geq 2$.  
Find neighbor forms for the unimodular lattices $U \in \mathcal{L}_{n+m}$ with root system $R \coprod {\bf D}_m$.}\ps
\ps\ps

\begin{lemma} 
\label{lem:existL2nei}
Assume $U \in \mathcal{L}_{n+m}$ contains $D:=\{0\} \times {\rm D}_m$ as a saturated subgroup {\rm (}$m\geq 2${\rm )}, and set $L_0=D^\perp \cap U$. 
Then there is a unique\footnote{Actually, we have ${\rm res}\, L_0 \simeq - {\rm res}\, {\rm D}_m$
so there is a unique unimodular lattice $L$ containing $L_0$ of index $2$, unless we have $m \equiv 0 \bmod 4$ in which case there are three.} $L \in \mathcal{L}_n$ containing $L_0$ with index $2$ and such that 
$U$ is a $2$-neighbor of $L \perp {\rm I}_m$. 
\end{lemma}

\begin{pf} By~\S\ref{notations} (iii), the lattice $U$ is the inverse image of some (order $4$) Lagrangian $I \subset {\rm res}\, L_0 \perp {\rm res}\, {\rm D}_m$ that is transversal to both summands. In particular, $I$ contains a unique element 
of the form $a+b$, where 
$b$ is a generator of ${\rm I}_m/{\rm D}_m \simeq \Z/2$ and $a\in {\rm res}\, L_0$, necessarily satisfying $a.a \equiv 0 \bmod \Z$ and $2a = 0$. The inverse image of $\Z/2 \,a$ in $L_0^\sharp$ is thus a unimodular lattice containing $L_0$ that we denote by $L$. By construction, $U \cap (L \perp {\rm I}_m)$ contains $L_0 \cap {\rm D}_m$ with index $2$, hence has index $2$ in $L \perp {\rm I}_m$.
\end{pf}

 \begin{comment}
For this purpose, we fix some $U \in \mathcal{L}_{n+m}$ containing some saturated subgroup $D$ isometric to ${\rm D}_m$ with $m\geq 2$. Up to replacing $U$ with some isometric lattice, we may and do assume $D = \{0\} \times {\rm D}_m \subset \R^n \times \R^m$. Set $L_0 = D^\perp \cap U$. We have then ${\rm res}\, L_0 \simeq - {\rm res}\, {\rm D}_m$, so by~\S\ref{notations} (iii) and Table~\ref{tab:venkov} there is a unique unimodular lattice $L \subset L_0^\sharp$ containing $L_0$ with index $2$, unless we have $m \equiv 0 \bmod 4$, in which case there are exactly three such $L$.
\end{comment}

\begin{prop} 
\label{prop:adddn} 
Assume $U \in \mathcal{L}_{n+m}$ contains $D:=\{0\} \times {\rm D}_m$ as a saturated subgroup {\rm (}$m\geq 2${\rm )}. Let $L$ be the associated rank $n$ unimodular lattice given by Lemma~\ref{lem:existL2nei}. Assume we have $L \simeq {\rm N}_d(x)$ for some odd integer $d$ and a $d$-isotropic 
$x \in \Z^n$. Then we have $U \simeq {\rm N}_{2d}(y)$ for some $2d$-isotropic $y \in \Z^{n+m}$ satisfying 
$y_i \equiv x_i \bmod d$ for $i=1,\dots,n$, and $y_i \equiv d \bmod 2d$ for $i=n+1,\dots,n+m$.
\end{prop}

\begin{pf} 
By Lemma~\ref{lem:existL2nei},  $U$ is isometric to some $2$-neighbor $U'$ of ${\rm N}_d(x) \perp {\rm I}_m$ 
in which the natural ${\rm D}_m$ is saturated. 
By Lemma~\ref{lem:ajoutin0}, we have ${\rm N}_d(x) \perp {\rm I}_m={\rm N}_d(\xi)$ with $\xi_i=x_i$ for $i=1,\dots,n$,
and $\xi_i=0$ for $i=n+1,\dots,n+m$. By Lemma~\ref{lem:voisdevois}, we have thus $U' ={\rm N}_{2d}(y)$ with 
$y_i \equiv \xi_i$ for $i=1,\dots,n+m$. As $U'$ contains the natural ${\rm D}_m$, we also have, for all $i,j>n$,
 $y_i \equiv  y_j \bmod 2d$ and $y_i \equiv 0 \bmod d$. If we have $y_i \equiv 0 \bmod 2d$ for $i>n$, then $U'$ contains the natural ${\rm I}_m$, contradicting the saturation of ${\rm D}_m$.
\end{pf}

\begin{remark} 
\label{rem:strict2nei}\label{rem:addn}
If $x$ as above is given, there are at most $2^n$ choices for $y \bmod 2d$, which may be a lot.
In practice, it is often useful to restrict the search for $y$ with visible root system $R^{\rm v} \coprod {\rm D}_m$,
where $R^{\rm v}$ denotes the visible root system of ${\rm N}_d(x)$. A concrete example of such a search was given in the introduction : this is how the last two lattices of Table~\ref{tab:10A1dim26} have been found {\rm (}case $m=2$, we have 
${\bf D}_2 \simeq 2\,{\bf A}_1$.{\rm )} We have used this method in very many cases during our proof of Theorem~\ref{X2627}.
\end{remark}

We end with a proposition the case $n \equiv 0 \bmod 4$ (see Proposition.~\ref{prop:defcompa}).

\begin{prop}
\label{prop:neicompanions} 
Assume $n \equiv 4 \bmod 8$ and $x \in \Z^n$ is $d$-isotropic with $d$ odd.
Assume $y \in \Z^n$ satisfies $y_i \equiv x_i \bmod d$ and $y_i\equiv 1 \bmod 2$ for all $i=1,\dots,n$.
Then the two companions of ${\rm N}_d(x)$ are the ${\rm N}_{2d}(y)^{\pm}$. 
\end{prop}

\begin{pf} 
Each companion $L$ of ${\rm N}_d(x)$ is a $2d$-neighbors of ${\rm I}_n$ by Lemma~\ref{lem:voisdevois},
so it is enough to show $L \supset {\rm M}_{2d}(y)$. But the even part of ${\rm I}_n$ is ${\rm M}_2(1^n)$, 
so that of ${\rm N}_d(x)$ (hence of $L$) is ${\rm N}_d(x) \cap {\rm M}_2(1^n)$ since $d$ is odd.
We deduce $L \supset {\rm N}_d(x) \cap {\rm M}_2(1^n) \supset  {\rm M}_d(x) \cap {\rm M}_2(1^n)={\rm M}_{2d}(y)$.
\end{pf}

\section{The Jordan-H\"older factors of ${\rm O}(L)$ for $L \in \mathcal{L}_n$ and $n\leq 27$}
\label{sect:newlatt}	
${}^{}$\indent In this section, we discuss the Jordan-H\"older factors of ${\rm O}(L)$ for $L \in \mathcal{L}_n$
and $n \leq 27$. It is enough to study the (usually smaller) reduced isometry groups ${\rm O}(L)^{\rm red}$ (see~\S\ref{def:redisogp}). For this, we view ${\rm O}(L)^{\rm red}$ as a finite permutation groups of suitable small norm vectors of $L$, of which the Plesken-Souvignier algorithm provides generators (see Remark~\ref{pleskensouvignier}), and then use the permutation groups algorithms in \cite{GAP}. For $n \leq 25$, the output is that the non-abelian simple groups appearing as a Jordan-H\"older factor of some ${\rm O}(L)^{\rm red}$ are
	$${\rm Co}_1, \,\,{\rm Co}_2, \,\,{\rm Co}_3, \,\,{\rm HS}, \,\,{\rm M}_{24}, \,\,{\rm M}_{22}, \,\,{\rm M}_{12}, \,\,{\rm U}_6(2), \,\,{\rm A}_8,\,\, {\rm A}_7, \,\,{\rm A}_6, \,\,{\rm A}_5\, \, {\rm and} \,\,{\rm L}_2(7).$$
\indent In dimension $26$, there are only $39$ lattices $L$ in ${\rm X}_{26}$ with ${\rm r}_1(L)=0$ and such that 
$|{\rm O}(L)^{\rm red}|$ is both $\equiv 0 \bmod 4$ and not of the form $p^a q^b$ with $p,q$ primes. We checked that in each of these cases ${\rm O}(L)^{\rm red}$ is indeed non solvable\footnote{With a unique non abelian simple factor, except for the lattice with root system $10\,{\bf A}_1$ and reduced mass $\scalebox{.8}{1/7372800}$, which contains two factors isomorphic to ${\rm A}_5$.}. For exactly $6$ lattices of these $39$ lattices, we obtain the following ``new'' Jordan-H\"older factors, not appearing for unimodular lattices of rank $<26$ : ${\rm A}_9, \,\,{\rm O}_5(5),\,\, {\rm L}_3(4), \,\,{\rm L}_3(3)\,\,\,{\rm and}\, \,  \,\,{\rm L}_2(8)\,\,\,({\rm twice})$.
\tabcolsep=3.5pt
\begin{table}[H]
{\scriptsize 
\renewcommand{\arraystretch}{1.8} \medskip
\begin{center}
\begin{tabular}{c||c|c|c|c|c}
${\rm R}_2(L)$ & $10{\bf A}_1$ & $\emptyset$ & $22{\bf A}_1\,{\bf D}_4$ &  $26{\bf A}_1$ & ${\bf A}_1\,9{\bf A}_2\, \, \&  \,\,9{\bf A}_1\,{\bf A}_2$ \\ \hline
$|{\rm O}(L)^{\rm red}|$ & $92897280$ & $18720000$ & $120960$ & $11232$ & $3024$ \\  \hline
${\rm O}(L)^{\rm red} $ &  \scalebox{.8}{${\rm C}_2.({\rm C}_2^8.{\rm A}_9)$} &  \scalebox{.8}{${\rm O}_5(5)\,:\,{\rm C}_4$} &  \scalebox{.8}{$({\rm L}_3(4) : {\rm C}_3) : {\rm C}_2$} &  \scalebox{.8}{${\rm L}_3(3) : {\rm C}_2$} &  \scalebox{.8}{${\rm C}_2 \times ({\rm L}_2(8) \,:\, {\rm C}_3)$}
\end{tabular} 
\end{center}
}
\caption{{\small The $6$ lattices in ${\rm X}_{26}$ whose reduced isometry group has a ``new'' Jordan-H\"older factor, as described by \texttt{GAP}'s \texttt{StructureDescription} function.}}
\label{tab:X26newsimpl}
\end{table}

	In a similar manner,  there are only $74$ lattices in ${\rm X}_{27}$ with no norm $1$ vectors and whose reduced isometry group has order both $\equiv 0 \bmod 4$ and not of the form $p^a q^b$ with $p,q$ primes. Their reduced isometry groups are indeed non solvable in all cases, with a unique non abelian factor. Again, for exactly $6$ lattices, a new Jordan-H\"older factor appears, namely : ${\rm U}_6(2)$, ${}^3{\rm D}_4(2)$, ${\rm U}_4(2)$, ${\rm O}_5(3)$, ${\rm U}_3(5)$, ${\rm M}_{11}$ and ${\rm PSL}_2(11)$. 
	
\tabcolsep=3pt
\begin{table}[H]
{\scriptsize 
\renewcommand{\arraystretch}{1.8} \medskip
\begin{center}
\begin{tabular}{c|c|c|c|c|c}
 ${\bf D}_4$ & $\emptyset$ & $\emptyset$ &  $3{\bf A}_1$ & ${\bf A}_1\,{\bf A}_3$ &  $11{\bf A}_1\,{\bf A}_4$ \\ \hline
 $55180984320$ & $634023936$ & $1658880$ & $756000$ & $7920$ & $660$ \\  \hline
 \scalebox{.8}{${\rm C}_2 \times (({\rm U}_6(2) : {\rm C}_3) : {\rm C}_2)$} & \scalebox{.8}{${\rm C}_2 \times ({}^{3}{\rm D}_4(2) : {\rm C}_3)$} & \scalebox{.8}{${\rm C}_2 \times ({\rm C}_2^6.{\rm O}_5(3))$} & \scalebox{.8}{${\rm C}_2 \times (({\rm U}_3(5) : {\rm C}_3) : {\rm C}_2)$} & \scalebox{.8}{${\rm C}_2 \times {\rm M}_{11}$} & \scalebox{.8}{${\rm C}_2 \times {\rm L}_2(11)$}
\end{tabular} 
\end{center}
}
\caption{{\small The $6$ lattices $L$ in ${\rm X}_{27}$ whose reduced isometry group has a ``new'' Jordan-H\"older factor (same format as Table~\ref{tab:X26newsimpl}).}}
\label{tab:X27newsimpl}
\end{table}

\section{Proof of Theorem~\ref{thmi:stat1}.}
\label{sect:pfthmb}

We now explain how to modify the proof of Thm. A in \cite{chestat} in order to deal as well with the case of 
cyclic $d$-neighbors with $d$ non necessarily prime. Note that the result is obvious for $n\leq 4$, since we know $|{\rm X}_n|=1$ in this case. To keep the discussion short we freely use the notations in \cite{chestat}. 

\begin{thm}\label{thmstatnotprime} Let $L$ be an integral lattice in $\R^n$ with $n>4$. 
Assume $\mathcal{G}={\rm Gen}(L)$ is a single spinor genus, choose $L'$ in $\mathcal{G}$, and for $d$
prime to $2\det L$ denote by ${\rm N}_d(L,L')$ the number of cyclic $d$-neighbors of $L$ which are isometric to $L'$.
Fix $\epsilon>0$.
Then for $d$ prime to $2 \det L$, we have 
\begin{equation} \label{ramprop} \,\,\,\frac{ {\rm N}_d(L,L') }{| {\rm C}_L(\Z/d) |}  =  \frac{1/|{\rm O}(L')|}{{\rm m}(\mathcal{G})} + {\rm O}( \frac{1}{d^{1-\epsilon}})\, \, \, {\rm when}\, \, d \rightarrow \infty.\end{equation} 
\end{thm}

For this, we first generalize the discussion in \cite[\S 4]{chestat} to $p^\alpha$-neighbors.
Fix $p$ an odd prime and $V$ a nondegenerate quadratic space over $\Q_p$. 
We assume $\dim V\geq 3$ and that the set $\mathcal{U}(V)$ of {\it unimodular} integral $\Z_p$-lattices in $V$ is non empty, fix $L \in \mathcal{U}(V)$. 
For any integer $\alpha \geq 1$, the cyclic $p^\alpha$-neighbors of $L$ are the lattices $N \in \mathcal{U}(V)$ such that $L/(L\cap N)$ is cyclic
of order $p^\alpha$. They form a subset $\mathcal{N}_{p^\alpha}(L)$ of $\mathcal{U}(V)$. A similar argument 
as in Corollary~\ref{cor:bijlinemap} shows that $\mathcal{N}_{p^\alpha}(L)$ is in natural bijection with ${\rm C}_L(\Z/p^\alpha)$,
and form a single ${\rm O}(L)$-orbit. Assuming $L \,=\, (\Z_p e \,\oplus \,\Z_p f) \,\perp \,M$ with $M$ unimodular, $e.e=f.f=0$ and $e.f=1$ (this is always possible), then 
%we have for instance
\begin{equation}
\label{eq:onepalphanei}
(\Z_p \,p^\alpha\,e\, \oplus \,\Z_p \,p^{-\alpha}\,f) \,\perp\, M\, \, \, \in \, \, \mathcal{N}_{p^\alpha}(L).
\end{equation}

We denote by ${\rm T}_{p^\alpha}$ the element of the Hecke ring ${\rm H}_V$ defined, for $L \in \mathcal{U}(V)$,
by ${\rm T}_{p^\alpha}L = \sum_{N \in \mathcal{N}_{p^\alpha}(L)} N$. The key estimate is the following.

\begin{prop} 
\label{prop:majoh}
Let $V$ be a non-degenerate quadratic space over $\Q_p$ of dimension $\geq 5$, $L \in \mathcal{U}(V)$, and $\alpha \geq 1$ an integer. Let $U$ a unitary irreducible unramified $\C[{\rm O}(V)]$-module, and $\lambda \in \C$ the eigenvalue of ${\rm T}_{p^\alpha}$ on the line $U^{{\rm O}(L_p)}$. If $\dim U>1$ then we have $|\lambda| \leq \, |{\rm C}_{L}(\Z/p^\alpha)| \, (\alpha+1)^2\,p^{-\alpha}$.
\end{prop}

\begin{pf} The proof is the same as that of \cite[Prop. 6.5]{chestat}, up to replacing the subset $C \subset {\rm O}(V)$ 
{\it loc. cit.} with the double coset ${\rm O}(L)c{\rm O}(L)$ of elements $g \in {\rm O}(V)$ such that $g(L)$ is a cyclic $p^\alpha$-neighbor of $L$. In the notation of that proof, and by Formula~\eqref{eq:onepalphanei}, we may take $c= \varepsilon_1^\ast(p^\alpha)$. 
Applying \cite[Thm. 1.2]{oh}, we obtain the inequality
$|\langle c e, e \rangle| \leq \Xi_{{\rm PGL}_2(\Q_p)}(p^\alpha)^2$
with {\small $\Xi_{{\rm PGL}_2(\Q_p)}(p^\alpha) = \frac{1}{p^{\alpha/2}}\frac{\alpha(p-1)+p+1}{p+1}$}.
\end{pf}

We now follow the arguments given in the beginning of \cite[\S 6]{chestat}.
For the application to Theorem~\ref{thmstatnotprime} (``lattice case'') the finite set $S$ {\it loc. cit.} is the
set of primes dividing $D:=2 \det L$. For each integer $d$ prime to $D$ we have a natural 
global Hecke operator ${\rm T}_d$ corresponding to the cyclic $d$-neighbors, and generalizing those above when $d$ is a prime power. For $d=d'd''$ with coprime $d'$ and $d''$, it satisfies ${\rm T}_d={\rm T}_{d'}{\rm T}_{d''}$, and its degree ${\rm c}_V(d):=|{\rm C}_L(\Z/d)|$ of course satisfies ${\rm c}_V(d)={\rm c}_V(d'){\rm c}_V(d'')$ as well.
These operators pairwise commute and act in a diagonalisable way on the space of automorphic forms denoted ${\rm M}(K)$ {\it loc. cit}. 
Fix a common eigenvector $v$ in the subspace ${\rm M}(K)$ and denote by $\lambda(d)$ the eigenvalue
of ${\rm T}_d / {\rm c}_V(d)$ on $v$. If $v$ is in the subspace ${\rm M}(K)^0$ {\it loc. cit.}\,, and $p$ is a prime not dividing $D$, 
Prop.~\ref{prop:majoh} implies that for all $\alpha\geq 1$ we have
$\lambda(p^\alpha) \,\leq \,\, (\alpha+1)^2 p^{-\alpha}$.
It follows that for each $d$ prime to $D$ we have $|\lambda(d)| \leq \sigma_0(d)^2 d^{-1}$, 
where $\sigma_0(n)$ denotes the number of divisors of the integer $n\geq 1$. 
It is well-known that, for any $\epsilon>0$, we have $\sigma_0(n) = {\rm O}(n^\epsilon)$ for $n \rightarrow \infty$.
This proves $\lambda(d) = {\rm O}(d^{2\epsilon-1})$ for $d \rightarrow \infty$, and the result follows. $\square$

\begin{remark}
\begin{itemize}
\item[(i)] If $\det L$ is odd, $L$ is even, and if we restrict to {\it even} cyclic $d$-neighbors, then the statement (and its proof) also holds with $2\,\det L$ replaced by $\det L$.\ps
\item[(ii)] The cases $n=3,4$ could be handled as well using similar methods as in the second proof of Theorem 6.3 in {\rm \cite{chestat}} (see the end of \S 6).
\end{itemize}
\end{remark}

\end{document}